\documentclass[11pt,reqno]{amsart}
\usepackage{float}
\usepackage[pdftex,hyperfootnotes=false,colorlinks=false,hypertexnames=false]{hyperref}

\usepackage[english]{babel}
\usepackage{amsmath,amssymb,anyfontsize,enumerate,stmaryrd,mathtools,thmtools,tikz,%txfonts,
nccmath}
\usepackage[oldstyle]{libertine}%
\usepackage[T1]{fontenc}
\usepackage[final]{microtype}
\usepackage[utf8]{inputenc}
\usepackage{csquotes}
\usepackage{comment}
\makeatletter
\renewcommand*\libertine@figurestyle{LF}
\makeatother
\makeatletter
\renewcommand*\libertine@figurestyle{OsF}
\makeatother
\usepackage{mathabx}
\usepackage{tikz-cd}
\tikzcdset{scale cd/.style={every label/.append style={scale=#1},
    cells={nodes={scale=#1}}}}
\usepackage[noabbrev,capitalise]{cleveref}
\usepackage{autonum}
\usetikzlibrary{%
  matrix,%
  calc,%
  arrows%
}

\usepackage{tabularx}
\usepackage{booktabs}
\usepackage[labelfont=bf,format=plain,justification=raggedright,singlelinecheck=false]{caption}

\usepackage[sanitizesort]{glossaries}
\newglossary[slg]{symbolslist}{syi}{syg}{Index of Notation}
\makenoidxglossaries

\newglossaryentry{TT}{
 name=\ensuremath{\mathbb{T}\mathcal{T}},
 description={The category of real tori with integral structure},
 type=symbolslist
}
\newglossaryentry{TA}{
 name=\ensuremath{\mathbb{T}\mathcal{A}},
 description={The category of tropical abelian varieties},
 type=symbolslist
}
\newglossaryentry{TC}{
 name=\ensuremath{\mathbb{T}\mathcal{C}},
 description={The category of tropical curves},
 type=symbolslist
}
\newglossaryentry{Ab}{
 name=\ensuremath{Ab},
 description={The category of abelian groups},
 type=symbolslist
}

\newglossaryentry{Sigma}{
 name=\ensuremath{\Sigma},
 description={Depending on the context: A real torus with integral structure or pptav},
 type=symbolslist
}
\newglossaryentry{widecheckSigma}{
 name=\ensuremath{\widecheck{\Sigma}},
 description={The dual of Sigma},
 type=symbolslist
}

\newglossaryentry{zeta}{
 name=\ensuremath{\zeta},
 description={Polarization on $\Sigma$},
 type=symbolslist
}
\newglossaryentry{fzeta}{
 name=\ensuremath{f_\zeta},
 description={Isogeny induced by \ensuremath{\zeta}},
 type=symbolslist
}

\newglossaryentry{f}{
 name=\ensuremath{f:\Sigma_1 \rightarrow \Sigma_2},
 description={A morphism in \ensuremath{\mathbb{T}\mathcal{T}}},
 type=symbolslist
}
\newglossaryentry{widecheckf}{
 name=\ensuremath{\widecheck{f}:\widecheck{\Sigma}_2 \rightarrow \widecheck{\Sigma}_1},
 description={The dual of \ensuremath{f}},
 type=symbolslist
}

\newglossaryentry{Ker}{
 name=\ensuremath{Ker(f)_0},
 description={Kernel of \ensuremath{f} in \ensuremath{\mathbb{T}\mathcal{T}} or \ensuremath{\mathbb{T}\mathcal{A}}},
 type=symbolslist
}

\newglossaryentry{ker}{
 name=\ensuremath{\ker(f)},
 description={Kernel of \ensuremath{f} in  \ensuremath{Ab}},
 type=symbolslist
}

\newglossaryentry{zetaGamma}{
 name=\ensuremath{\zeta_\Gamma},
 description={Natural principal polarization on Jac\ensuremath{(\Gamma)}},
 type=symbolslist
}

\newglossaryentry{PhiP0}{
 name=\ensuremath{ \Phi_{P_0}},
 description={The tropical Abel-Jacobi map with reference point \ensuremath{P_0}},
 type=symbolslist
}
\newglossaryentry{TE}{
 name=\ensuremath{\mathbb{T}E},
 description={Tropical elliptic curve},
 type=symbolslist
}

\theoremstyle{plain}
    \newtheorem{theorem}{Theorem}
    \newtheorem{construction/theorem}[theorem]{Construction/Theorem}
    \newtheorem{corollary}[theorem]{Corollary}
    \newtheorem{lemma}[theorem]{Lemma}
    \newtheorem{proposition}[theorem]{Proposition}
    
\theoremstyle{definition}
    \newtheorem{remark}[theorem]{Remark}
    \newtheorem{remark/reference}[theorem]{Remark/Reference}
    
    \newtheorem{example}[theorem]{Example}
    \newtheorem{definition}[theorem]{Definition}
    \newtheorem{construction}[theorem]{Construction}
    \newtheorem{algorithm}[theorem]{Algorithm}
     
     \newtheorem{convention}[theorem]{Convention}

\newcounter{diagram}
\crefname{diagram}{Diagram}{Diagrams}
\newcounter{savedequation}
 \newenvironment{diagram}[1][]
 {%
  \setcounter{savedequation}{\value{equation}}%
  \setcounter{equation}{\value{diagram}}%
  \crefalias{equation}{diagram}%
  \begin{equation}\begin{tikzcd}[#1]%
 }
 {%
  \end{tikzcd}\end{equation}%
  \setcounter{diagram}{\value{equation}}%
  \setcounter{equation}{\value{savedequation}}%
  \ignorespacesafterend % <-- no space after the diagram environment
 }

\newcommand\sbullet[1][.5]{\mathbin{\vcenter{\hbox{\scalebox{#1}{$\bullet$}}}}}

\DeclareMathOperator{\Coker}{Coker}
\DeclareMathOperator{\Hom}{Hom}
\DeclareMathOperator{\Prin}{Prin}
\DeclareMathOperator{\Ker}{Ker}
\DeclareMathOperator{\im}{im}
\DeclareMathOperator{\Div}{Div}
\DeclareMathOperator{\rk}{rk}
\DeclareMathOperator{\Jac}{Jac}
\DeclareMathOperator{\Pic}{Pic}
\DeclareMathOperator{\H1}{H_1}
\DeclareMathOperator{\Co}{C_0}
\DeclareMathOperator{\C1}{C_1}
\title{Tropical Split Jacobians of genus 2 and optimal covers}
\author[L.~Cobigo]{Lou-Jean Leila Cobigo}
\address{L.~Cobigo: Eberhard Karls Universität Tübingen, Fachbereich Mathematik, Auf der Morgenstelle 10, 72076 Tübingen, Germany}
\email{lou-jean-leila.cobigo-bihavan@uni-tuebingen.de}

\begin{document}
\begin{abstract}
We explore connections between the category of \emph{tropical abelian varieties (tav)}, $\mathbb{T}\mathcal{A}$, and the the category of \emph{tropical curves}, $\mathbb{T}\mathcal{C}$, first in a broader context and then specifically by studying the phenomenon of \emph{tropical split Jacobians}. Jacobians of genus $2$ curves are two-dimensional tav and as such more complicated than their one-dimensional cousins. Whenever $\Gamma$, however, is a covering of an elliptic curve, it so happens that $\Jac(\Gamma)$ splits into simpler objects, the direct sum of elliptic curves. This relation is pathological in essentially two ways, the splitting of $\Jac(\Gamma)$ is not unique, and it is a priori not clear how to compute it. 
Similar to algebraic geometry, optimal coverings offer a remedy for both: They resolve indeterminacy as they provide us with a \emph{canonical} choice. They resolve indeterminability as they provide us with an algorithmic approach. Our methods build on theory developed by Len, Mikhalkin, Röhrle, Ulirsch, Zakharov, Zharkov and many more. In accordance with this heritage, we want to put forth tropical geometry as a setting in which we can see abstraction at work. This means pairing "abstract machinery" with a constructive/algorithmic approach.
\end{abstract}
\maketitle
\tableofcontents
\glsaddall
\printnoidxglossary[type=symbolslist,sort=def]
\section{Introduction}
We explore the following connection between the category of \emph{tropical abelian varieties (tav)}, $\mathbb{T}\mathcal{A}$, and the the category of \emph{tropical curves}, $\mathbb{T}\mathcal{C}$:
Let $\Gamma$, $\tilde{\Gamma}$ be tropical curves, $\varphi: \Gamma \rightarrow \tilde{\Gamma}$ a cover. Denote by $\Jac(\Gamma)$ and $\Jac(\tilde{\Gamma})$ the respective Jacobians and consider the assignment $\mathcal{F}$ sending $\Gamma \in \mathbb{T}\mathcal{C}$ to $\Jac(\Gamma)\in \mathbb{T}\mathcal{A}$, and sending $\varphi$ to the push-forward $\varphi_*: \Jac(\Gamma) \rightarrow \Jac(\tilde{\Gamma})$.

As $\mathcal{F}$ provides a way to communicate between $\mathbb{T}\mathcal{C}$ and $\mathbb{T}\mathcal{A}$, it becomes essential to appreciate the nature of this communication:
Given a phenomenon $P$ in $\mathbb{T}\mathcal{C}$ (an object that satisfies a certain property, a collection of interrelated objects,...), one may wonder at \emph{how} it is expressed as phenomenon in $\mathbb{T}\mathcal{A}$, or vice versa. 
The present works embeds in this framework by studying the phenomenon of \emph{tropical split Jacobians}.
\paragraph{\emph{Tropical split Jacobians in }  \ensuremath{\mathbb{T}\mathcal{A}} \emph{and in} \ensuremath{\mathbb{T}\mathcal{C}}.}
Jacobians of genus $2$ curves are two-dimensional tav and as such more complicated than their one-dimensional cousins. Whenever $\Gamma$, however, is a covering of an elliptic curve, it so happens that $\Jac(\Gamma)$ splits into simpler objects, the direct sum of elliptic curves (Theorem \ref{theorem_introduction1}). This decrease in complexity comes with a multitude of advantages, but suffers from two major problems
\begin{itemize}
    \item Indeterminacy: Non-uniqueness of the splitting.
    \item Indeterminability: Non-computability of the splitting.
\end{itemize}
Similar to algebraic geometry, optimal coverings (see Definition \ref{definition_optimalmap}) offer a remedy for both: They resolve indeterminacy as they provide us with a \emph{canonical} choice. They resolve indeterminability as they provide us with an algorithmic approach.
\subsection{Context and Background}
Attempts at unification of various constructions in the context of tropical curves and tavs (e.g. \cite{MR2457739}, \cite{MR3375652}) , using category theory as common language, have recently been made in \cite{röhrle2024tropicalngonalconstruction} continuing the work of \cite{MR4261102}. We follow up on this development and  extend the theory in tropical geometry, see Sections \ref{section_catofrealtori}, \ref{section_catoftav} and \ref{section_crossingbridges}, to broaden its scope of action and application, and to provide a transparent and structured framework to its organic growth.
There is a large body of work on tavs and tropical covers. Jacobians and Prym varieties are two classes of tavs that have been studied  extensively, bringing forth connections between curves and tavs; connections that are fruitful at both levels, the level of objects and the level of moduli spaces (see e.g. \cite{MR2457739}, \cite{MR2772537}, \cite{MR2739784}, \cite{MR2968636}, \cite{MR4261102}, \cite{MR4382460}, \cite{röhrle2024tropicalngonalconstruction}).
Prym varieties provide a link to unramified covers (see \cite{MR3782424}). Ramified tropical covers, on the other hand, have made a notable appearance in the context of various enumerative problems (see e.g. \cite{MR2661417},\cite{MR3717092}), as they open the way to a combinatorial approach.  

The study of connections between ramified covers and tavs, initiated by \cite{röhrle2024tropicalngonalconstruction}, is new. Our paper provides a new approach to study the interplay of ramified covers and tavs, here in the context of genus 2 curves covering elliptic curves and their Jacobians.
\subsubsection{Inspiration from Algebraic Geometry}
Once the tropical framework is established, we realize that classical questions like:
\begin{itemize}
    \item When is the Jacobian of a curve $C$ of genus $2$ isogenous to a product of two elliptic curves, $E$ and $E'$?
    \item Is such a decomposition unique and can it be made explicit?  
\end{itemize}
that go back to Jacobi, Legendre and others are just as natural in the tropical world. Moreover, the kinship of both areas provides guidance for the tropical investigation; many constructions behave analogously. This lets us oftentimes import ideas from algebraic geometry. We give a short overview to make the parallels visible, for details we refer to (\cite{MR0936803}, \cite{MR1085258}):
The answer to the first question is positive if and only if there exists a finite cover $f: C \rightarrow E $. The answers to the remaining two questions are "in general it is not unique" and "making it explicit is not trivial". 
Whenever $f$ is \emph{optimal}, however, i.e. $f$ does not factor through a non-trivial isogeny, the situation changes: There exists another elliptic curve $E'$, an optimal covering $f': C \rightarrow E'$ and an isogeny $\phi: \Jac(C) \rightarrow E \oplus E'$ whose kernel is isomorphic to the group of $deg(f)$-torsion points of $E$. Explicit descriptions of $f'$ and $E'$ can already be found in the works of Jacobi and Legendre for degree $\leqslant 4$, Kuhn, Frey and Kani with a more modern approach obtain results for various cases of degree $\leqslant 11$ (\cite{MR0936803},\cite{MR1085258}).

We want to put forth tropical geometry as a setting that allows abstraction and concreteness to coexist: The classical treatment uses a powerful, rather abstract machinery that was developed over a long period of time. What makes a tropical investigation so interesting is that with tropical geometry we have a setting that makes both possible: We will work at the same level of abstraction, at the same time, make abstract concepts accessible by providing concrete algorithms to work with them.
\subsection{Results}
Our first result identifies a phenomenon in $\mathbb{T}\mathcal{C}$ that relates to split Jacobians as follows.
\begin{theorem}\label{theorem_introduction1}(Theorem \ref{theorem_Jacsplitsiffcoverexists})
    Let $\Gamma$ be a tropical curve of genus $2$. Then $\Jac(\Gamma)$ splits if and only if $\Gamma$ covers an elliptic curve. 
\end{theorem}
This relation, however, is not a bijection in the sense that a cover fixes the splitting of $\Jac(\Gamma)$ and a splitting the cover. 
Our next results show that these failings can be remedied by introducing the concept of optimal covers (see Definition \ref{definition_optimalmap}), which is motivated by a similar notion in algebraic geometry.
\begin{theorem}\label{theorem_introduction2.1}(Theorem \ref{theorem_Jacsplits1})
    If $\varphi: \Gamma \rightarrow \mathbb{T}E$ is an optimal cover of degree $d$, then $\mathbb{T}E':=\ker(\varphi_{*})$ is a tropical elliptic curve and $\Jac(\Gamma)$ splits, i.e. there exists an isogeny $\phi:\mathbb{T}E' \bigoplus \mathbb{T}E \rightarrow \Jac(\Gamma)$ whose kernel is isomorphic to $\Jac_d(\mathbb{T}E)$, the group of $d$-torsion points of $\mathbb{T}E$.
\end{theorem}
Theorem \ref{theorem_introduction2.1} should be understood as follows: Optimal coverings allow us to fix a "canonical" representative for a splitting of $\Jac(\Gamma)$. The complement of $\mathbb{T}E$, $\mathbb{T}E'$, enjoys a similar relation with $\Gamma$:
\begin{theorem}\label{theorem_introduction2.2}(Theorem \ref{theorem_Jacsplits2})
   In the setting of Theorem \ref{theorem_introduction2.1}: There exists another optimal cover $\varphi': \Gamma \rightarrow \mathbb{T}E'$ of degree $d$ such that $\varphi'$ interacts "nicely" with $\varphi$.\\
    The term "nicely" is made precise in Theorem \ref{theorem_Jacsplits2}. Informally, $\varphi'$ gives rise to the same splitting $\phi$ whose kernel also satisfies $\ker(\phi)\cong \Jac_d(\mathbb{T}E')$.
\end{theorem}
Taken together, Theorem \ref{theorem_introduction2.1} and Theorem \ref{theorem_introduction2.2} describe the situation in $\mathbb{T}\mathcal{C}$ and $\mathbb{T}\mathcal{A}$, i.e. split Jacobians appear either as a pair of optimal covers in $\mathbb{T}\mathcal{C}$ or as an isogeny satisfying $\Jac_d(\mathbb{T}E) \cong \ker(\phi)\cong \Jac_d(\mathbb{T}E')$ in $\mathbb{T}\mathcal{A}$. Informal, but illustrative is Figure \ref{figure_introduction}.
\begin{figure}[H]
\centering
    \tikzset{every picture/.style={line width=0.75pt}} %set default line width to 0.75pt        

\begin{tikzpicture}[x=0.75pt,y=0.75pt,yscale=-1,xscale=1]

%Straight Lines [id:da6357125880641057] 
\draw    (145.33,94) -- (151.25,102.74) -- (165.88,124.34) ;
\draw [shift={(167,126)}, rotate = 235.9] [color={rgb, 255:red, 0; green, 0; blue, 0 }  ][line width=0.75]    (10.93,-3.29) .. controls (6.95,-1.4) and (3.31,-0.3) .. (0,0) .. controls (3.31,0.3) and (6.95,1.4) .. (10.93,3.29)   ;
%Straight Lines [id:da29831389032130184] 
\draw    (117.33,95) -- (97.48,123.36) ;
\draw [shift={(96.33,125)}, rotate = 304.99] [color={rgb, 255:red, 0; green, 0; blue, 0 }  ][line width=0.75]    (10.93,-3.29) .. controls (6.95,-1.4) and (3.31,-0.3) .. (0,0) .. controls (3.31,0.3) and (6.95,1.4) .. (10.93,3.29)   ;
%Shape: Rectangle [id:dp8392961366255509] 
\draw   (59,31.5) -- (361.33,31.5) -- (361.33,170.5) -- (59,170.5) -- cycle ;
%Straight Lines [id:da2707590075986266] 
\draw    (59,61) -- (360.33,61) ;
%Straight Lines [id:da8737790266579251] 
\draw    (210,31) -- (210.33,171) ;

% Text Node
\draw (128,71.4) node [anchor=north west][inner sep=0.75pt]    {$\Gamma $};
% Text Node
\draw (80,135.4) node [anchor=north west][inner sep=0.75pt]    {$\mathbb{T} E$};
% Text Node
\draw (160,135.4) node [anchor=north west][inner sep=0.75pt]    {$\mathbb{T} E'$};
% Text Node
\draw (90.3,89.93) node [anchor=north west][inner sep=0.75pt]    {$\varphi $};
% Text Node
\draw (158.3,89.93) node [anchor=north west][inner sep=0.75pt]    {$\varphi '$};
% Text Node
\draw (222,101) node [anchor=north west][inner sep=0.75pt]    {$\Jac( \Gamma ) \ \simeq \mathbb{T} E'\ \oplus \mathbb{T} E$};
% Text Node
\draw (108,37.4) node [anchor=north west][inner sep=0.75pt]    {$In\ \mathbb{T}\mathcal{C} :$};
% Text Node
\draw (259,39.4) node [anchor=north west][inner sep=0.75pt]    {$In\ \mathbb{T}\mathcal{A} :$};

\end{tikzpicture}

    \caption{Split Jacobians in \ensuremath{\mathbb{T}\mathcal{C}} and in \ensuremath{\mathbb{T}\mathcal{A}}.}
    \label{figure_introduction}
\end{figure}
Breaking down the right-hand side of Figure \ref{figure_introduction} further, reveals the symmetry which the left-hand side already suggests. We have the diagram (see Diagram \ref{diagram_fromctocultimately}):
\begin{diagram}
     %\begin{tikzcd}[row sep=huge]
     \mathbb{T}E' \ar[dd,"m_d"] \ar[dr,"\varphi^{'*}",sloped] \ar[r,"\iota_1", hook] & \mathbb{T}E' \bigoplus \mathbb{T}E \ar[d,dashed,"{\phi}" description] &\mathbb{T}E \ar[dl,"\varphi^{*}",sloped] \ar[l,"\iota_1", hook', swap] \ar[dd,"m_d"]\\
        & \ar[dl,"\varphi'_{*}",sloped] \Jac(\Gamma) \ar[d,dashed,"{\tilde{\phi}}" description] \ar[dr,"\varphi_{*}" ] &  \\
        \mathbb{T}E'  & \ar[l, "p_1"] \mathbb{T}E' \bigoplus \mathbb{T}E \ar[r, "p_2"]  &\mathbb{T}E  
   % \end{tikzcd}
\end{diagram}
 where $\iota_i$ and $p_i$ are the canonical injections, respectively projections, and $m_d$ is the componentwise multiplication-by-$d$ map.
 
We dedicate a separate analysis to $\varphi_*$ and its dual $\varphi^*$: For covers of the form $\varphi: \Gamma \rightarrow \mathbb{T}E$, where $\Gamma$ is of genus 2, see Subsection \ref{subsection_push-forward} and for covers without any restriction on the genus of $\Gamma$ see Subsection \ref{subsection_pull-back}. This analysis is detached from split Jacobians, but still focuses on covers of elliptic curves. One reason for this is that the target $\mathbb{T}E$ may be viewed as an object of $\mathbb{T}\mathcal{A}$ as well. This "double-identity" allows us to invoke the universal property of the Jacobian (\cite{röhrle2024tropicalngonalconstruction}, Proposition 4.14) and relate factorizations of $\varphi$ to factorizations of $\varphi_*$.

An effort is made throughout to balance abstract and algorithmic techniques: Computational results concerning optimal covers, $\varphi_*$, and $\varphi^*$ are discussed in Subsection \ref{subsection_push-forward} to \ref{subsection_CriteriaforOptimality}. Concerning split Jacobians, Subsection \ref{subsection_algorithmforcomplementarycover} supports Theorems \ref{theorem_introduction2.1} and \ref{theorem_introduction2.2} with an algorithmic point of view. 

The present work relies on technology developed in \cite{MR2457725}, \cite{röhrle2024tropicalngonalconstruction}, \cite{MR4261102}, \cite{MR2275625} and many more for handling tavs efficiently, and which we develop further in Subsections \ref{subsection_morecatconstr}, \ref{subsection_factorizationinTt} and \ref{subsection_exactnessanddualization} to meet our needs.
\subsection{Plan of the paper}
Split Jacobians are studied in $\mathbb{T}\mathcal{A}$ and $\mathbb{T}\mathcal{C}$. This exposes them to two different theories. We treat these separately at first: See Sections \ref{section_catofrealtori} and \ref{section_catoftav} for the theory of tavs and Section \ref{section_catoftc} for the theory of curves. Section \ref{section_crossingbridges} brings them together and sets the stage for the remainder of the paper.
 
We begin Section \ref{section_catofrealtori} with a reminder of real tori with integral structure. 
These form a category, denoted by $\mathbb{T}\mathcal{T}$, that $\mathbb{T}\mathcal{A}$ is built on. 
Still in $\mathbb{T}\mathcal{T}$, Subsections \ref{subsection_morecatconstr} and \ref{subsection_factorizationinTt} extend the existing theory with further concepts relevant to split Jacobians.
 Tropical abelian varieties are discussed in Section \ref{section_catoftav}, starting with preliminaries in Subsection \ref{subsection_preliminariescatofabelvar}. In \ref{subsection_exactnessanddualization} we transfer the notions developed in \ref{subsection_morecatconstr} and \ref{subsection_factorizationinTt} to $\mathbb{T}\mathcal{A}$. Here, we also define exact sequences and study their behaviour under dualization. This concludes the part on tavs. 
 
 Section \ref{section_catoftc} is devoted to the category of tropical curves, recalling background on curves and covers in Subsection \ref{subsection_preliminariescatoftc}. We define optimal tropical covers in Subsection \ref{subsection_optimalcovers}.
 
 Section \ref{section_crossingbridges} establishes a connection between $\mathbb{T}\mathcal{C}$ and $\mathbb{T}\mathcal{A}$, one that we explore in \ref{subsection_curvesofg2coverg1} in the setting of curves of genus 2 covering curves of genus 1. Subsections \ref{subsection_push-forward} and \ref{subsection_pull-back} are a blend of computational and abstract results on $\varphi_*$ and $\varphi^*$, merged together in Subsection \ref{subsection_CriteriaforOptimality} to obtain criteria for optimality. 
 
 Section \ref{section_tropicalsplitJacobians}, finally, is devoted to the phenomenon of tropical split Jacobians. Here we prove Theorem \ref{theorem_introduction1}, \ref{theorem_introduction2.1} and \ref{theorem_introduction2.2}, our main abstract results, which we accompany with examples and make algorithmically accessible in Subsection \ref{subsection_algorithmforcomplementarycover}.
 
\subsection{Future work} In an upcoming paper we take a different approach. Instead of looking at the phenomenon of split Jacobian as a whole, we analyze its building blocks, a pair of elliptic curves together with a finite subgroup of their product, and how to reassemble them into a Jacobian.
\paragraph{\emph{Acknowledgments}.} I am indebted to Hannah Markwig for numerous discussions, advice, and support. Significant improvements of the present work are due to her. I am also grateful for the works of Len, Mikhalkin, Röhrle, Ulirsch, Zakharov, Zharkov and many more and the inspiration they have provided. I thank Felix Röhrle for many interesting discussions and for his valuable feedback on an earlier version of this draft. I thank Thomas Blomme whose original and canny method for computing floor multiplicities naturally led to a question about Jacobians of curves of genus 2. Finally, thanks are also due to the TRR195 for partial support of this work. 
\section{The category of real tori with integral structure}\label{section_catofrealtori}
\subsection{Preliminaries}\label{subsection_preliminariescatofrealtori}
 We start with a description of a surrounding category $\gls{TT}$, the category of \emph{real tori with integral structure}.
 
 \emph{Objects}: More precisely, we start with a description of its objects (see \cite{MR4382460}, Section 2.3 and \cite{röhrle2024tropicalngonalconstruction}, Section 4): The data $\gls{Sigma}:=(\Lambda,\Lambda',[\cdot,\cdot])$, where
\begin{itemize}
\item $\Lambda$ and $\Lambda'$ are finitely generated free abelian groups of the same rank,
    \item $[\cdot,\cdot]: \Lambda \times \Lambda' \rightarrow \mathbb{R}$ is a non-degenerate pairing,
\end{itemize}
 contains all the information needed for building a \emph{real torus with integral structure} and will therefore be identified as such. The torus is given by the quotient $\Hom(\Lambda,\mathbb{R})/\Lambda'$, where the pairing $[\cdot,\cdot]$ encodes the way in which $\Lambda'$ is to be viewed as a lattice in $\Hom(\Lambda,\mathbb{R})$, namely via the embedding $ \Lambda' \rightarrow \Hom(\Lambda,\mathbb{R}): \lambda' \mapsto [\cdot, \lambda']$. Note that we also have an embedding $ \Lambda \rightarrow \Hom(\Lambda',\mathbb{R}): \lambda \mapsto [\lambda,\cdot]$ which turns $\Lambda$  into a lattice in $\Hom(\Lambda',\mathbb{R})$. This allows for the construction of the \emph{dual torus} $\gls{widecheckSigma}$ which is realized as the quotient $\Hom(\Lambda',\mathbb{R})/\Lambda$ and identified by the data $\widecheck{\Sigma}:=(\Lambda',\Lambda,[\cdot,\cdot]^t)$, where  $[\cdot,\cdot]^t$ denotes the transposed pairing. The \emph{dimension} of $\Sigma$ is the $\mathbb{R}$-vector space dimension of $\Hom(\Lambda,\mathbb{R})$ and equal to $\rk(\Lambda)$ (equivalently  equal to $\rk(\Lambda')$).
 
 \emph{Maps}: Now that we have defined the objects of $\mathbb{T}\mathcal{T}$, we turn to structure-preserving maps.
For $i=1,2$ let $\Sigma_i:=(\Lambda_i,\Lambda_i',[\cdot,\cdot]_i)$ be real tori with integral structure. A pair of group homomorphisms $(f^\#: \Lambda_2 \rightarrow \Lambda_1,f_\#: \Lambda'_1 \rightarrow \Lambda'_2)$ is called a \emph{morphism of integral tori} $\gls{f}$ and denoted by $f:=(f^\#,f_\#)$, if 
\begin{diagram}\label{diagram_equationcompatibilitycondition}
    [f^\#(\lambda_2),\lambda'_1]_1=[\lambda_2,f_\#(\lambda'_1)]_2
\end{diagram}
is satisfied for all $\lambda'_1\in \Lambda'_1$ and $\lambda_2\in \Lambda_2$. From this data we construct a homomorphism of real tori as follows: Consider the $\mathbb{R}-$linear map $\Hom(f^\#): \Hom(\Lambda_1,\mathbb{R}) \rightarrow \Hom(\Lambda_2,\mathbb{R})$ and note that the restriction of $\Hom(f^\#)$ to $\Lambda_1'$ (viewed as lattice via the embedding given by $[\cdot,\cdot]_1$) is exactly $f_\#$. Hence, the compatibility condition (\ref{diagram_equationcompatibilitycondition}) guarantees that $\Hom(f^\#)$ takes $\Lambda_1'$ to $\Lambda_2'$. Passing to the quotients yields the requested morphism of tori.

Just as the data for one torus actually yields two tori, $\Sigma$ and $\widecheck{\Sigma}$, we obtain a second morphism from the data underlying $f$ by considering the transposed pair, $(f_\#,f^\#)=:\gls{widecheckf}$, called the \emph{dual morphism} .

Since properties of the induced map of quotients are encoded as properties of the pair $(f^\#,f_\#)$ (see \cite{röhrle2024tropicalngonalconstruction}, Definition 4.8), we call $f$
\begin{itemize}
    \item \emph{surjective}, if $f^\#$ is injective.
    \item \emph{finite}, if $f_\#$ is injective (equivalently if $[\Lambda_1 : f^\#(\Lambda_2)] < \infty$).
    \item \emph{injective},  is $f$ is finite and $f_\#(\Lambda'_1)$ is saturated in $\Lambda'_2$.
    \item an \emph{isogeny}, if $f$ is surjective and finite.
\end{itemize}
The following categorical constructions have been introduced in e.g. \cite{röhrle2024tropicalngonalconstruction} and \cite{MR4382460}.
\begin{definition}\label{definition_(Co-)KernelImageTori}
    To a morphism $f: \Sigma_1 \rightarrow \Sigma_2$ we associate the integral tori
    \begin{enumerate}
    \item $\gls{Ker}:=(\Lambda_1/\im(f^\#)^{\text{sat}},\ker(f_\#),[\cdot,\cdot]_K)$, where $[\cdot,\cdot]_K$ is the pairing induced by $[\cdot,\cdot]_1$.
    \item $\Coker(f):=(\ker(f^\#),\Lambda'_2/\im(f_\#)^{\text{sat}},[\cdot,\cdot]_C)$, where $[\cdot,\cdot]_C$ is the pairing induced by $[\cdot,\cdot]_2$.
    \item  Im$(f):=\Ker(q)_0$, where $q: \Sigma_2 \rightarrow \Coker(f)$ is the morphism induced by the natural maps on the lattices.
    \end{enumerate}
    As the notation suggests, these are Kernels and Cokernels in the sense of category theory (see \cite{röhrle2024tropicalngonalconstruction}, Proposition 4.7).
    \end{definition}
    \begin{remark}
        Note that the torus built from the data $\Ker(f)_0$ is naturally identified with the connected component of the identity of the kernel of $f$ (viewed as map on the quotients) (see \cite{MR4382460}, Section 2.3). The same applies to $\Ker(q)_0$ and 
        the image of $f$, i.e. $\frac{\Hom(f^{\#})(\Hom(\Lambda_1,\mathbb{R}))}{\Hom(f^{\#})(\Hom(\Lambda_1,\mathbb{R}))\cap \Lambda'_2}$.
    \end{remark}
    \subsection{More categorical constructions}\label{subsection_morecatconstr}
    Here, we provide additional categorical constructions needed in Sections \ref{section_crossingbridges} and \ref{section_tropicalsplitJacobians}.
\begin{lemma}\label{lemma_quotientsoftori}
 Let $\Sigma_2$ be a real torus with integral structure and $\Sigma_1 \xhookrightarrow{i} \Sigma_2$ a subtorus. The quotient group $\Sigma_2/\Sigma_1$ has an integral structure given by $(\ker(i^\#),\Lambda'_2/\im(i_\#),[\cdot,\cdot]_Q)$, where $[\cdot,\cdot]_Q: \ker(i^\#) \times \Lambda'_2/\im(i_\#) \rightarrow \mathbb{R} $ is the pairing induced by $[\cdot,\cdot]_2.$ It is a torus of dimension $\rk(\Lambda_2)-\rk(\Lambda_1)$. 
\end{lemma}
\begin{proof}
    Note that $\Sigma^Q:=(\ker(i^\#),\Lambda'_2/\im(i_\#),[\cdot,\cdot]_Q)$ is an integral torus:
    \begin{itemize}
        \item $\ker(i^\#)$ and $\Lambda'_2/\im(i_\#)$ are lattices: $i$ being injective implies that $\im(i_\#)$ is saturated. Hence, $\Lambda'_2/\im(i_\#)$ is one as well as a direct summand of the lattice $\Lambda_2$.
        \item $\rk(\Lambda'_2/\im(i_\#))=\rk(\Lambda_2)-\rk(\Lambda_1)=\rk(\ker(i^\#))$.
        \item $[\cdot,\cdot]_Q$ is well-defined and non-degenerate since $[\cdot,\cdot]_2$ is. 
    \end{itemize}
    We show that $\Sigma^Q$ can be identified with the group theoretic quotient, i.e. that $\Sigma_2/\Sigma_1 \cong \Sigma^Q $ as groups: Let $V_2:=\Hom(\Lambda_2,\mathbb{R})/\Lambda'_2$ and $V_1:=\Hom(\Lambda_1,\mathbb{R})/\Lambda'_1$. We have 
    \begin{align}
        \Sigma_2/\Sigma_1=V_2/i(V_1) \overset{\Phi_1}{\cong} \frac{V_2/\Hom (i^\#)(V_1)}{\Lambda'_2/\im(i_\#)},
    \end{align}
   where $ \Phi_1$ is constructed from $V_2 \rightarrow \frac{V_2/\Hom (i^\#)(V_1)}{\Lambda'_2/\im(i_\#)}, f \mapsto \bar{f}$  by applying two times the universal property of the quotient group. It is easy to see that $\Phi_1$ is bijective.
   The morphism
   \begin{align}
       V_2 \rightarrow \Hom(\ker(i^\#),\mathbb{R}), f \mapsto f_{|\ker(i^\#)}
   \end{align}
   satisfies $f \mapsto 0$ for $f\in \Hom (i^\#)(V_1)$ and hence, gives rise to a unique morphism
   \begin{align}
      \Phi_2: V_2/\Hom (i^\#)(V_1) \rightarrow \Hom(\ker(i^\#),\mathbb{R} ).
   \end{align}
   Note that $\Phi_2$ is $\mathbb{R}$-linear and injective: Suppose $f_{|\ker(i^\#)}=0$. Since $i^\#$ is surjective we can complete 
   \begin{diagram}
        %\begin{tikzcd}
\Lambda_1  \arrow[dr, dashed]
& \arrow[l, "i^\#", swap]\Lambda_2 \arrow[d, "f"]\\
 & \mathbb{R}
%\end{tikzcd
    \end{diagram}
   to a commutative diagram, hence $f\in \Hom(\ker(i^\#),\mathbb{R} )$. For surjectivity, we take advantage of the additional structure of $\Phi_2$ and conclude the proof with a dimension argument:
   \begin{align}
     &  \dim_{\mathbb{R}}(\Hom(\ker(i^\#),\mathbb{R} ))=\rk(\ker(i^\#))=\rk(\Lambda_2)-\rk(\Lambda_1)\\
     &  \dim_{\mathbb{R}}(V_2/\Hom (i^\#)(V_1))=\dim_{\mathbb{R}}(V_2) -\dim_{\mathbb{R}}(\Hom (i^\#)(V_1))=\rk(\Lambda_2)-\rk(\Lambda_1).
   \end{align}
   
\end{proof}
The quotient $\Sigma_2 \xrightarrow{\pi} \Sigma_2/\Sigma_1$, where $\pi=(\pi^\#,\pi_\#)$ is given by the inclusion, $\pi^\#$, and the quotient map, $\pi_\#$, satisfies a universal property that can be deduced from universal properties of the lattices it is built from.
\begin{lemma}\label{lemma_universalpropertyquotientoftori}
    Let $f:\Sigma_2 \rightarrow \Sigma_3$ be a homomorphism of integral tori such that $f$ vanishes on a subtorus $\Sigma_1  \xhookrightarrow{i} \Sigma_2$. Then there exists a unique homomorphism $g:\Sigma_2/\Sigma_1 \rightarrow \Sigma_3 $ such that 
    \begin{diagram}
        %\begin{tikzcd}
\Sigma_2 \arrow[r, "f"] \arrow[d, "\pi"]
& \Sigma_3 \\
 \Sigma_2/\Sigma_1 \arrow[ur, "g",swap]
%\end{tikzcd
    \end{diagram}
    commutes.
\end{lemma}
\begin{proof}
    We can express $f\circ i=0$ in terms of lattices as $i^\#\circ f^\#=0$ and $f_\# \circ i_\#=0$. By the universal property of the kernel $\ker(i^\#) \xhookrightarrow{\pi^\#} \Lambda_2$ and universal property of the quotient $\Lambda'_2 \xrightarrow{\pi_\#} \Lambda'_2/\im(i_\#)$ there exists unique group homomorphisms $g^\#$ and $g_\#$ such that
    \begin{diagram}
     %\begin{tikzcd}[row sep=huge]
     \Lambda_2    &[5em] \ar[l,"f^\#",swap] \Lambda_3 \ar[dl,"{g^\#}"] &  & \Lambda'_2 \ar[d,"\pi_\#", swap] \ar[r,"f_\#"] &[5em] \Lambda'_3  \\ [5ex]
         \ker(i^\#) \ar[u,"\pi^\#"]  &  & &  \Lambda'_2/\im(i_\#) \ar[ur,"{g_\#}", swap]
   % \end{tikzcd}
\end{diagram}
commute. Then $g:=(g^\#,g_\#)$ is a homomorphism of tori
\begin{equation}
    [\lambda_3,g_\#(\bar{\lambda'_2})]_3=[\lambda_3,f_\#(\lambda'_2)]_3=[f^\#(\lambda_3),\lambda'_2]_2=[g^\#(\lambda_3),\bar{\lambda'_2}]_Q
\end{equation}
that satisfies $f=g\circ \pi$.
\end{proof}
\begin{remark}
    The proof of Lemma \ref{lemma_quotientsoftav} provides inspiration for determining the integral structure of the quotient $\Sigma_2/\Sigma_1$. It is an adaptation of methods from the complex world to the tropical setting. Based on the premise that $\mathbb{T}\mathcal{A}$ behaves like its algebraic counterpart, $(\ker(i^\#),\Lambda'_2/\im(i_\#),[\cdot,\cdot]_Q)$ is therefore lattice representation of $\Sigma_2/\Sigma_1$ that is to be expected. It can also be seen to agree with $\Coker(i)$, where $i: \Sigma_1 \hookrightarrow \Sigma_2$ is the inclusion, which would be the definition expected from a tropical categorical point of view.    
\end{remark}
Since free abelian groups are the building blocks of integral tori, it is not surprising that the notion of a product and a coproduct exists and strongly relies on the corresponding notion in the category of abelian groups.
\begin{definition}\label{definition_productandcoproductsinTt}
Given integral tori $\Sigma_1 $ and $ \Sigma_2$, we define their \emph{product} $\Sigma_1 \otimes \Sigma_2$ as follows: 
\begin{itemize}
    \item The underlying discrete data is $(\Lambda_1 \times \Lambda_2,\Lambda'_1 \times \Lambda'_2,[\cdot,\cdot]_1 + [\cdot,\cdot]_2)$, where $\Lambda_1 \times \Lambda_2$, respectively $\Lambda'_1 \times \Lambda'_2$, denotes the direct product of groups and $[\cdot,\cdot]_1+[\cdot,\cdot]_2$ is given by $((\lambda_1,\lambda_2),(\lambda'_1,\lambda'_2)) \mapsto [\lambda_1,\lambda'_1]_1 + [\lambda_2,\lambda'_2]_2$.
    \item The object  $\Sigma_1 \otimes \Sigma_2$ is equipped with a pair of morphisms, $\pi_1:\Sigma_1 \otimes \Sigma_2 \rightarrow \Sigma_1 $ and $\pi_2:\Sigma_1 \otimes \Sigma_2 \rightarrow \Sigma_2 $, induced by the natural projection and inclusion maps between the lattices. 
\end{itemize}
The resulting object is an integral torus, that is: 
\begin{itemize}
    \item The pairing $[\cdot,\cdot]_1+[\cdot,\cdot]_2$ is non-degenerate.
    \item The projection maps $\pi_1$ and $\pi_2$ are morphism of tori.
\end{itemize}
 The \emph{coproduct} of $\Sigma_1 $ and $ \Sigma_2$ is defined analogously and will be denoted by $\Sigma_1 \oplus \Sigma_2$.
   
\end{definition}
\begin{proof}
The objects, $\Sigma_1 \otimes \Sigma_2$ and $\Sigma_1 \oplus \Sigma_2$, are both integral tori. We still need to verify, though, whether they deserve to be called product, respectively coproduct, in the sense of category theory. For the purpose of illustration we only examine the case $\Sigma_1 \otimes \Sigma_2$, arguments for the coproduct are completely analogous: \\
Let $\Sigma_3$ be a torus with morphisms $f_1$ and $f_2$ into each factor of $\Sigma_1 \otimes \Sigma_2$.  We can complete the diagram
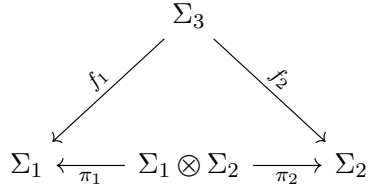
\begin{center} 
    \begin{tikzcd}[row sep=huge]
        & \Sigma_3 \ar[dl,"f_1",sloped] \ar[dr,"f_2",sloped] & \\
        \Sigma_1  & \Sigma_1 \otimes \Sigma_2 \ar[l,"\pi_1"] \ar[r,"\pi_2", swap]  &\Sigma_2 
    \end{tikzcd}
         
 \end{center}
to a commutative diagram by connecting the torus $\Sigma_3$ to $\Sigma_1 \otimes \Sigma_2$ through
\begin{align}
    f: \Sigma_3 \rightarrow \Sigma_1 \otimes \Sigma_2
\end{align} defined by the pair $(f^\#, f_\#)$, where
\begin{itemize}
    \item $f_\#$ is the \emph{unique} group homomorphism that makes 
    
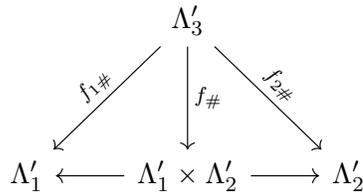
\begin{center} 
    \begin{tikzcd}[row sep=huge]
        & \Lambda'_3 \ar[dl,"f_{1\#}",sloped] \ar[dr,"f_{2\#}",sloped] \ar[d,"f_{\#}"] & \\
        \Lambda'_1  & \Lambda'_1 \times \Lambda'_2 \ar[l] \ar[r,swap]  &\Lambda'_2 
    \end{tikzcd}
         
 \end{center}
    commute. Existence and uniqueness of  $f_\#$ are an immediate consequence of the universal property of the direct products of groups.
    \item $f^\#$ is the \emph{unique} group homomorphism that makes 
    
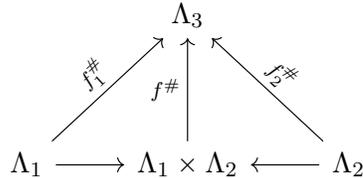
\begin{center} 
    \begin{tikzcd}[row sep=huge]
        & \Lambda_3  & \\
        \Lambda_1 \ar[ur,"f_1^{\#}",sloped] \ar[r,swap] & \Lambda_1 \times \Lambda_2   \ar[u,"f^{\#}"]  &\Lambda_2 \ar[ul,"f_2^{\#}",sloped] \ar[l]
    \end{tikzcd}
    \end{center}
    commute. Existence and uniqueness follow from viewing $\Lambda_1 \times \Lambda_2 $ as coproduct (in $Ab$) and applying the respective universal property.
\end{itemize}
%Since uniqueness of $f$ is obvious, we only need to check whether $f$ itself is a morphism
Since $f$ is clearly unique, we only need to check whether it is a morphism to conclude the argument. But this follows directly from the compatibility of the pair $(f^\#,f_\#)$: We have
\begin{align}
    ([\cdot,\cdot]_1+[\cdot,\cdot]_2)( (\lambda_1,\lambda_2),f_\#(\lambda'_3))&= [\lambda_1,f_{1\#}(\lambda'_3)]_1+ [\lambda_2,f_{2\#}(\lambda'_3)]_2\\
    &= [f_1^{\#}(\lambda_1),\lambda'_3]_3 + [f_2^{\#}(\lambda_2),\lambda'_3]_3 \\
    &= [f_1^{\#}(\lambda_1)+f_2^{\#}(\lambda_2),\lambda'_3]_3\\
    &= [f^\#(\lambda_1,\lambda_2),\lambda'_3]_3
\end{align}
for all $\lambda'_3\in \Lambda'_3$ and $(\lambda_1,\lambda_2) \in \Lambda_1 \times \Lambda_2$. Hence, $\Sigma_1 \otimes \Sigma_2$  has the universal property of the product.
\end{proof}
\begin{remark}
 Although coproduct and product coincide (for finitely many factors), just as it is the case in the category of abelian groups, we insist on using different notations to indicate to the reader when we are using which universal property.
\end{remark}
We conclude by defining \emph{equalizers} and \emph{coequalizers} in $\mathbb{T}\mathcal{T}$: For morphisms \begin{tikzcd}
      \Sigma_1 \arrow[r,shift left,"f"] \arrow[r,shift right,swap,"g"] & \Sigma_2,
\end{tikzcd}
we consider the following pairs $(L,\phi_L)$ and $(C,\phi_C)$, where %$(\phi^\#_L,\phi_{L \#})$, where $\phi^\#_L$ is the projection and \phi_{L \#} inclusion\\
\begin{itemize}
    \item $L:=(\Lambda_L,\Lambda'_L,[\cdot,\cdot]_L)$ and $C:=(\Lambda_C,\Lambda'_C,[\cdot,\cdot]_{CE})$ are integral tori with
\begin{align}
   & \Lambda_L:=\Lambda_1/\im(f^\#-g^\#)^{sat}, \Lambda'_L:=\{ \lambda'_1 \in \Lambda'_1\thinspace:\thinspace f_{\#}(\lambda'_1)=g_{\#}(\lambda'_1) \}\\
   & \Lambda_C:=\{ \lambda_2 \in \Lambda_2\thinspace:\thinspace f^{\#}(\lambda_2)=g^{\#}(\lambda_2) \}, \Lambda'_C:=\Lambda'_2/\im(f_\#-g_\#)^{sat}
\end{align}
and with pairings $[\cdot,\cdot]_L$ and $[\cdot,\cdot]_{CE}$ induced by $[\cdot,\cdot]_1$ and $[\cdot,\cdot]_2$.
\item $\phi_L: L \rightarrow \Sigma_1$ and $\phi_C: \Sigma_2 \rightarrow C$ are the morphisms induced by the natural maps on the lattices.
\end{itemize}
A word on well-definedness, exemplary for $[\cdot,\cdot]_L$: Let $\lambda_1,\lambda_2\in \Lambda_1$ with $\lambda_1-\lambda_2 \in \im(f^\#-g^\#)$ and $\lambda'\in \Lambda'_L$. Using that $f-g$ is a morphism of tori we have:
\begin{align}
    [\underbrace{\lambda_1-\lambda_2}_{=f^\#-g^\#(\lambda_3)}, \lambda']_1=[\lambda_3,f_\#-g_\#(\lambda')]_2=0.
\end{align}
\begin{lemma}\label{lemma_(co)equalizer}
    The category of integral tori $\mathbb{T}\mathcal{T}$ has equalizers and coequalizers. It is finitely complete and cocomplete.
\end{lemma}
\begin{proof}
By the existence theorem for limits (colimits), $\mathbb{T}\mathcal{T}$ is finitely complete (finitely cocomplete), if it has binary products (binary coproducts) and binary equalizers (binary coequalizers). By Definition \ref{definition_productandcoproductsinTt}, $\mathbb{T}\mathcal{T}$ has binary product and coproducts, the existence of equalizers and coequalizers remains to be shown: 
  Let \begin{tikzcd}
      \Sigma_1 \arrow[r,shift left,"f"] \arrow[r,shift right,swap,"g"] & \Sigma_2
\end{tikzcd} be a parallel pair. We claim that their equalizer is given by $(L,\phi_L)$ as constructed above, i.e. that 
\begin{enumerate}
    \item $f$ and $g$ become equal when pulled back to $L$: $f \circ \phi_L= g \circ \phi_L$.
    \item $L$ is universal with this property.
\end{enumerate}
We work on the level of lattices. Note that $(\Lambda'_L, \phi_{L \#})$ is the equalizer of the pair $(f_\#,g_\#)$ in $\gls{Ab}$ and $(\Lambda_L, \phi^\#_{L})$ is related to the coequalizer $(\tilde{\Lambda}_L, \tilde{\phi}_{L})$ of $(f^\#,g^\#)$ in $\gls{Ab}$ as follows: $\Lambda_L$ is the torsion-free quotient group of $\tilde{\Lambda}_L$ and $\phi^\#_{L}$ the composition of $\tilde{\phi}_{L}$ with the natural projection $p:\tilde{\Lambda}_L \rightarrow \Lambda_L$. This yields 
\begin{align}
    f_\# \circ \phi_{L \#} = g_\# \circ \phi_{L \#} \text{ and } \phi^\#_{L} \circ f^\# = \phi^\#_{L} \circ  g^\# 
\end{align}
proving point (1) of the claim. For point (2), let $(\tilde{L}, \phi_{\tilde{L}})$ be another pair satisfying (1). We need to construct a unique arrow $\tilde{L}\xrightarrow{\gamma} L$ with $\phi_{\tilde{L}}=\phi_{L}\circ \gamma.$ % Let $\gamma_\#$ be the unique group homomorphism such that
Since $(\Lambda'_L, \phi_{L \#})$ is the equalizer of $(f_\#,g_\#)$, there exists a unique group homomorphism $\gamma_\#$ such that the left diagram of Figure \ref{figure_equalizerfigure} commutes.
\begin{figure}[H]
    \centering
    \begin{tikzcd}
        \Lambda'_{\tilde{L}} \arrow[d,"\gamma_\#", dashed, swap]  \arrow[r, "\phi_{\tilde{L} \#}"] & \Lambda'_1 \arrow[r,bend left,"f_\#"]  \arrow[r,bend right,swap,"g_\#"] & \Lambda'_2\\
     \Lambda'_L  \arrow[ur,"\phi_{L \#}", swap] & &
    \end{tikzcd}
    \begin{tikzcd}
      \Lambda_{\tilde{L}}  & \arrow[l, "\phi^\#_{\tilde{L}}",swap] \arrow[d,"\tilde{\phi}_{L}"] \Lambda_1  & \arrow[l,bend right,swap,"g^\#"] \arrow[l,bend left,"f^\#"]  \Lambda_2\\
     \Lambda_L  \arrow[u,"\gamma^\#", dashed]  & \arrow[ul,"\tilde{\gamma}^\#",dashed] \tilde{\Lambda}_L \arrow[l,"p"] &
    \end{tikzcd}
    \caption{Applying the universal property of the equalizer (on the left) and the universal property of the quotient and of the coequalizer (on the right) in the category of abelian groups.}
    \label{figure_equalizerfigure}
\end{figure}
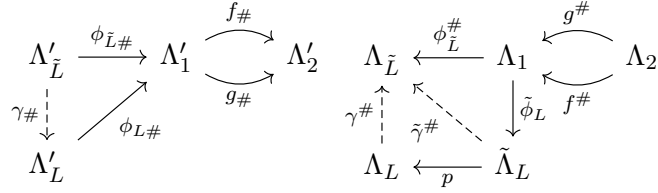
Similarly, we complete the righthand side of Figure \ref{figure_equalizerfigure} to a commutative diagram: The universal property of the coequalizer gives rise to a unique map $\tilde{\gamma}^\#$ whose kernel $\ker(\tilde{\gamma}^\#)$ contains the torsion subgroup of $\tilde{\Lambda}_L$ as $\Lambda_{\tilde{L}}$ is torsion-free. Applying the universal property of the quotient $\Lambda_L$ then yields $\gamma^\#$ as desired. Moreover, $\gamma:=(\gamma^\#,\gamma_\#)$ is a morphism of tori since
\begin{align}
   [\bar{\lambda},\underbrace{\gamma_\#(\lambda')}_{=\phi_{\tilde{L}\#}(\lambda')}]_L=[\lambda,\phi_{\tilde{L}\#}(\lambda') ]_1=
    [\underbrace{\phi^\#_{\tilde{L}}(\lambda)}_{=\gamma^\#(\phi^\#_L(\lambda))},\lambda' ]_{\tilde{L}}=[\gamma^\#(\underbrace{\phi^\#_L(\lambda)}_{=\bar{\lambda}}),\lambda' ]_{\tilde{L}}
\end{align}
holds for all $\lambda \in \Lambda_1 $ and $\lambda'\in \Lambda'_{\tilde{L}}$ concluding point (2).\\
The statements for coequalizers, analogous to points (1) and (2), are verified by the pair $(C,\phi_C)$. The proof is similar and will therefore be omitted.
\end{proof}
\subsection{Factorization}\label{subsection_factorizationinTt}
We use factorizations to gain insight into the structure of a morphism. To that end, we consider factorizations with special properties.
 
\begin{definition}\label{definition_minimalmaximalmaps}
    Let $f: \Sigma_1 \rightarrow \Sigma_2$ be a surjective morphism of integral tori. We call a map $g_1: \Sigma_1 \rightarrow \Sigma_3$ \emph{minimal} for $f$, if $f$ factors through $g_1$ (i.e. $f=g_2\circ g_1$) and for all other factorizations $\tilde{g}_2 \circ \tilde{g}_1$ of $f$, where $\tilde{g}_1: \Sigma_1 \rightarrow \tilde{\Sigma}_3$ and $\dim_\mathbb{R}(\tilde{\Sigma}_3)\leqslant \dim_\mathbb{R}(\Sigma_3)$, there exists a unique map $h_1$ such that the diagram on the left 
     
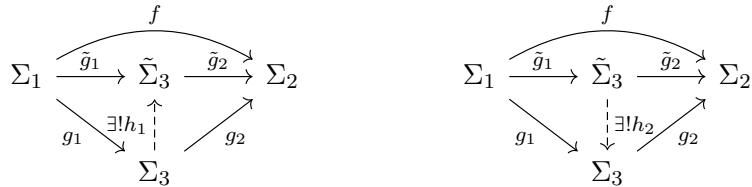
\begin{diagram}%\label{diagram_minimalandmaximalmaps}
     %\begin{tikzcd}[row sep=huge]
     \Sigma_1  \ar[dr,"g_1",swap] \ar[rr, bend left,"f"] \ar[r,"\tilde{g}_1"] & \tilde{\Sigma}_3 \ar[r,"{\tilde{g}_2}" ] & \Sigma_2 & &\Sigma_1  \ar[dr,"g_1",swap] \ar[rr, bend left,"f"] \ar[r,"\tilde{g}_1"] & \tilde{\Sigma}_3 \ar[d,dashed,"\exists ! h_2"] \ar[r,"{\tilde{g}_2}" ] & \Sigma_2 \\
     & \Sigma_3 \ar[ur,"g_2",swap ] \ar[u,dashed,"\exists ! h_1"]& & & & \Sigma_3 \ar[ur,"g_2",swap ] &
   % \end{tikzcd}
\end{diagram}
commutes. In complete analogy, we call $g_1$ \emph{maximal} for $f$, if $f$ factors through $g_1$ and for all other factorizations of $f$ with $\dim_\mathbb{R}(\tilde{\Sigma}_3)\geqslant \dim_\mathbb{R}(\Sigma_3)$ there exists a unique map $h_2$ such that the diagram on the right commutes.
\end{definition}
Let us now bring these theoretical constructs to life by considering the \emph{tropical Stein factorization} (see \cite{MR2062673} for an analytique analogue): 
 Let $f: \Sigma_1 \rightarrow \Sigma_2$ be surjective and consider 
 \begin{diagram}\label{diagram_Steinfactorization}
     %\begin{tikzcd}[row sep=huge]
     \Sigma_1 \ar[rr, bend right,"f",swap] \ar[r,"\pi_f"] & \Sigma_1/\Ker(f)_0 \ar[r,"{\phi_f}" ] & \Sigma_2 
   % \end{tikzcd}
\end{diagram}
    Since $\Sigma_1/\Ker(f)_0$ is a torus of dimension $\dim_\mathbb{R}(\Sigma_2)$ (Lemma \ref{lemma_quotientsoftori}), the projection map $\pi_f$ is a morphism with $\ker(\pi_f)=\Ker(f)_0$ and $\phi_f$ an isogeny. 
\begin{lemma}[Tropical Stein Factorization]\label{lemma_tropicalSteinfactorization}
    Every surjective morphism of tori $f: \Sigma_1 \rightarrow \Sigma_2$ factors canonically through an isogeny and a morphism whose group-theoretic kernel is connected as in (\ref{diagram_Steinfactorization}), such that $\pi_f$ is minimal for $f$ (see Definition \ref{definition_minimalmaximalmaps}).
    \end{lemma}
    Explicitly, (\ref{diagram_Steinfactorization}) satisfies the following universal property: For all factorizations $g' \circ \pi'$ of $f$ such that $g'$ is an isogeny there exists a unique isogeny $\phi'$ such that 
    \begin{diagram}\label{diagram_UnivPropSteinfactorization}
     %\begin{tikzcd}[row sep=huge]
     \Sigma_1  \ar[dr,"\pi_f",swap] \ar[rr, bend left,"f"] \ar[r,"\pi'"] & \Sigma_3 \ar[r,"{g'}" ] & \Sigma_2 \\
     & \Sigma_1/\Ker(f)_0 \ar[ur,"{\phi_f}",swap ] \ar[u,dashed,"\exists ! \phi'"]&
   % \end{tikzcd}
\end{diagram}
commutes.
\begin{proof}
For such a factorization $g' \circ \pi'$ of $f$ we have
\begin{align}
    \gls{ker}=\ker(g' \circ \pi')= \ker(\pi') + (\pi')^{-1}(\ker(g')\setminus \{0\}),
\end{align}
where $\pi'^{-1}(\ker(g')\setminus \{0\})$ is a finite set since $g'$ is an isogeny and $\ker(\pi')$ a (not necessarily connected) group. This implies $\Ker(f)_0 \subset \ker(\pi')$ as $\Ker(f)_0$ is the connected component containing the identity. By the universal property of the quotient (Lemma \ref{lemma_universalpropertyquotientoftori}) there exists a unique morphism of integral tori $\phi'$ such that (\ref{diagram_UnivPropSteinfactorization}) commutes, i.e. $\pi'$ factors through $\pi_f$ and $\phi_f$ through $g'$. Moreover, $\phi'$ is surjective (since $\pi'$ is) with finite kernel (since $\ker(\phi')\subset \ker(\phi_f)$ holds), hence an isogeny.
\end{proof}
\begin{remark}
    In \cite{röhrle2024tropicalngonalconstruction}, Röhrle and Zakharov prove that an isogeny canonically factors into a free isogeny and a dilation. This factorization enjoys a universal property similar to the tropical Stein factorization (see \cite{röhrle2024tropicalngonalconstruction}, Lemma 4.9). Note that both factorizations are complementary: The tropical Stein factorization decomposes \emph{surjective} morphisms, when specialized to isogenies, however, it is trivial.
\end{remark}
As a first application of the tropical Stein factorization we rediscover a fact whose analogue for complex tori is well-known (see \cite{MR1767634}, Proposition 4.6):% and Proof):
\begin{lemma}\label{lemma_c(f)=c(dualf)}
    Let $f:\Sigma \rightarrow \Sigma'$ be a morphism of integral tori. The number of connected components of $\ker(f)$ and $\ker(\widecheck{f})$ agree.
\end{lemma}
\begin{proof}
The tropical Stein factorizations of the surjective morphisms $\Sigma  \xrightarrow{f} \im(f)$ and $\widecheck{\Sigma}'\xrightarrow{\widecheck{f}} \im(\widecheck{f})$ give rise to factorizations of $f$ and $\widecheck{f}$:
\begin{align}
        & \Sigma \xrightarrow{\pi_f} \Sigma/{\Ker(f)}_0 \xrightarrow{\phi_f} \im(f) \hookrightarrow \Sigma'\\
       & \widecheck{\Sigma'} \xrightarrow{\pi_{\widecheck{f}}} \Sigma'/{\Ker(\widecheck{f})}_0 \xrightarrow{\phi_{\widecheck{f}}} \im(\widecheck{f})  \hookrightarrow \Sigma.
    \end{align}
     Note that the number of connected components of each kernel, $c(\ker(f))$ and $c(\ker(\widecheck{f}))$, is given by $|\ker(\phi_f)|$ and $|\ker(\phi_{\widecheck{f}})|$) (i.e. by the geometric degree of the respective isogeny, $\phi_f$ and $\phi_{\widecheck{f}}$  (see \cite{röhrle2024tropicalngonalconstruction})). In particular, it is finite. We can dualize the first sequence to obtain another factorization of $\widecheck{f}$
    \begin{align}
        \widecheck{\Sigma'} \twoheadrightarrow \widecheck{\im}(f) \xrightarrow{\widecheck{\phi}_f} \widecheck{(\Sigma/{\Ker(f)}_0)} \rightarrow \widecheck{\Sigma},
    \end{align}
    where the last map factors through $\im(\widecheck{f})$ (simply because it is a factorization of \nolinebreak $\widecheck{f}$)
    \begin{align}
        \widecheck{\Sigma}' \twoheadrightarrow \widecheck{\im}(f) \xrightarrow{\widecheck{\phi}_f} \widecheck{(\Sigma/{\Ker(f)}_0)} \xrightarrow{g} \im(\widecheck{f}) \hookrightarrow \widecheck{\Sigma}.
    \end{align}
    Since $g$ is finite as (the restriction of the) dual of the surjective map $\pi_f$ and surjective (artificially made so by restricting the target), it is an isogeny. With $g$ being one, the composition $g \circ \widecheck{\phi}_f$ is one as well. So minimality of the tropical Stein factorization (Lemma \ref{lemma_tropicalSteinfactorization}) forces $\phi_{\widecheck{f}}$ to factor through $g \circ \widecheck{\phi}_f$ and therefore $\ker(\phi_{\widecheck{f}})\supset \ker(\widecheck{\phi}_f)$ holds. Using $\ker(\widecheck{\phi}_f)\cong \ker(\phi_f)$ for isogenies, we can conclude that $c(\ker(f))\leqslant c(\ker(\widecheck{f}))$. Finally, "$\geqslant$" follows from $\widecheck{\widecheck{f}}=f.$
\end{proof}
\section{The category of tropical abelian varieties}\label{section_catoftav}
\subsection{Preliminaries}\label{subsection_preliminariescatofabelvar}
By enhancing an integral torus with the data of a polarization, we create a new object:
\begin{definition}\label{definition_tav}
A \emph{tropical abelian variety (tav)} $ \Sigma=(\Lambda,\Lambda',[\cdot,\cdot])$ is a real torus with integral structure together with a \emph{polarization}, i.e. a group homomorphism $\gls{zeta}: \Lambda' \rightarrow \Lambda$  such that the bilinear form $[\zeta(\cdot),\cdot]: \Lambda' \times \Lambda' \rightarrow \mathbb{R}$ is symmetric and positive definite. If $\zeta$ is a bijection, we call $\zeta$ a \emph{principal polarization} (pp) and $ \Sigma$ a \emph{principally polarized tropical abelian variety (pptav)}. In any case, the \emph{dimension} of $\Sigma$ is the $\mathbb{R}$-vector space dimension of $\Hom(\Lambda,\mathbb{R})$ and equal to $\rk(\Lambda)$ (equivalently  equal to $\rk(\Lambda')$).
\end{definition}
A polarization $\zeta$ defines an isogeny $f_\zeta:=(\zeta,\zeta)$ between $\Sigma$ and its dual $ \widecheck{\Sigma}$, which, endowed with the so-called \emph{dual polarization} $\widecheck{\zeta}$ (previously defined by Röhrle and Zakharov, see first version of \cite{röhrle2024tropicalngonalconstruction}), is a tav as well. We can characterize its kernel, $\ker(\gls{fzeta})$, using the \emph{type} of $\zeta$, which is given by the invariant factors $(\alpha_1,...,\alpha_n)$ (where $n:=\rk(\Lambda)$) of its Smith normal form. We have:
\begin{align}
    \ker(\gls{fzeta})\cong \mathbb{Z}/\alpha_1\mathbb{Z} \times ... \times \mathbb{Z}/\alpha_n\mathbb{Z}.
\end{align}
If $(\lambda'_i)^n_{i=1}$ and $(\lambda_i)^n_{i=1}$ are lattice basis of $\Lambda'$ and $\Lambda$, respectively, such that $\zeta(\lambda'_i)= \alpha_i \lambda_i$, then $\widecheck{\zeta}$ is defined by $\lambda_i \mapsto \frac{\alpha_1\alpha_n}{\alpha_i}\lambda'_i$.

\emph{Morphism of tav} are not required to satisfy any additional conditions, i.e. the properties of morphisms of real tori with integral structures are, in a sense, sufficient: For example, they are, as is generally the case, a source of new objects (see Definition \ref{definition_(Co-)KernelImage}). To this end, let us first look at how polarizations can be transported along morphisms.

\begin{definition}\label{definition_induced/pf/pbpolarizationandpolarizedisogeny}
Let $\Sigma_2$ be a tav with polarization $\zeta_2$, $\Sigma_1$ a real torus with integral structure, and $\phi: \Sigma_1 \rightarrow \Sigma_2$ an isogeny. Then $\phi^* \zeta_2:=\phi^\#\circ \zeta_2 \circ \phi_\#$ is a polarization on $\Sigma_1 $ and called the \emph{induced polarization} or alternatively \emph{the pull-back} of $\zeta_2$ by $\phi$. Conversely, suppose $\Sigma_1$ carries a polarization $\zeta_1$. We can define the \emph{push-forward} of $\zeta_1$ by $\phi$ as $\phi_* \zeta_1:=\widecheck{\zeta}$, where $\zeta:= \widecheck{\phi}^* \widecheck{\zeta_1}$. We say that an isogeny $\phi: \Sigma_1 \rightarrow \Sigma_2$ is \emph{polarized} with respect to polarizations $\zeta_1$ on $\Sigma_1$ and $\zeta_2$ on $\Sigma_2$, if $\zeta_1$ is the polarization induced by $\phi$ and $\zeta_2$, in other words, if the diagram 
 \begin{center}
  
\begin{tikzcd}
\Sigma_1 \arrow[r, "\phi"] \arrow[d, "{f_{\zeta_1}}"]
& \Sigma_2 \arrow[d, "{f_{\zeta_2}}"] \\
\widecheck{\Sigma}_1
& \arrow[l, "\widecheck{\phi}"] \widecheck{\Sigma}_2
\end{tikzcd}
 \end{center}
commutes.
    
\end{definition}
Note that pulling back polarizations always turns $\phi$ into a polarized isogeny. The same is not true for the push-forward.
\begin{remark}
\hfill
\begin{itemize}
    \item The pull-back of a polarization has been introduced in \cite{MR4382460} and \cite{röhrle2024tropicalngonalconstruction}. Another notion of push-forward can be extracted from Proposition 4.11 in \cite{röhrle2024tropicalngonalconstruction} by using the identification
    \begin{align}
        \Hom(\Lambda'_1, \Lambda_1) \cong \Hom(\Lambda_1,\mathbb{Z})\otimes \Lambda'_1 \big{(}\cong \text{H}_{1,1}(\Sigma_1) \big{)}
    \end{align}
    via $\zeta_1 \longmapsto \sum^n_{i=1}\Hom(\zeta_1(\lambda_i'),\mathbb{Z})\otimes \lambda_i'$, where $(\lambda'_i)_i$ is a lattice basis of $\Lambda'_1$, and using the natural notion of push-forward in this setting, which is given by $f_*:= \Hom(f^{\#},\mathbb{Z}) \otimes f_{\#}$.
    \item Isogenies are special in that sense that they allow for bidirectional transport of polarizations. However, a surjective morphism is enough for the push-forward and a finite one for the pull-back. 
\end{itemize}
\end{remark}

\begin{definition}\label{definition_(Co-)KernelImage}(\cite{MR4261102}, Proposition 1.1.2)
    To a morphism $f: \Sigma_1 \rightarrow \Sigma_2$ we associate the following tropical abelian varieties:
    \begin{enumerate}
    \item $\Ker(f)_0$ (Definition \ref{definition_(Co-)KernelImageTori}) with polarization $\zeta_K:=i^*\zeta_1$, where $i:\Ker(f)_0 \hookrightarrow \Sigma_1$ is the inclusion.
    \item $\Coker(f)$ (Definition \ref{definition_(Co-)KernelImageTori}) with polarization $\zeta_C:=q_*\zeta_2$, where $q: \Sigma_2 \rightarrow \Coker(f) $ is the quotient map. 
      \end{enumerate}
\end{definition}
Similarly, we can transfer the notions of (co-)products and (co-)equalizers to $\mathbb{T}\mathcal{A}$: To a pair of morphism \begin{tikzcd}
      \Sigma_1 \arrow[r,shift left,"f"] \arrow[r,shift right,swap,"g"] & \Sigma_2,
\end{tikzcd} we associate the following tropical abelian varieties:
    \begin{enumerate}
    \item $(L,\phi_L)$ (Lemma  \ref{lemma_(co)equalizer}) with polarization $\zeta_L:=\phi_L^*\zeta_1$, the \emph{equalizer} of the pair $(f,g)$.
    \item $\Coker(f)$ (Lemma  \ref{lemma_(co)equalizer}) with polarization $\zeta_C:=\phi_{C*}\zeta_2$, the \emph{coequalizer} of the pair $(f,g)$.
      \end{enumerate}
\begin{definition}\label{definition_productandcoproducts}
Given tavs $\Sigma_1 $ and $ \Sigma_2$, we define their \emph{product} $\Sigma_1 \otimes \Sigma_2$ as follows: 
\begin{itemize}
    \item The underlying real torus with integral structure is given by the product $\Sigma_1 \otimes \Sigma_2$ in $\mathbb{T}\mathcal{T}$ (see Definition \ref{definition_productandcoproductsinTt}).
    \item  The polarization  $\zeta_1 \times \zeta_2$ is defined component-wise.
\end{itemize}
The resulting object is a tav, that is the group homomorphism $\zeta_1 \times \zeta_2$ satisfies the condition described in Definition \ref{definition_tav}. The \emph{coproduct} of $\Sigma_1 $ and $ \Sigma_2$ is defined analogously and will be denoted by $\Sigma_1 \oplus \Sigma_2$.
   
\end{definition}

Putting all these pieces together we see that the category of tavs, $\gls{TA}$, has a certain amount of structure. It makes it the perfect setting for several categorical constructions. %opens the doorway for constructions
\begin{lemma}
    The category of polarized tropical abelian varieties, $\mathbb{T}\mathcal{A}$, is additive, i.e. 
    \begin{enumerate}
        \item it contains a zero object given by the trivial object: $(\{0\},\{0\},[\cdot,\cdot])$.
        \item it has an abelian group structure on the Hom-Sets, which is given by the usual addition.
        \item for any two objects, the product and coproduct are defined as in Definition \ref{definition_productandcoproducts}.
    \end{enumerate}
    Moreover, for every morphism $f$
    \begin{itemize}
        \item cokernel and kernel exist (Definition \ref{definition_(Co-)KernelImage}).
        \item such that if $\Ker(f)_0=0$, then $f$ is the kernel of its cokernel.
        \item such that if $\Coker(f)_0=0$, then $f$ is the cokernel of its kernel.
    \end{itemize}
    Hence, $\mathbb{T}\mathcal{A}$ is abelian.
    
\end{lemma}
This has already been remarked in \cite{MR4261102}. The existence of binary coequalizers and equalizers further shows: $\mathbb{T}\mathcal{A}$ is finitely complete, respectively finitely cocomplete. 
\subsection{Exactness and Dualization}\label{subsection_exactnessanddualization}
Although $\mathbb{T}\mathcal{A}$ has a notion of kernels, these are not suitable for defining exact sequences. We use group-theoretic kernels instead and define:
\begin{definition}
    We call a sequence 
     \begin{align}\label{shortexactsequence}
        0 \rightarrow \Sigma_1 \xrightarrow{f} \Sigma_2 \xrightarrow{g} \Sigma_3 \rightarrow 0,
    \end{align} of tropical abelian varieties \emph{exact}, if $f$ is injective, $g$ is surjective, and $\im(f)=\ker(g)$ holds. Note that in particular $\ker(g)$ is connected and itself a tropical abelian variety.  
\end{definition}
Warning: In contrast to the category of abelian groups not every morphism of tavs gives rise to a short exact sequence.
\begin{lemma}\label{lemma_dualizingexactsequences}
   Let   
    \begin{align}
        0 \rightarrow \Sigma_1 \xrightarrow{f} \Sigma_2 \xrightarrow{g} \Sigma_3 \rightarrow 0.
    \end{align}
    be a short exact sequence of tropical abelian varieties.
    Then the dual sequence
    \begin{align}
        0 \rightarrow \widecheck{\Sigma}_3 \xrightarrow{\widecheck{g}} \widecheck{\Sigma}_2 \xrightarrow{\widecheck{f}} \widecheck{\Sigma}_1 \rightarrow 0
    \end{align}
    is exact. In other words 
    \begin{align}
        \widecheck{\cdot}:\mathbb{T}\mathcal{A} \rightarrow \mathbb{T}\mathcal{A}
    \end{align}
    is an exact functor.
\end{lemma}
\begin{proof}
    Suppose $\Sigma_i=(\Lambda_i,\Lambda'_i,[,]_i )$ and denote by $V_i$ the universal cover $\Hom(\Lambda_i, \mathbb{R})$ of $\Sigma_i$  for $i=1,2,3$. We construct the commutative diagram in Figure \ref{figure_snakelemma} as follows:
   \begin{figure}[ht]
        \centering
         \begin{tikzpicture}[>=triangle 60]
\matrix[matrix of math nodes,column sep={60pt,between origins},row
sep={60pt,between origins},nodes={asymmetrical rectangle}] (s)
{
|[name=03]| 0 &|[name=ka]| \Lambda'_1 &|[name=kb]| \Lambda'_2 &|[name=kc]| \Lambda'_3 \\
|[name=04]| 0 &|[name=A]| V_1 &|[name=B]| V_2 &|[name=C]| V_3 &|[name=01]| 0 \\
|[name=02]| 0 &|[name=A']|\Sigma_1 &|[name=B']| \Sigma_2 &|[name=C']| \Sigma_3 &|[name=no]| 0 \\
&|[name=ca]| 0 &|[name=cb]| 0 &|[name=cc]| 0 \\
};
\draw[->] 
          (ka) edge (A)
          (kb) edge (B)
          (kc) edge (C)
          (A) edge node[auto] {$F$} (B)
          (B) edge node[auto] {$G$} (C)
          (C) edge (01)
          (A) edge node[auto] {$\pi_1$} (A')
          (B) edge node[auto] {$\pi_2$} (B')
          (C) edge node[auto] {$\pi_3$} (C')
          (02) edge (A')
          (A') edge node[auto] {$f$} (B')
          (B') edge node[auto] {$g$} (C')
          (A') edge (ca)
          (B') edge (cb)
          (C') edge (cc)
          (C') edge (no)
          (03) edge (ka)
          (04) edge (A)
;
\draw[->,gray] (ka) edge (kb)
               (kb) edge (kc)
               (ca) edge (cb)
               (cb) edge (cc)
;
\draw[->,gray,rounded corners] (kc) -| node[auto,text=black,pos=.7]
{\(\partial\)} ($(01.east)+(.5,0)$) |- ($(B)!.35!(B')$) -|
($(02.west)+(-.5,0)$) |- (ca);
\end{tikzpicture}
\caption{Applying the snake lemma in the proof of Lemma \ref{lemma_dualizingexactsequences}.}
\label{figure_snakelemma}
        \end{figure} 
    
   Lift the composition $V_1\xrightarrow{\pi_1} \Sigma_1 \xrightarrow{f} \Sigma_2$ to a map $F$ into the universal covering $V_2$ of $\Sigma_2$ such that the respective square commutes and $F(0)=0$. In fact $F$ is given by $\Hom(f^\#)$ and $\mathbb{R}$-linear.  Hence, $\ker(F)$ is a vector space that satisfies $\ker(F)\subset \Lambda_1$. It follows $\ker(F)=\{0\}$. Proceed analogously for $g$ to fill in the second row with an exact sequence $0 \rightarrow V_1 \rightarrow V_2 \rightarrow V_3 \rightarrow 0$.\\
    Since the category of abelian groups is abelian, we can use the snake Lemma to get an exact sequence of kernels and cokernels (see Figure \ref{figure_snakelemma}): 
    \begin{align}
        0 \rightarrow \Lambda'_1 \xrightarrow{f_\#} \Lambda'_2 \xrightarrow{g_\#} \Lambda'_3 \rightarrow 0.
    \end{align}
    Applying the exact functor $\Hom(\cdot, \mathbb{R})$ yields an exact sequence between the universal covers of the dual varieties 
    \begin{align}\label{exactsequenceof universalcovers}
        0 \rightarrow \Hom(\Lambda'_3, \mathbb{R}) \xrightarrow{\Hom(g_\#)} \Hom(\Lambda'_2, \mathbb{R}) \xrightarrow{\Hom(f_\#)} \Hom(\Lambda'_1, \mathbb{R}) \rightarrow 0.
    \end{align}
    The maps factor through the quotient. By pushing (\ref{exactsequenceof universalcovers}) down we get 
    \begin{align} 
    \widecheck{\Sigma}_3 \xrightarrow{\widecheck{g}} \widecheck{\Sigma}_2 \xrightarrow{\widecheck{f}} \widecheck{\Sigma}_1 \rightarrow 0.
    \end{align}
    We still need to check whether $\widecheck{g}$ is injective, i.e. $g^\#(\Lambda_3)$ is saturated in $\Lambda_2$: 
    Apply $\Hom(\cdot, \mathbb{R})$ to the second row of \ref{figure_snakelemma} to obtain
    \begin{align}
        0 \rightarrow \Lambda_3 \rightarrow \Lambda_2 \rightarrow \Lambda_1 \rightarrow 0.
    \end{align}
    Since the $\Lambda_i$ are free abelian groups, the sequence splits and $\Lambda_3$ is a direct summand of $\Lambda_2$. This finishes the proof. 
\end{proof}
 \paragraph{\emph{Quotients of tavs.}}
We can always extend $0 \rightarrow \Sigma' \rightarrow \Sigma$ to a short exact sequence $0 \rightarrow \Sigma' \rightarrow \Sigma \rightarrow \Sigma'' \rightarrow 0$ in the category of \emph{abelian groups}. Whether we can do so in the category of \emph{tropical abelian varieties} is a priori not clear. We address this question in the following Lemmas, adapting methods from the complex world to the tropical setting. 
\begin{lemma}\label{lemma_quotientsoftav}
    Let $\Sigma$ be an abelian variety and suppose $\Sigma' \xrightarrow{i} \Sigma$ is an abelian subvariety. Then there exists an abelian variety $\Sigma''$ and a morphism $ \Sigma \twoheadrightarrow \Sigma''$ whose kernel is $\Sigma'$.
       
\end{lemma}
\begin{proof}
Recall that we have a relatively good handle on when $ \Sigma_1 \rightarrow \Sigma_2 \rightarrow 0$  extends (in the category of tori!) to the left. This is precisely the case when the kernel is connected.\\
Seizing the opportunity, we make our lives easier by working with the dual instead: The map $\widecheck{i}: \widecheck{\Sigma} \rightarrow \widecheck{\Sigma}'$ is surjective and gives rise to a short exact sequence in the category of \emph{groups}
    \begin{align}\label{sequenceinlemma_quotientsoftav}
      0 \rightarrow \ker(\widecheck{i}) \rightarrow \widecheck{\Sigma} \rightarrow \widecheck{\Sigma}' \rightarrow 0.
    \end{align}
     We claim that (\ref{sequenceinlemma_quotientsoftav}) lives in the category of \emph{tropical abelian varieties}, i.e. that $\ker(\widecheck{i})$ is connected. Denote by $c(\ker(i))$ the number of connected components of $\ker(i)$. As $i$ is injective, $\ker(i)$ is trivial and $1=c(\ker(i))=c(\ker(\widecheck{i}))$ by Lemma \ref{lemma_c(f)=c(dualf)}. All that remains to be done is the translation back to our original setting: Apply the dualization functor to obtain a short exact sequence of tavs (Lemma \ref{lemma_dualizingexactsequences})
     \begin{align}
      0 \rightarrow  \Sigma' \rightarrow \Sigma \rightarrow \widecheck{\ker}(\widecheck{i}) \rightarrow 0
    \end{align}
    and set $ \Sigma'':=\widecheck{\ker}(\widecheck{i})$. 
    \end{proof}
    \begin{remark}\label{remark_latticerepresentationofquotient}
    From the preceding proof we can extract a concrete description of the quotient in terms of lattices. We have $\Sigma''=(\ker(i^\#),\Lambda'_2/\im(i_\#),[\cdot,\cdot]^t_K)$, where $i$ denotes the inclusion map, $\Sigma'=(\Lambda'_1,\Lambda'_2,[\cdot,\cdot]')$, and $[\cdot,\cdot]_K$ is the pairing induced by the pairing $[\cdot,\cdot]$ on $\Sigma$, and recognize the quotient of integral tori from Lemma \ref{lemma_quotientsoftori}, but now automatically equipped with a polarization.
    \end{remark}
    There is another kind of quotient that creates an object in $\mathbb{T}\mathcal{A}$, the group-theoretic quotient of a tav by a finite group.
\begin{lemma}\label{lemma_quotientoftavbyfinitesubgroup}
Let $\Sigma$ be a tav and suppose $G \subset \Sigma$ is a finite subgroup. Then there exists a tav $\Sigma_G$ and a free isogeny $ \Sigma \twoheadrightarrow \Sigma_G$ whose kernel is $G$.
\end{lemma}
\begin{proof}
    Let $\Hom(\Lambda,\mathbb{R}) \xrightarrow{\pi} \Sigma$ be the universal covering of $\Sigma$. Then $\pi^{-1}(G)$ is a lattice containing $\Lambda'$ such that $[\pi^{-1}(G):\Lambda']< \infty$. The embedding $\pi^{-1}(G) \xhookrightarrow{i} \Hom(\Lambda,\mathbb{R})$ induces a non-degenerate pairing
    \begin{align}
       [\cdot,\cdot]_G:  \Lambda \times \pi^{-1}(G) \rightarrow \mathbb{R}, (\lambda,\lambda') \rightarrow [\lambda,\lambda']_G:=i(\lambda')(\lambda)
    \end{align}
        which identifies the quotient $\Sigma_G:=\Sigma/G$ as real torus with integral structure built from $(\Lambda,\pi^{-1}(G), [\cdot,\cdot]_G)$. Upgrading $\Sigma_G$ to a tav requires the data of a polarization. To that extend, note that the quotient map $q: \Sigma \rightarrow \Sigma_G $ is an isogeny so its dual $\widecheck{q}$ is one, too. Hence, we can use $\widecheck{q}$ to define a polarization on $\widecheck{\Sigma}_G$ as the induced polarization $\widecheck{\zeta}_G:=\widecheck{q}^*\widecheck{\zeta}$ and set $\zeta_G:=\widecheck{\widecheck{\zeta}}_G$, i.e. $\zeta_G=q_*(\zeta)$.
\end{proof}
\begin{remark}
Note that $q$ is not necessarily polarized with respect to $\zeta$ and $\zeta_G$ since ${q}^*(\zeta_G)$ and $\zeta$ do not agree in general. This will be discussed further in a forthcoming paper.
\end{remark}
\paragraph{\emph{Factorization of morphisms of abelian varieties.}}
Since morphisms of tavs are not required to satisfy any additional properties, we can transfer definitions and results of Subsection \ref{subsection_factorizationinTt} to $\mathbb{T}\mathcal{A}$. In particular, we have that the tropical Stein factorization of a surjective morphism of tavs
 \begin{diagram}
     %\begin{tikzcd}[row sep=huge]
     \Sigma_1 \ar[rr, bend right,"f",swap] \ar[r,"\pi_f"] & \Sigma_1/\Ker(f)_0 \ar[r,"{\phi_f}" ] & \Sigma_2 
   % \end{tikzcd}
\end{diagram}
is legitimate in  $\mathbb{T}\mathcal{A}$ as $\Sigma_1/\Ker(f)_0$ is a tav by Lemma \ref{lemma_quotientsoftav}.

\section{The category of Tropical Curves}\label{section_catoftc}
\subsection{Preliminaries}\label{subsection_preliminariescatoftc}
In this section we describe the \emph{category of tropical curves}, $\mathbb{T}\mathcal{C}$. We base our exposition on \cite{MR3375652} and \cite{MR2772537}, drawing on the foundational work of Mikhalkin in \cite{MR2275625} and of Mikhalkin and Zharkov in \cite{MR2457739}. Objects of $\mathbb{T}\mathcal{C}$ are \emph{tropical curves}. These are, in analogy to the complex setting, (connected) topological spaces  homeomorphic to a locally finite 1-dimensional simplicial complex and carry a tropical structure (\cite{MR2457739}, Definition 3.1).
An equivalent way to package the information of a tropical curve is by means of a metric graph (\cite{MR2457739}, Section 3.3). This is the data we want to work with.
\begin{definition}\label{definition_tropicalcurve}
A tropical curve $\Gamma$ is the geometric realization of a metric graph $(G,l)$, i.e. a finite graph $G$ (with no legs/ends) together with a function $l: E(G) \rightarrow \mathbb{R}_{>0}$: It is the topological quotient
    \begin{align}
     \Gamma := \bigcup^{\sbullet}_{e\in E(G)} [0,l(e)]/ \sim    \end{align}
   endowed with the path metric and with equivalence relation $\sim$ coming from the incidence relations of $G$. 
   Any $(G',l')$ that is obtained from $(G,l)$ by adding or deleting vertices of valence $2$ (and adapting the length function accordingly) is called a \emph{model} for $\Gamma$ and $G$ its \emph{combinatorial type}. The \emph{genus} of $\Gamma$ is the genus of a model and is given by the number $|E(G)| - |V(G)| +1$, since we do \emph{not} allow genus at vertices. 
\end{definition}
For us, the structure preserving maps of $\mathbb{T}\mathcal{C}$ will be \emph{harmonic morphisms} (see e.g. \cite{MR2525845} or \cite{MR3278571}, Section 2.1): Let $\Gamma$ and $\tilde{\Gamma}$ be tropical curves. A continuous and surjective map $\varphi: \Gamma \rightarrow \tilde{\Gamma} $ is called a harmonic morphism if there exist models $(G,l)$ and $(\tilde{G},\tilde{l})$ such that
\begin{itemize}
    \item $\varphi(V(G))\subset V(\tilde{G})$ and $\varphi^{-1}(E(\tilde{G}))\subset E(G)$.
    \item $\varphi$ is locally integer affine linear: On each edge $e\in E(G)$, $\varphi$ restricts to an affine function with integer slope $d_e(\varphi)$ (possibly $0$), called the \emph{weight} or \emph{expansion factor} of $\varphi$ at $e$.
    \item $\varphi$ is harmonic/balanced at every $P\in \Gamma$: For any $\tilde{v}\in T_{\varphi(P)} \tilde{\Gamma}$
    \begin{align}
        d_P(\varphi):=\sum_{v\in T_P\Gamma, v \mapsto \tilde{v}} d_{v}(\varphi)
    \end{align}
    is independent of $\tilde{v}$, where $T_P\Gamma$ ($T_{\varphi(P)} \tilde{\Gamma}$) is the set of tangent directions emanating from $P$ ($\varphi(P)$) and $d_{v}(\varphi)$ is the directional derivative of $\varphi$ in the direction of $v$ (i.e. $d_{v}(\varphi):=d_e(\varphi)$ for the edge $e$ in direction of $v$).
    % d_{v}(\varphi) habe ich nicht extra definiert (auch nicht im ulirsch paper, definition wäre im Lifting harmonic... zu finden)
\end{itemize}
In analogy to the complex setting, we refer to harmonic morphisms as \emph{tropical covers}. A tropical cover $\varphi$ is \emph{finite}, if $d_e(\varphi)>0$ for all edges $e$, and \emph{non-finite} else. Its \emph{degree} is given by the number $deg(\varphi):=\sum_{P\in \Gamma, P \mapsto \tilde{P}} d_{P}(\varphi)$, where $\tilde{P}\in \tilde{\Gamma}$ is an arbitrary point. We will work with tropical curves through choice of a model and by abuse of notation identify $\Gamma$ with $(G,l)$.
\begin{remark}
    The ramification index of a point $P\in \Gamma$ is given by
     \begin{align}
     R_P(\varphi):= 2d_P(\varphi) - 2 - \sum_{v\in T_P\Gamma} (d_{v}(\varphi) -1).
        % R_P(\varphi):= d_P(\varphi)(2-2g(\varphi(P))) - (2-2g(P)) - \sum_{v\in T_P\Gamma} (d_{v}(\varphi) -1).
     \end{align}
     One may then distinguish between \emph{unramified covers}, i.e. covers that satisfy $R_P(\varphi)=0$ for all $P\in \Gamma$, and \emph{ramified} ones. The map $\varphi$ in Figure \ref{figure_examplenotoptimalcover} (b), for example, is ramified, while $\phi$ is not. We mostly consider ramified covers.
     
\end{remark}
\subsection{Optimal tropical covers}\label{subsection_optimalcovers}
As in algebraic geometry (see e.g. \cite{MR0936803}), there is a distinguished class of so-called \emph{optimal covers} among the morphisms of $\mathbb{T}\mathcal{C}$. We introduce them here.

\begin{definition}\label{definition_optimalmap}
Let $\Gamma$ and $\mathbb{T}E$ be curves of genus $2$ and $1$. We call a harmonic map $\varphi: \Gamma \rightarrow \mathbb{T}E$ \emph{optimal} if it does not factor through a non-trivial cover, i.e. if there exists a curve $\mathbb{T}\tilde{E}$ and maps $\tilde{\varphi}: \Gamma \rightarrow \mathbb{T}\tilde{E}$, $\phi: \mathbb{T}\tilde{E} \rightarrow \mathbb{T}E $ such that 
\begin{center}
       \begin{tikzcd}
\Gamma \arrow[r,"\tilde{\varphi}"] \arrow[dr,"\varphi",swap]
& \mathbb{T}\tilde{E} \arrow[d,"\phi"]\\
& \mathbb{T}E
\end{tikzcd}
\end{center}
 
commutes, then $\phi$ is an isomorphism (i.e. $deg(\phi)=1$).
\end{definition}
\begin{figure}[hbt!]
    \centering
\tikzset{every picture/.style={line width=0.75pt}} %set default line width to 0.75pt        

\begin{tikzpicture}[x=0.75pt,y=0.75pt,yscale=-0.82,xscale=0.74]

%Shape: Ellipse [id:dp08479369437105211] 
\draw   (366.3,223.14) .. controls (366.32,216.25) and (397.69,210.68) .. (436.37,210.69) .. controls (475.05,210.7) and (506.39,216.3) .. (506.37,223.19) .. controls (506.35,230.08) and (474.98,235.66) .. (436.3,235.64) .. controls (397.62,235.63) and (366.28,230.03) .. (366.3,223.14) -- cycle ;
%Straight Lines [id:da7175927960258821] 
\draw    (432.33,129.67) -- (432.56,173.58) ;
\draw [shift={(432.57,175.58)}, rotate = 269.7] [color={rgb, 255:red, 0; green, 0; blue, 0 }  ][line width=0.75]    (10.93,-3.29) .. controls (6.95,-1.4) and (3.31,-0.3) .. (0,0) .. controls (3.31,0.3) and (6.95,1.4) .. (10.93,3.29)   ;
%Shape: Arc [id:dp10620459597211007] 
\draw  [draw opacity=0] (358.4,86.31) .. controls (344.51,82.64) and (336.13,77.9) .. (336.13,72.72) .. controls (336.15,61.2) and (377.69,51.91) .. (428.92,51.97) .. controls (480.15,52.04) and (521.67,61.43) .. (521.65,72.94) .. controls (521.65,77.83) and (514.18,82.31) .. (501.68,85.85) -- (428.89,72.83) -- cycle ; \draw   (358.4,86.31) .. controls (344.51,82.64) and (336.13,77.9) .. (336.13,72.72) .. controls (336.15,61.2) and (377.69,51.91) .. (428.92,51.97) .. controls (480.15,52.04) and (521.67,61.43) .. (521.65,72.94) .. controls (521.65,77.83) and (514.18,82.31) .. (501.68,85.85) ;  
%Shape: Ellipse [id:dp4761551420106953] 
\draw   (357.91,86.18) .. controls (357.94,78.22) and (390.36,71.77) .. (430.34,71.79) .. controls (470.32,71.8) and (502.72,78.27) .. (502.69,86.23) .. controls (502.67,94.19) and (470.24,100.63) .. (430.26,100.62) .. controls (390.28,100.6) and (357.89,94.14) .. (357.91,86.18) -- cycle ;
%Shape: Ellipse [id:dp5645337041060137] 
\draw   (495.72,155.04) .. controls (495.75,147.63) and (525,141.63) .. (561.05,141.64) .. controls (597.11,141.66) and (626.31,147.67) .. (626.29,155.08) .. controls (626.26,162.49) and (597.01,168.49) .. (560.96,168.48) .. controls (524.9,168.47) and (495.7,162.45) .. (495.72,155.04) -- cycle ;
%Straight Lines [id:da45594270427107997] 
\draw    (519.3,107.46) -- (538.61,128.53) ;
\draw [shift={(539.96,130)}, rotate = 227.48] [color={rgb, 255:red, 0; green, 0; blue, 0 }  ][line width=0.75]    (10.93,-3.29) .. controls (6.95,-1.4) and (3.31,-0.3) .. (0,0) .. controls (3.31,0.3) and (6.95,1.4) .. (10.93,3.29)   ;
%Straight Lines [id:da04372987177593202] 
\draw    (517.33,183) -- (498.71,202.55) ;
\draw [shift={(497.33,204)}, rotate = 313.6] [color={rgb, 255:red, 0; green, 0; blue, 0 }  ][line width=0.75]    (10.93,-3.29) .. controls (6.95,-1.4) and (3.31,-0.3) .. (0,0) .. controls (3.31,0.3) and (6.95,1.4) .. (10.93,3.29)   ;
%Straight Lines [id:da5063719659940809] 
\draw    (357.91,86.18) ;
\draw [shift={(357.91,86.18)}, rotate = 0] [color={rgb, 255:red, 0; green, 0; blue, 0 }  ][fill={rgb, 255:red, 0; green, 0; blue, 0 }  ][line width=0.75]      (0, 0) circle [x radius= 3.35, y radius= 3.35]   ;
%Straight Lines [id:da4428912648958778] 
\draw    (502.13,85.72) -- (502.69,86.23) ;
\draw [shift={(502.69,86.23)}, rotate = 42.05] [color={rgb, 255:red, 0; green, 0; blue, 0 }  ][fill={rgb, 255:red, 0; green, 0; blue, 0 }  ][line width=0.75]      (0, 0) circle [x radius= 3.35, y radius= 3.35]   ;
%Shape: Ellipse [id:dp07553804582367651] 
\draw   (13.21,223.28) .. controls (13.23,215.18) and (45.28,208.63) .. (84.79,208.64) .. controls (124.3,208.66) and (156.31,215.23) .. (156.29,223.33) .. controls (156.27,231.43) and (124.22,237.98) .. (84.71,237.97) .. controls (45.2,237.96) and (13.18,231.38) .. (13.21,223.28) -- cycle ;
%Straight Lines [id:da21674156651342869] 
\draw    (118.33,136.67) -- (101.29,167.91) ;
\draw [shift={(100.33,169.67)}, rotate = 298.61] [color={rgb, 255:red, 0; green, 0; blue, 0 }  ][line width=0.75]    (10.93,-3.29) .. controls (6.95,-1.4) and (3.31,-0.3) .. (0,0) .. controls (3.31,0.3) and (6.95,1.4) .. (10.93,3.29)   ;
%Shape: Arc [id:dp6689407299667709] 
\draw  [draw opacity=0] (89.27,91.31) .. controls (75.38,87.64) and (67,82.9) .. (67.01,77.72) .. controls (67.02,66.2) and (108.56,56.91) .. (159.79,56.97) .. controls (211.02,57.04) and (252.54,66.43) .. (252.52,77.94) .. controls (252.52,82.83) and (245.05,87.31) .. (232.55,90.85) -- (159.77,77.83) -- cycle ; \draw   (89.27,91.31) .. controls (75.38,87.64) and (67,82.9) .. (67.01,77.72) .. controls (67.02,66.2) and (108.56,56.91) .. (159.79,56.97) .. controls (211.02,57.04) and (252.54,66.43) .. (252.52,77.94) .. controls (252.52,82.83) and (245.05,87.31) .. (232.55,90.85) ;  
%Shape: Ellipse [id:dp9095982856105416] 
\draw   (88.78,91.18) .. controls (88.81,83.22) and (121.23,76.77) .. (161.21,76.79) .. controls (201.19,76.8) and (233.59,83.27) .. (233.56,91.23) .. controls (233.54,99.19) and (201.11,105.63) .. (161.13,105.62) .. controls (121.15,105.6) and (88.76,99.14) .. (88.78,91.18) -- cycle ;
%Shape: Ellipse [id:dp17121665561017285] 
\draw   (174.29,222.06) .. controls (174.31,214.66) and (204.23,208.68) .. (241.11,208.69) .. controls (277.98,208.7) and (307.86,214.71) .. (307.83,222.11) .. controls (307.81,229.51) and (277.89,235.5) .. (241.01,235.48) .. controls (204.13,235.47) and (174.26,229.46) .. (174.29,222.06) -- cycle ;
%Straight Lines [id:da740710675048627] 
\draw    (205.33,137.67) -- (220.39,165.9) ;
\draw [shift={(221.33,167.67)}, rotate = 241.93] [color={rgb, 255:red, 0; green, 0; blue, 0 }  ][line width=0.75]    (10.93,-3.29) .. controls (6.95,-1.4) and (3.31,-0.3) .. (0,0) .. controls (3.31,0.3) and (6.95,1.4) .. (10.93,3.29)   ;
%Straight Lines [id:da5454245143495668] 
\draw    (88.78,91.18) ;
\draw [shift={(88.78,91.18)}, rotate = 0] [color={rgb, 255:red, 0; green, 0; blue, 0 }  ][fill={rgb, 255:red, 0; green, 0; blue, 0 }  ][line width=0.75]      (0, 0) circle [x radius= 3.35, y radius= 3.35]   ;
%Straight Lines [id:da318138791783102] 
\draw    (233,90.72) -- (233.56,91.23) ;
\draw [shift={(233.56,91.23)}, rotate = 42.05] [color={rgb, 255:red, 0; green, 0; blue, 0 }  ][fill={rgb, 255:red, 0; green, 0; blue, 0 }  ][line width=0.75]      (0, 0) circle [x radius= 3.35, y radius= 3.35]   ;

% Text Node
\draw (428.57,243.83) node [anchor=north west][inner sep=0.75pt]    {$\mathbb{T}\tilde{E}$};
% Text Node
\draw (530.88,55.63) node [anchor=north west][inner sep=0.75pt]    {$\Gamma $};
% Text Node
\draw (402.57,136.98) node [anchor=north west][inner sep=0.75pt]    {$\tilde{\varphi }$};
% Text Node
\draw (344.31,40.49) node [anchor=north west][inner sep=0.75pt]    {$8$};
% Text Node
\draw (365.97,99.92) node [anchor=north west][inner sep=0.75pt]    {$4$};
% Text Node
\draw (365.97,59.91) node [anchor=north west][inner sep=0.75pt]    {$4$};
% Text Node
\draw (571.46,178.12) node [anchor=north west][inner sep=0.75pt]    {$\mathbb{T} E$};
% Text Node
\draw (543.43,96.93) node [anchor=north west][inner sep=0.75pt]    {$\textcolor[rgb]{0.49,0.83,0.13}{\varphi }$};
% Text Node
\draw (601.17,126.4) node [anchor=north west][inner sep=0.75pt]    {$\textcolor[rgb]{0.49,0.83,0.13}{4}$};
% Text Node
\draw (522,195.4) node [anchor=north west][inner sep=0.75pt]    {$\textcolor[rgb]{0.49,0.83,0.13}{\phi }$};
% Text Node
\draw (464.8,101.72) node [anchor=north west][inner sep=0.75pt]    {$\textcolor[rgb]{0.49,0.83,0.13}{1}$};
% Text Node
\draw (500.99,43.47) node [anchor=north west][inner sep=0.75pt]    {$\textcolor[rgb]{0.49,0.83,0.13}{2}$};
% Text Node
\draw (462.96,57.27) node [anchor=north west][inner sep=0.75pt]    {$\textcolor[rgb]{0.49,0.83,0.13}{1}$};
% Text Node
\draw (37.75,47.63) node [anchor=north west][inner sep=0.75pt]    {$\Gamma $};
% Text Node
\draw (80.45,133.98) node [anchor=north west][inner sep=0.75pt]    {$\tilde{\varphi }$};
% Text Node
\draw (67.18,50.49) node [anchor=north west][inner sep=0.75pt]    {$1$};
% Text Node
\draw (156.84,85.92) node [anchor=north west][inner sep=0.75pt]    {$1$};
% Text Node
\draw (156.84,59.91) node [anchor=north west][inner sep=0.75pt]    {$1$};
% Text Node
\draw (232.33,248.12) node [anchor=north west][inner sep=0.75pt]    {$\mathbb{T} E$};
% Text Node
\draw (238.3,139.93) node [anchor=north west][inner sep=0.75pt]    {$\varphi $};
% Text Node
\draw (21,191.4) node [anchor=north west][inner sep=0.75pt]    {$12$};
% Text Node
\draw (292,193.4) node [anchor=north west][inner sep=0.75pt]    {$3$};
% Text Node
\draw (70.57,241.83) node [anchor=north west][inner sep=0.75pt]    {$\mathbb{T}\tilde{E}$};
% Text Node
\draw (286,21.4) node [anchor=north west][inner sep=0.75pt]    {$( b)$};
% Text Node
\draw (6,21.4) node [anchor=north west][inner sep=0.75pt]    {$( a)$};

\end{tikzpicture}

    \caption{A curve of genus $2$ covering two elliptic curves. The numbers in Figure \ref{figure_examplenotoptimalcover} (a) are edge lengths and in Figure \ref{figure_examplenotoptimalcover} (b) edge weights.}
    \label{figure_examplenotoptimalcover}
\end{figure}
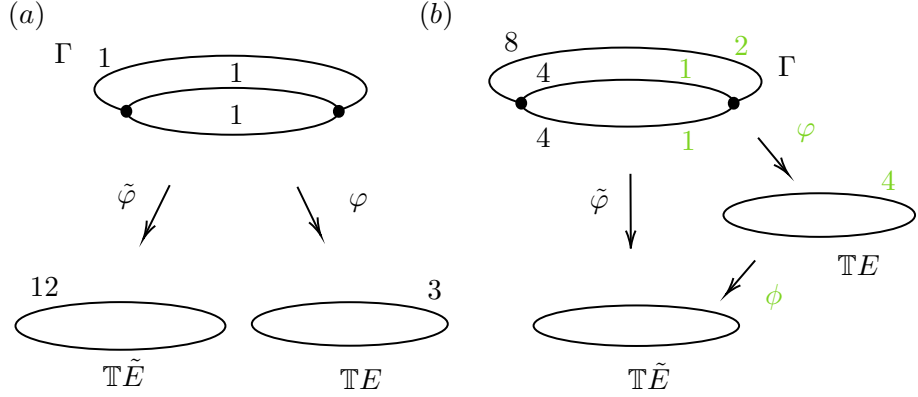
\begin{example}\label{example_optimalandnotoptimal}
Figure \ref{figure_examplenotoptimalcover} (a) shows a cover of degree $2$ on the right, which for degree reasons must be optimal. The cover on the left, however, is not as it factors, for example, through the first, giving rise to a cover of degree $4$. This factorization is shown in Figure \ref{figure_examplenotoptimalcover} (b), with numbers corresponding to edge weights, the green ones to those of $\varphi$ and $\phi$, and the blacks ones to those of $\tilde{\varphi}$.
\end{example}
\section{Crossing Bridges}\label{section_crossingbridges}
\subsection{Preliminaries}
We establish a functorial connection between $\mathbb{T}\mathcal{C}$ and $\mathbb{T}\mathcal{A}$. To an object $\Gamma$ of $\mathbb{T}\mathcal{C}$ we assign an object of $\mathbb{T}\mathcal{A}$ through the following construction (see \cite{MR2772537} or \cite{MR4382460}): 
\begin{construction}\label{construction_theJacobian}
    Let $(G,l)$ be an oriented model of $\Gamma$ and $s,t: E(G) \rightarrow V(G)$ the source and target maps. Then $(G,l)$ comes with two lattices that are related by a non-degenerate pairing:
    \begin{itemize}
        \item \emph{The lattice of harmonic} $1$-\emph{forms}, $\Omega^1_G(\mathbb{Z})$: For each oriented edge $e$ we introduce a formal symbol $de$ called a \emph{basic} $1$-form on $G$ and set $\Omega^1_G(\mathbb{Z})$ to be 
        \begin{align}
            \{ \omega:=\sum_{e} \omega_e de: \omega_e \in \mathbb{Z}, \thinspace \sum_{e:t(e)=V} \omega_e = \sum_{e:s(e)=V} \omega_e \thinspace \forall \thinspace V\in V(G) \}.
        \end{align} 
        It is the subgroup of the free group over $\{ de: e\in E(G)\} $ consisting of harmonic $1-$forms on $G$.
        \item \emph{The lattice of integral} $1$-\emph{cycles}, $\H1(G,\mathbb{Z})$: It is the first simplicial homology group of $G$ given by $\ker(\partial)$, where 
        \begin{align}
            \partial: \C1(G,\mathbb{Z})\rightarrow \Co (G,\mathbb{Z}), e\mapsto t(e)-s(e)
        \end{align} is the boundary operator.
        \item \emph{The integration pairing}, $\int_\cdot \cdot$: We can integrate a basic $1$-form
    \begin{align}
        \int_e de':=\begin{cases}
            l(e), \text{ if } e=e'\\
            0, \text{ else}
        \end{cases}
    \end{align}
    and extend linearly to obtain a perfect pairing
    \begin{align}
       \int_\cdot \cdot: \Omega^1_G(\mathbb{Z}) \times \H1(G,\mathbb{Z}) \rightarrow \mathbb{R}, (\omega, c) \mapsto \int_c \omega.
    \end{align}
    \end{itemize}
    These building blocks are independent of the choice of model (see also \cite{MR2772537}). This means that lattices that arise from different models (that have compatible orientations) are related by isomorphisms, that leave the integration pairing invariant. We will write $\Omega^1_\Gamma(\mathbb{Z})$ and $\H1(\Gamma,\mathbb{Z})$, instead, and complete Construction \ref{construction_theJacobian} by assigning a pptav to $\Gamma$.
\begin{definition}\label{definition_JacobianandAbelJacobimap}
    The \emph{Jacobian of} $\Gamma$ is the pptav built from $(\Omega^1_\Gamma(\mathbb{Z}),\H1(\Gamma,\mathbb{Z}),\int_{\cdot} \cdot)$ with principal polarization $\gls{zetaGamma}: \H1(\Gamma,\mathbb{Z}) \rightarrow \Omega^1_\Gamma(\mathbb{Z}), \sum a_e e \mapsto \sum a_e de $. It is related to $\Gamma$ by the \emph{tropical Abel-Jacobi map}:
\begin{align}
    \gls{PhiP0}: \Gamma \rightarrow \Jac(\Gamma), P \mapsto \int_{\gamma_{P}}\cdot,
\end{align}
where $P_0\in V(\Gamma)$ is a fixed vertex and $\gamma_P\in \C1(\Gamma,\mathbb{Z})$ is any path connecting $P_0$ to $P$ in $\Gamma$.
\end{definition}
\end{construction}
\begin{remark}\label{remark_JacobianandPicardvariety}
Note that the construction of the tropical Jacobian works on exactly the same principle as the classical one over $\mathbb{C}$. In fact, the analogy goes one step further: We have an isomorphism $\Jac(\Gamma)\overset{\Phi_\Gamma}{\cong} \Pic^0(\Gamma)$ (\cite{MR2457739}, Theorem 6.2), where $\Pic^0(\Gamma):=\Div^0(\Gamma)/\Prin(\Gamma)$ is the group-theoretic quotient of $\Div^0(\Gamma)$, the group of divisors of degree $0$, by $\Prin(\Gamma)$, the subgroup of principal divisors (\cite{MR2457739}, Section 4.2). This setting makes it easy to specify two homomorphsims that are induced by a cover $\varphi: \Gamma \rightarrow \tilde{\Gamma}$: The push-forward $\varphi_*: \Pic^0(\Gamma) \rightarrow \Pic^0(\tilde{\Gamma})$ and the pull-back $\varphi^*: \Pic^0(\tilde{\Gamma}) \rightarrow \Pic^0(\Gamma)$ of divisors defined by
\begin{align}
    \varphi_* \Big{(}\sum_i n_iP_i\Big{)}:=\sum_i n_i \varphi(P_i) \text{ and } \varphi^*\Big{(}\sum_j \tilde{n}_j \tilde{P}_j\Big{)}:=\sum_j\sum_{ P_{ij} \mapsto \tilde{P}_j} \tilde{n}_j d_{P_{ij}}(\varphi) P_{ij}.
\end{align}
\end{remark}
Now, Remark \ref{remark_JacobianandPicardvariety} suggests two natural candidates for playing the counterpart of a cover $\varphi$ in $\mathbb{T}\mathcal{A}$. We start by describing these as morphisms of tavs and understanding their relationship: As the notation suggests, $\varphi_*$ and $\varphi^*$ are dual to each other.
\begin{lemma}\label{lemma_pushandpullaredual}
Let $\varphi: \Gamma \rightarrow \tilde{\Gamma} $ be a cover. We consider the morphism $\psi_{*}: \Jac(\Gamma) \rightarrow \Jac(\tilde{\Gamma})$ induced by the push-forward $\varphi_*$ under the identification  $\Pic^0(\Gamma) \rightarrow \Jac(\Gamma)$, respectively $\Pic^0(\tilde{\Gamma}) \rightarrow \Jac(\tilde{\Gamma})$, given by the tropical Abel-Jacobi map. Then the dual homomorphism $\widecheck{\psi}_{*}:  \widecheck{\Jac}(\tilde{\Gamma}) \rightarrow \widecheck{\Jac}(\Gamma)$ is induced by the pull-back on divisors $\varphi^*: \Pic^0(\tilde{\Gamma}) \rightarrow \Pic^0(\Gamma)$. 
\end{lemma}
Lemma \ref{lemma_pushandpullaredual} is most certainly known, although a proof does not seem to appear in the literature. We include it here for the sake of completeness and as an illustration of Remark \ref{remark_JacobianandPicardvariety}.
\begin{proof} 
 Let us fix reference points $Q_\Gamma\in \Gamma$ and $\varphi(Q_\Gamma)\in \tilde{\Gamma}$ and consider the map $\psi_{*}$ that makes the following diagram commute:\\
 \begin{center}
     \begin{tikzcd}
 \Pic^0(\Gamma) \arrow[r, "\varphi_* "] \arrow[d, "\phi_\Gamma"]
& \Pic^0(\tilde{\Gamma}) \arrow[d, "\phi_{\tilde{\Gamma}}"] \\
\Jac(\Gamma) \arrow[r,"\psi_*"] & \Jac(\tilde{\Gamma})\\
\end{tikzcd}.
 \end{center}
   In doing so, we may figure out the assignment rule of $\psi_{*}$ by "diagram-chasing" elements: As the vertical arrows are the isomorphisms induced by the Abel-Jacobi map (\cite{MR2457739}, Theorem 6.2), i.e. defined by
 \begin{align}
 \phi_{\Gamma}\Big{(}\sum_i n_iP_i\Big{)}=\sum_i n_i \int^{P_i}_{Q_\Gamma}(\cdot) \text{ and } \phi_{\tilde{\Gamma}}\Big{(}\sum_i n_i\tilde{P}_i\Big{)}=\sum_i n_i \int^{\tilde{P}_i}_{\varphi(Q_\Gamma)}(\cdot),
 \end{align}
 we get
  \begin{align}
     \psi_{*}: \Jac(\Gamma) \rightarrow \Jac(\tilde{\Gamma}): \sum_i n_i \int^{P_i}_{Q_\Gamma} (\cdot) \mapsto \sum_i n_i \int^{\varphi(P_i)}_{\varphi(Q_\Gamma)} (\cdot).
  \end{align}
  We claim that $\psi_{*}$ is given by the pair $(\psi^\#_{*},\psi_{*\#})$, where $\psi^\#_{*}$ is the pull-back of $1$-forms and $\psi_{*\#}$ the push-forward of cycles, i.e. that for all $\sum_i n_i \int^{P_i}_{Q_\Gamma} (\cdot)\in \Jac(\Gamma)$
  \begin{align}
      \Hom (\psi^\#_{*},\mathbb{Z})\Big{(} \sum_i n_i \int^{P_i}_{Q_\Gamma} (\cdot) \Big{)}=\sum_i n_i \int^{P_i}_{Q_\Gamma}  \psi^\#_{*} (\cdot) 
      \text{ and } \sum_i n_i \int^{\varphi(P_i)}_{\varphi(Q_\Gamma)}(\cdot) 
  \end{align}
  agree on $\Omega^1_{\tilde{\Gamma}}(\mathbb{Z})$. Inserting an arbitrary $\omega=\sum_k a_k d \tilde{e}_k \in \Omega^1_{\tilde{\Gamma}}(\mathbb{Z})$ and using linearity of the integral we see that we only need to verify the equality by summands, i.e. that for fixed tuple $(k,i)$ we have
  \begin{align}
      \int^{\varphi(P)}_{\varphi(Q_\Gamma)}(d\tilde{e}_k)=\sum_{e_j \mapsto \tilde{e}_k}\int^{P}_{Q_\Gamma}d_{e_j}(\varphi) de_j, \text{ where } P:=P_i.
  \end{align}
 Choose a model of $\Gamma$ ($\tilde{\Gamma}$) containing $Q_\Gamma$ and $P$ ($\varphi(Q_\Gamma)$ and $\varphi(P)$) as vertices and such that the set $E$ ($\tilde{E}$) of edges, that cannot be subdivided further, satisfies the following conditions: No $e\in E$ ($\tilde{E}$) is a cycle, nor does $\varphi(e)$ contain any cycle. If $\gamma:=\sum_j m_j e_j$ is a path (with $e_j \in E$) in $\Gamma$ from $Q_\Gamma$ to $P$, then $\tilde{\gamma}:=\sum_j m_j \varphi(e_j)$ is a path connecting $\varphi(Q_\Gamma)$ to $\varphi(P)$ in $\tilde{\Gamma}$. Since integration is path-independent in the respective Jacobians, we use $\gamma$ and $\tilde{\gamma}$ to compute the integrals above:
 \begin{align}
    &\int^{\varphi(P)}_{\varphi(Q_\Gamma)}(d\tilde{e}_k)= \Big{(}\sum_{e_j \mapsto \tilde{e}_k} m_j\Big{)}\cdot l(\tilde{e}_k) \text{ and}
    \sum_{e_j \mapsto \tilde{e}_k}\int^{P}_{Q_\Gamma}d_{e_j}(\varphi) de_j= \sum_{e_j \mapsto \tilde{e}_k} m_j \thinspace d_{e_j}(\varphi) \thinspace l(e_j).
 \end{align}
 Recalling how $\varphi$ behaves with respect to edge lengths (i.e. $d_{e_j}(\varphi)l(e_j)=l(\tilde{e}_k)$ holds since by our choice of models we excluded winding) finally proves the claim.\\
 The transposed pair $(\psi_{*\#},\psi^\#_{*},)$ is the dual of $\psi_*$. Under the natural isomorphisms $\Omega^1_\Gamma \overset{\gls{zetaGamma}}{\cong} \H1(\Gamma,\mathbb{Z})$ and $\Omega^1_{\tilde{\Gamma}} \cong \H1(\tilde{\Gamma},\mathbb{Z})$ we can interpret it as morphism between the Jacobians we denote by $\psi^{*}=(\psi^{*\#}, \psi^*_{\#} )$, where $\psi^{*\#}$ is the push-forward of $1$-forms and $\psi^*_{\#}$ the pull-back of cycles. As above one shows that ${\psi}^{*}$ is induced by the pull-back of divisors $\varphi^*: \Pic^0(\tilde{\Gamma}) \rightarrow \Pic^0(\Gamma)$, thereby completing the proof. 
\end{proof}
From now on the identification $\Jac(\Gamma)\cong \Pic^0(\Gamma)$ will be implicit and, justified by Lemma \ref{lemma_pushandpullaredual}, we use the notation $\varphi^*$ and $\varphi_*$ for both settings.
\subsection{Curves of genus 2 covering curves of genus 1.}\label{subsection_curvesofg2coverg1}
We explore this connection in the setting of curves of genus $2$ covering curves of genus $1$: As Definition \ref{definition_optimalmap} suggests we will be interested in factorizing tropical coverings. A daunting task made easier by the following Lemma that allows us to transfer the factorization question from the category of tropical curves to the category of tropical abelian varieties.

Note that the tropical Abel-Jacobi map is an isomorphism in the case of genus $1$. It can be used to endow $\mathbb{T}E$ with the structure of a group. %, i.e. by choosing a reference point we can turn $\mathbb{T}E$ into pptav.
In analogy to the classical case we call $\gls{TE}$ an \emph{elliptic curve} and depending on the context, regard $\mathbb{T}E$ as an object of $\mathbb{T}\mathcal{C}$ or, by abuse of notation, as an object of $\mathbb{T}\mathcal{A}$ (the identification of $\Jac(\mathbb{T}E)$ with $\mathbb{T}E$ will be used implicitly).

\begin{lemma}\label{lemma_facorizationgtropicalmaps}
   Let $\varphi: \Gamma \rightarrow \mathbb{T}E$ and $\tilde{\varphi}: \Gamma \rightarrow \mathbb{T}\tilde{E}$ be harmonic maps and $\varphi_*: \Jac(\Gamma) \rightarrow \mathbb{T}E$ and $\tilde{\varphi}_*: \Jac(\Gamma) \rightarrow \mathbb{T}\tilde{E}$ the respective push-forwards. Then the following are equivalent:
    
    \begin{itemize}
       \item  $\varphi=\phi \circ \tilde{\varphi}$ for an isogeny $\phi:\mathbb{T}\tilde{E} \rightarrow \mathbb{T}E $.
      \item $\varphi_*=\hat{\phi} \circ \tilde{\varphi}_*$ for an isogeny $\hat{\phi}:\mathbb{T}\tilde{E} \rightarrow \mathbb{T}E $.
    \end{itemize}
    
\end{lemma}
\begin{remark}\label{remark_muispushforward}
The proof of Lemma \ref{lemma_facorizationgtropicalmaps} relies on a universal property of the Jacobian (\cite{röhrle2024tropicalngonalconstruction}, Proposition 4.14). Applied to our setting, the map $\mu: \Jac(\Gamma) \rightarrow \mathbb{T}E$ in \cite{röhrle2024tropicalngonalconstruction} induced by $\varphi: \Gamma \rightarrow \mathbb{T}E$ corresponds to $\varphi_*$. To make the connection to \cite{röhrle2024tropicalngonalconstruction} clearer we write $\mu$ for $\varphi_*$ and $\tilde{\mu}$ for $\tilde{\varphi}_*$ in the proof of Lemma \ref{lemma_facorizationgtropicalmaps}.
\end{remark}
 
\begin{proof}[Proof of Lemma \ref{lemma_facorizationgtropicalmaps}.]
 The fact that factorization on the level of tropical covers corresponds to factorization on the level of Jacobians is essentially a consequence of Proposition 4.14. in \cite{röhrle2024tropicalngonalconstruction}.\\
   To see this, note that $\varphi$ and $\tilde{\varphi}$ are morphisms of rational polyhedral spaces whose targets are also integral tori. Hence, there exists unique homomorphisms $\mu$ and $\tilde{\mu}$ such that the following diagram commutes for all $P_0\in \Gamma$: \\
   \begin{center}
   \begin{tikzpicture}
\matrix(m)[matrix of math nodes,
row sep=3em, column sep=2.5em,
text height=1.5ex, text depth=0.25ex]
{\mathbb{T}E& \Gamma & \mathbb{T}\tilde{E}\\
\mathbb{T}E&\Jac(\Gamma)&\mathbb{T}\tilde{E}\\};
\path[->,font=\scriptsize]
(m-1-2) edge node[auto] {$\Phi_{P_0}$} (m-2-2);
\path[->,font=\scriptsize]
(m-1-2) edge node[above] {$\tilde{\varphi}$} (m-1-3)
(m-1-3) edge node[auto] {$t_{-\tilde{\varphi}(P_0)}$} (m-2-3)
(m-1-1) edge node[left] {$t_{-\varphi(P_0)}$} (m-2-1)
(m-2-2) edge node[below] {$\tilde{\mu}$} (m-2-3)
;
\path[->,font=\scriptsize]
(m-2-2) edge node[auto] {$\mu$} (m-2-1)
(m-1-2) edge node[above] {$\varphi$} (m-1-1)
;
\end{tikzpicture}.
      
   \end{center}

   If $\varphi=\phi \circ \tilde{\varphi}$, choose $P_0 \in \tilde{\varphi}^{-1}(0)$ and observe that $t_{-\varphi(P_0)}=t_{-\tilde{\varphi}(P_0)}=id$. Then $\hat{\phi}:=\phi$ satisfies $\mu=\hat{\phi} \circ \tilde{\mu}$. Conversely, suppose  $\mu$ factors through an isogeny $\hat{\phi}$. Let $P_0 \in \varphi^{-1}(0)$ and obtain a factorisation of $\varphi$ through an isogeny  $\phi:=\hat{\phi}$ and a harmonic map given by post-composing  $\tilde{\varphi}$ with the translation $t_{-\tilde{\varphi}(P_0)}$:
   \begin{center}

   \begin{tikzpicture}
\matrix(m)[matrix of math nodes,
row sep=3em, column sep=2.8em,
text height=1.5ex, text depth=0.25ex]
{\Gamma& \mathbb{T}\tilde{E} & \mathbb{T}\tilde{E} & \mathbb{T}E.\\};
\path[->]
(m-1-1) edge node[above] {$\tilde{\varphi}$} (m-1-2)
(m-1-2) edge node[above] {$t_{-\tilde{\varphi}(P_0)}$} (m-1-3)
(m-1-3) edge node[above] {$\hat{\phi}$} (m-1-4)
(m-1-1) edge [bend right] node[below] {$\varphi$} (m-1-4)
;
\end{tikzpicture}
      
   \end{center}
\end{proof}
Lemma \ref{lemma_facorizationgtropicalmaps} defines our approach: Translating covering language into the language of pptav. This gives us the opportunity to address cover-related questions with a new set of tools and guides our next steps: We study 
\begin{itemize}
    \item isogenies of elliptic curves
    \item the push-forward and pull-back morphism
    \end{itemize}
separately, and merge our results in Subsection \ref{subsection_CriteriaforOptimality} to develop criteria for verifying optimality.
\subsection{Isogenies of Elliptic Curves.}
We start by retrieving (the analogy of) an algebraic result in the tropical world: Recall that multiplication-by-$a$ maps, where $a$ is a complex number, yield isogenies between complex elliptic curves (i.e. $1$-dimensional complex abelian varieties) and what is more, all possibilities are fully exhausted by this class of mappings. Essentially the same holds for $1$-dimensional integral tori, that is tropical elliptic curves.
\begin{lemma}\label{lemma_isogeniesbetweenellipticcurves}
    Let $\phi$ be an isogeny between elliptic curves, $\mathbb{T}E$ and $\mathbb{T}\tilde{E}$, of length $l$ and $\tilde{l}$.\\
    Under the identifications $\mathbb{T}E\cong \mathbb{R}/l\mathbb{Z}$ and $\mathbb{T}\tilde{E}\cong \mathbb{R}/\tilde{l}\mathbb{Z}$ we have: $\phi$ is induced by an $\mathbb{R}-$linear map $\tilde{\phi}: \mathbb{R} \rightarrow \mathbb{R}, z \mapsto a^\#\cdot z$, where $a^\#\in \mathbb{Z}$ and $\tilde{\phi}(\tilde{l}\mathbb{Z})\subset l\mathbb{Z}$ holds.\\
    Under the identifications $\mathbb{T}E\cong \Jac(\mathbb{T}E)$ and $\mathbb{T}\tilde{E}\cong \Jac(\mathbb{T}\tilde{E})$ we have: $\phi=(\phi^\#,\phi_\#)$ is defined by a tuple $(a^\#,a_\#)\in \mathbb{Z}^2$ such that $\tilde{l}a^\#=la_\#$ holds and $\phi$ is an isomorphism precisely when $(a^\#,a_\#)=(\pm 1,\pm 1)$.
\end{lemma}
We see: On the quotient, $\phi$ corresponds to a multiplication-by-$a^\#$ map, where $a^\#\in \mathbb{Z}$. The tuple $(\phi^\#,\phi_\#)$, of course, provides two multiplication maps, that characterize the winding and dilation behaviour of $\phi$ separately (see \cite{röhrle2024tropicalngonalconstruction}). %siehe geometric und dilation degree
\begin{proof}
    We choose to work in the category of tropical abelian varieties. Thus, we will identify all elliptic curves involved with their respective Jacobians. An isogeny $\phi$ is then given by a pair of group homomorphisms 
    \begin{align}
      & \phi^\#:  \Omega^{1}_{\mathbb{T}E} \rightarrow \Omega^{1}_{\mathbb{T}\tilde{E}}: \omega \mapsto a^\# \tilde{\omega}\\
      & \phi_\# : \H1(\mathbb{T}\tilde{E},\mathbb{Z}) \rightarrow  \H1(\mathbb{T}E,\mathbb{Z}): \tilde{B} \mapsto a_\# B
    \end{align}
    where $B$ ($\tilde{B}$) is a homology basis of  $\mathbb{T}E$ ($\mathbb{T}\tilde{E}$) and $\omega$ ($\tilde{\omega}$) the associated basis of tropical 1-forms. The compatibility of $(\phi^\#,\phi_\#)$ with the integration pairings on $\mathbb{T}E$ and $\mathbb{T}\tilde{E}$ poses additional requirements on the pair of integers $(a^\#,a_\#)$: 
    \begin{align}
        \int_{\tilde{B}}\phi^\#(\omega)=\int_{\phi_\#(\tilde{B})}\omega.
    \end{align}
    Specifying the metric on  $\mathbb{T}E$ and $\mathbb{T}\tilde{E}$ by the real numbers $l$ and $\tilde{l}$ finally leads to relation, $\tilde{l}a^\#=la_\#$, as claimed. The actual isogeny $\phi$ is given by $z \mapsto a^\# \cdot z$.
\end{proof}
\begin{remark}
   In algebraic geometry elliptic curves (over algebraically closed fields) are characterized by their $j$-invariant. The appropriate tropical analogue is the cycle length (see \cite{MR2457725}, \cite{MR2275625}). It is a geometric property that allows us to distinguish between non-isomorphic curves. Computations in Lemma \ref{lemma_kernelofpushforwardinTA}, \ref{lemma_kernelpushforwardinTA_DB} and \ref{lemma_explicitdescriptionofgamma} should be seen in this light.
\end{remark}
Continuing in the vein of Lemma \ref{lemma_isogeniesbetweenellipticcurves}, we obtain:
\begin{lemma}\label{lemma_factorizingisogeniesofellipticcurves}
    For $i=1,2,3$ let $\mathbb{T}E_i$ be an elliptic curve of length $l_i$. Let $\phi_1: \mathbb{T}E_1 \rightarrow \mathbb{T}E_2$ and  $\phi_2: \mathbb{T}E_1 \rightarrow \mathbb{T}E_3$ be two isogenies given by tuples $(a_i^\#,a_{i \#})\in \mathbb{Z}^2$, $i=1,2$ as in Lemma \ref{lemma_isogeniesbetweenellipticcurves}. Then $\phi_1$ factors through $\phi_2$ if and only the following holds:
    \begin{align}
        a_3^\#:= \frac{a_1^\#}{a_2^\#}\in \mathbb{Z} \text{ and } a_{3 \#}:=a_3^\# \cdot \frac{l_3}{l_2} \in \mathbb{Z}.
    \end{align}
    In this case, $(a_3^\#,a_{3, \#})$ defines an isogeny $\phi_3$ that satisfies $\phi_1=\phi_3 \circ \phi_2$.
\end{lemma}
\begin{proof}
Suppose $\phi_1$ factors through $\phi_2$ via an isogeny $\phi_3$, i.e. $\phi_3$ satisfies $\phi_{3 \#} \circ \phi_{2 \#}=\phi_{1 \#}$ and $\phi_2^ \#\circ \phi_3^ \#=\phi_1^ \#$. Expressed in terms of the defining tuples $(a_i^\#,a_{i \#})$ for $i=1,2,3$ as in Lemma \ref{lemma_isogeniesbetweenellipticcurves} we have, $a_{3 \#} \cdot a_{2 \#}=a_{1 \#}$ and $a_2^\# \cdot a_3^\#=a_1^\#$, or equivalently, using $l_1 a_1^\#=l_2 a_{1\#}$ and $l_1 a_2^\#=l_3 a_{2\#}$:
\begin{diagram}\label{diagram_infactorizingisogeniesofellipticcurves}
    a_3^\#= \frac{a_1^\#}{a_2^\#}\in \mathbb{Z} \text{ and } a_{3 \#}=a_3^\# \cdot \frac{l_3}{l_2} \in \mathbb{Z}.
\end{diagram}
    Conversely, suppose (\ref{diagram_infactorizingisogeniesofellipticcurves}) holds. Then $(a_3^\#,a_{3, \#})$ gives rise to a well-defined morphism of tori $\phi_3$, which satisfies $\phi_1=\phi_3 \circ \phi_2$ by construction.
\end{proof}
\subsection{Procedure to address a factorization problem.}
Here, we describe a procedure to solve the factorization problem described below. How to perform specific computations mentioned in step 2 is covered in Subsection \ref{subsection_push-forward}. Let $\Gamma$ be of genus 2 and suppose $\Gamma$ covers elliptic curves $\mathbb{T}E$ and $\mathbb{T}\tilde{E}$ via $\varphi$, $\tilde{\varphi}$ respectively (see Lemma \ref{lemma_facorizationgtropicalmaps}):
\begin{align}
    \text{When does the morphism } \varphi_* \text{ factor through } \tilde{\varphi}_* ?
\end{align}
\paragraph{\emph{Step 1: Reduction to a lower dimensional factorization problem.}} This factorization problem can be simplified further: Let $\mu:=\varphi_*$ and $\tilde{\mu}:=\tilde{\varphi}_*$ and suppose $\ker(\tilde{\mu}) \subset \ker(\mu)$. Otherwise conclude immediately: $\mu$ does not factor through $\tilde{\mu}$. Let $\mu=\phi_{\mu} \circ \pi_{\mu}$ and  $\tilde{\mu}=\phi_{\tilde{\mu}} \circ \pi_{\tilde{\mu}}$ be the Stein factorization of $\mu$ and $\tilde{\mu}$. Since 
\begin{align}
     \Ker(\mu)_0 + \pi^{-1}_{\mu} (\ker(\phi_\mu))=\ker(\mu) \supset \ker(\tilde{\mu})= \Ker(\tilde{\mu})_0 + \pi^{-1}_{\tilde{\mu}}(\ker(\phi_{\tilde{\mu}}))
\end{align}
we have $\Ker(\tilde{\mu})_0 \subset \Ker(\mu)_0 $. This yields a well-defined map $pr$ such that the left part of
 \begin{diagram}
     %\begin{tikzcd}[row sep=huge] \ar[rr, bend left,"\mu"]
    \Jac(\Gamma)  \ar[dr,"{\pi_{\tilde{\mu}}}",swap]  \ar[r,"\pi_\mu"] & \Jac(\Gamma)/\Ker(\mu)_0 \ar[r,"{\phi_\mu}" ] & \mathbb{T}E \\
    & \Jac(\Gamma)/\Ker(\tilde{\mu})_0 \ar[r,"{\phi_{\tilde{\mu}}}",swap ] \ar[u,," pr"] & \mathbb{T}\tilde{E}  \ar[u,dashed," ?"]\\
   % \end{tikzcd}
\end{diagram}
commutes. As $\Ker(\tilde{\mu})_0$ is a subgroup of finite index of $\Ker(\mu)_0$ since $\Jac(\Gamma)/\Ker(\tilde{\mu})_0 $ and $\Jac(\Gamma)/\Ker(\mu)_0 $ are both 1-dimensional tori, $pr$ is an isogeny (and thus $\phi_\mu \circ pr$ as well). Moreover, $\ker(\phi_\mu \circ pr)\supset \ker(\phi_{\tilde{\mu}})$ holds. That is, the initial factorization problem reduces to the lower dimensional problem of factorizing $\phi_\mu \circ pr$ through $\phi_{\tilde{\mu}}$.
\paragraph{\emph{Step 2: Solve lower dimensional factorization problem.}}
We use Lemma \ref{lemma_factorizingisogeniesofellipticcurves}: Set $\mathbb{T}E_1:=\Jac(\Gamma)/\Ker(\tilde{\mu})_0 $, $\mathbb{T}E_2:=\mathbb{T}E$, $\mathbb{T}E_3:=\mathbb{T}\tilde{E}$, $\phi_1:=\phi_\mu \circ pr$ and $\phi_2:=\phi_{\tilde{\mu}}$. If $\Gamma$ is of type "theta", see Figure \ref{figure_generalgenus2cover}, compute $l_1$ and $(a_i^\#,a_{i \#})$ (for $i=1,2$) using Lemma \ref{lemma_explicitdescriptionofgamma}. For the other maximal combinatorial type of a curve of genus 2, the dumbbell graph (see Figure \ref{figure_dumbbellcover}), we have a formula for $l_1$ in  Lemma \ref{lemma_kernelpushforwardinTA_DB}. The pair $(a_i^\#,a_{i \#})$ is then computed as in the proof of Lemma  \ref{lemma_explicitdescriptionofgamma}.
Then $\mu$ factors through $\tilde{\mu}$ if and only if the following holds:
\begin{itemize}
    \item $\ker(\tilde{\mu}) \subset \ker(\mu)$,
    \item $a_3^\#:= \frac{a_1^\#}{a_2^\#}\in \mathbb{Z} \text{ and } a_{3 \#}:=a_3^\# \cdot \frac{l_3}{l_2} \in \mathbb{Z}.$
\end{itemize}
In this case, $\hat{\phi}:=\phi_3$ solves the initial factorization problem.
\subsection{The push-forward}\label{subsection_push-forward}
As before, let $\varphi: \Gamma \rightarrow \mathbb{T}E$ be a cover and $g(\Gamma)=2$. Then $\varphi_*$ is clearly surjective making its kernel the first object of interest. Note, however, that $\varphi_*$ belongs to two different worlds, the category of tav ($\mathbb{T}\mathcal{A}$) and the category of abelian groups ($Ab$). This means that its kernel will have two different identities as well.

We study both, first in detail for the case where $\Gamma$ is of type "theta" as in Figure \ref{figure_generalgenus2cover}. For the case where the combinatorial type of $\Gamma$ is the dumbbell-graph, we provide a collection of results in Subsection \ref{subsubsection_typedb}. 
\subsubsection{Type theta}\label{subsubsection_typetheta}
\begin{convention}\label{convention_genus2coversofE}
This next part calls for consistency in naming that allows us to refer to individual parts of a cover $\varphi: \Gamma \rightarrow \mathbb{T}E$ explicitly. We introduce the labeling we want to use for the graphs underlying $\Gamma$ and $\mathbb{T}E$ in Figure \ref{figure_graphlabelling}.
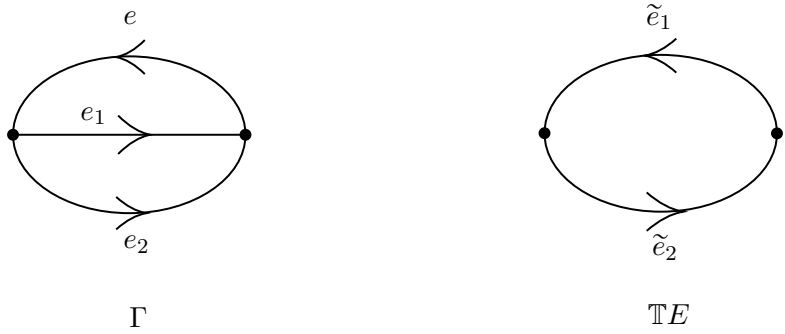
\begin{figure}[H]
    \centering
 \tikzset{every picture/.style={line width=0.75pt}} %set default line width to 0.75pt        

\begin{tikzpicture}[x=0.7pt,y=0.7pt,yscale=-0.8,xscale=0.8]

%Shape: Ellipse [id:dp34235573607566194] 
\draw   (50.33,112) .. controls (50.33,82.73) and (85.48,59) .. (128.83,59) .. controls (172.19,59) and (207.33,82.73) .. (207.33,112) .. controls (207.33,141.27) and (172.19,165) .. (128.83,165) .. controls (85.48,165) and (50.33,141.27) .. (50.33,112) -- cycle ;
%Shape: Ellipse [id:dp7462550369305541] 
\draw   (409.33,110.93) .. controls (409.36,81.66) and (444.53,57.96) .. (487.88,58) .. controls (531.23,58.04) and (566.36,81.8) .. (566.33,111.07) .. controls (566.31,140.34) and (531.14,164.04) .. (487.79,164) .. controls (444.43,163.96) and (409.31,140.2) .. (409.33,110.93) -- cycle ;
%Straight Lines [id:da584019106667631] 
\draw    (50.33,112) -- (207.33,112) ;
%Straight Lines [id:da6966224157773571] 
\draw    (50.33,112) ;
\draw [shift={(50.33,112)}, rotate = 0] [color={rgb, 255:red, 0; green, 0; blue, 0 }  ][fill={rgb, 255:red, 0; green, 0; blue, 0 }  ][line width=0.75]      (0, 0) circle [x radius= 3.35, y radius= 3.35]   ;
%Straight Lines [id:da3432115084468623] 
\draw    (207.33,112) ;
\draw [shift={(207.33,112)}, rotate = 0] [color={rgb, 255:red, 0; green, 0; blue, 0 }  ][fill={rgb, 255:red, 0; green, 0; blue, 0 }  ][line width=0.75]      (0, 0) circle [x radius= 3.35, y radius= 3.35]   ;
%Straight Lines [id:da49659900159263204] 
\draw    (566.33,111.07) ;
\draw [shift={(566.33,111.07)}, rotate = 0] [color={rgb, 255:red, 0; green, 0; blue, 0 }  ][fill={rgb, 255:red, 0; green, 0; blue, 0 }  ][line width=0.75]      (0, 0) circle [x radius= 3.35, y radius= 3.35]   ;
%Straight Lines [id:da8803718960670524] 
\draw    (409.33,110.93) ;
\draw [shift={(409.33,110.93)}, rotate = 0] [color={rgb, 255:red, 0; green, 0; blue, 0 }  ][fill={rgb, 255:red, 0; green, 0; blue, 0 }  ][line width=0.75]      (0, 0) circle [x radius= 3.35, y radius= 3.35]   ;
\draw   (139.33,48) .. controls (132.89,54.39) and (126.44,58.22) .. (120,59.5) .. controls (126.44,60.78) and (132.89,64.61) .. (139.33,71) ;
\draw   (120,154) .. controls (126.78,160.11) and (133.56,163.78) .. (140.33,165) .. controls (133.56,166.22) and (126.78,169.89) .. (120,176) ;
\draw   (478.33,151) .. controls (486.56,158.22) and (494.78,162.56) .. (503,164) .. controls (494.78,165.44) and (486.56,169.78) .. (478.33,177) ;
\draw   (121.33,99) .. controls (128.56,106.22) and (135.78,110.56) .. (143,112) .. controls (135.78,113.44) and (128.56,117.78) .. (121.33,125) ;
\draw   (498.33,47) .. controls (491.22,53.39) and (484.11,57.22) .. (477,58.5) .. controls (484.11,59.78) and (491.22,63.61) .. (498.33,70) ;

% Text Node
\draw (123,26.4) node [anchor=north west][inner sep=0.75pt]    {$e$};
% Text Node
\draw (94,91.4) node [anchor=north west][inner sep=0.75pt]    {$e_{1}$};
% Text Node
\draw (123,178.4) node [anchor=north west][inner sep=0.75pt]    {$e_{2}$};
% Text Node
\draw (476,21.4) node [anchor=north west][inner sep=0.75pt]    {$\widetilde{e}_{1}$};
% Text Node
\draw (480,177.4) node [anchor=north west][inner sep=0.75pt]    {$\widetilde{e}_{2}$};
% Text Node
\draw (127,226.4) node [anchor=north west][inner sep=0.75pt]    {$\Gamma $};
% Text Node
\draw (477,224.4) node [anchor=north west][inner sep=0.75pt]    {$\mathbb{T} E$};

\end{tikzpicture}

    \caption{Labeling for the graphs underlying $\Gamma$ and $\mathbb{T}E$. }
    \label{figure_graphlabelling}
\end{figure}
Next, let us fix concrete homology basis (with edges oriented as in Figure \ref{figure_graphlabelling}):
\begin{itemize}
    \item $(B_1,B_2):=(e +e_2,e_2-e_1)$ of $\H1(\Gamma,\mathbb{Z})$.
    \item $\tilde{B}:=\tilde{e}_1 + \tilde{e}_2$ of $\H1(\mathbb{T}E,\mathbb{Z})$.
    \end{itemize}
    We denote by $(\omega_1,\omega_2)$ and $\tilde{\omega}$ the canonical basis of tropical $1$-forms associated to  $(B_1,B_2)$, respectively $\tilde{B}$ (i.e. $\omega_i=\zeta_\Gamma(B_i)$ and $\tilde{\omega}=\zeta_{\mathbb{T}E}(\tilde{B})$). 
    We want to stress that the following arguments apply to \emph{any} choice of basis, the formulas below, however, will differ.
    
    An arbitrary cover of degree $d$ is captured by the following data:
\begin{itemize}
    \item a triple of \emph{winding numbers} $(n,n_1,n_2)\in \mathbb{N}_0^3$, where $n_i$ ($n$) counts how many times the edge $e_i$ ($e$) passes over the edge $\widetilde{e}_{2}$ ($\widetilde{e}_{1}$),
    \item a triple of \emph{dilation factors} $(d_e(\varphi),d_{e_1}(\varphi),d_{e_2}(\varphi))\in \mathbb{N}_0^3$,
    
\end{itemize}
that satisfy:
\begin{align}
   & nd_e(\varphi)+(n_1-1)d_{e_1}(\varphi)+(n_2-1)d_{e_2}(\varphi)=d.\\
   & (n-1)d_e(\varphi)+n_1d_{e_1}(\varphi)+n_2d_{e_2}(\varphi)=d.
\end{align}
\begin{figure}[H]
    \centering
    \tikzset{every picture/.style={line width=0.75pt}} %set default line width to 0.75pt        

\begin{tikzpicture}[x=0.75pt,y=0.75pt,yscale=-1,xscale=1]

%Shape: Ellipse [id:dp8660390501968123] 
\draw   (220.34,260.43) .. controls (220.37,247.14) and (269.65,236.37) .. (330.4,236.4) .. controls (391.15,236.42) and (440.37,247.21) .. (440.33,260.51) .. controls (440.3,273.81) and (391.02,284.57) .. (330.27,284.55) .. controls (269.52,284.53) and (220.3,273.73) .. (220.34,260.43) -- cycle ;
%Straight Lines [id:da5945040424239858] 
\draw    (440.33,260.5) ;
\draw [shift={(440.33,260.5)}, rotate = 0] [color={rgb, 255:red, 0; green, 0; blue, 0 }  ][fill={rgb, 255:red, 0; green, 0; blue, 0 }  ][line width=0.75]      (0, 0) circle [x radius= 3.35, y radius= 3.35]   ;
%Straight Lines [id:da5453948742270707] 
\draw    (220.33,260.44) ;
\draw [shift={(220.33,260.44)}, rotate = 0] [color={rgb, 255:red, 0; green, 0; blue, 0 }  ][fill={rgb, 255:red, 0; green, 0; blue, 0 }  ][line width=0.75]      (0, 0) circle [x radius= 3.35, y radius= 3.35]   ;
\draw   (317.02,278.65) .. controls (328.54,281.93) and (340.06,283.9) .. (351.59,284.56) .. controls (340.06,285.21) and (328.54,287.18) .. (317.02,290.46) ;
\draw   (345.05,231.39) .. controls (335.08,234.29) and (325.12,236.04) .. (315.15,236.62) .. controls (325.12,237.2) and (335.08,238.94) .. (345.05,241.84) ;
%Straight Lines [id:da5231446032186506] 
\draw    (221.33,80.71) ;
\draw [shift={(221.33,80.71)}, rotate = 0] [color={rgb, 255:red, 0; green, 0; blue, 0 }  ][fill={rgb, 255:red, 0; green, 0; blue, 0 }  ][line width=0.75]      (0, 0) circle [x radius= 3.35, y radius= 3.35]   ;
%Straight Lines [id:da7728429300432185] 
\draw    (449.33,81.71) ;
\draw [shift={(449.33,81.71)}, rotate = 0] [color={rgb, 255:red, 0; green, 0; blue, 0 }  ][fill={rgb, 255:red, 0; green, 0; blue, 0 }  ][line width=0.75]      (0, 0) circle [x radius= 3.35, y radius= 3.35]   ;
\draw   (285.33,53.67) .. controls (279,57) and (272.67,59) .. (266.33,59.67) .. controls (272.67,60.33) and (279,62.33) .. (285.33,65.67) ;
\draw   (286.33,150.67) .. controls (292.62,153.61) and (298.9,155.37) .. (305.18,155.96) .. controls (298.9,156.55) and (292.62,158.31) .. (286.33,161.25) ;
\draw   (375.33,83.67) .. controls (382.25,86.9) and (389.17,88.84) .. (396.09,89.49) .. controls (389.17,90.14) and (382.25,92.08) .. (375.33,95.32) ;
%Curve Lines [id:da686234217221418] 
\draw    (221.33,80.71) .. controls (219.87,131.13) and (264.93,152.17) .. (301.37,156.32) .. controls (314.77,157.84) and (327,157.09) .. (335.33,154.67) .. controls (366.33,145.67) and (255.33,143.67) .. (299.33,138.67) .. controls (343.33,133.67) and (358.33,133.67) .. (295.33,127.67) .. controls (232.33,121.67) and (465.33,134.67) .. (449.33,81.71) ;
%Curve Lines [id:da21987253917421623] 
\draw    (221.33,80.71) .. controls (229.33,88.67) and (246.33,100.67) .. (263.33,101.67) .. controls (280.33,102.67) and (322.33,104.36) .. (339.33,99.67) .. controls (356.34,94.97) and (270.33,101.67) .. (304.33,93.67) .. controls (338.33,85.67) and (366.33,74.67) .. (306.33,78.67) .. controls (246.33,82.67) and (397.33,97.67) .. (449.33,81.71) ;
%Curve Lines [id:da832796419608693] 
\draw    (221.33,80.71) .. controls (228.33,59.67) and (295.33,55.67) .. (330.33,55.67) .. controls (365.33,55.67) and (258.33,38.67) .. (297.33,36.67) .. controls (336.33,34.67) and (351.33,15.67) .. (299.33,19.67) .. controls (247.33,23.67) and (407.33,43.67) .. (449.33,81.71) ;
%Straight Lines [id:da24824002703052162] 
\draw    (320.33,174.67) -- (320.33,215.67) ;
\draw [shift={(320.33,217.67)}, rotate = 270] [color={rgb, 255:red, 0; green, 0; blue, 0 }  ][line width=0.75]    (10.93,-3.29) .. controls (6.95,-1.4) and (3.31,-0.3) .. (0,0) .. controls (3.31,0.3) and (6.95,1.4) .. (10.93,3.29)   ;

% Text Node
\draw (200,70) node [anchor=north west][inner sep=0.75pt]    {$P_0$};
% Text Node
\draw (453,72) node [anchor=north west][inner sep=0.75pt]    {$P_1$};
% Text Node
\draw (283.76,211.16) node [anchor=north west][inner sep=0.75pt]    {$\widetilde{e}_{1}$};
% Text Node
\draw (375.37,253.04) node [anchor=north west][inner sep=0.75pt]    {$\widetilde{e}_{2}$};
% Text Node
\draw (139.97,263.86) node [anchor=north west][inner sep=0.75pt]    {$\mathbb{T} E$};
% Text Node
\draw (245.58,39.45) node [anchor=north west][inner sep=0.75pt]    {$e$};
% Text Node
\draw (397.79,66.01) node [anchor=north west][inner sep=0.75pt]    {$e_{1}$};
% Text Node
\draw (392.27,122.24) node [anchor=north west][inner sep=0.75pt]    {$e_{2}$};
% Text Node
\draw (146.44,36.14) node [anchor=north west][inner sep=0.75pt]    {$\Gamma $};
% Text Node
\draw (334,185.4) node [anchor=north west][inner sep=0.75pt]    {$\varphi $};

\end{tikzpicture}

    \caption{Schematic picture of a finite cover of an elliptic curve by a curve of genus 2 of type "theta". }
    \label{figure_generalgenus2cover}
\end{figure}
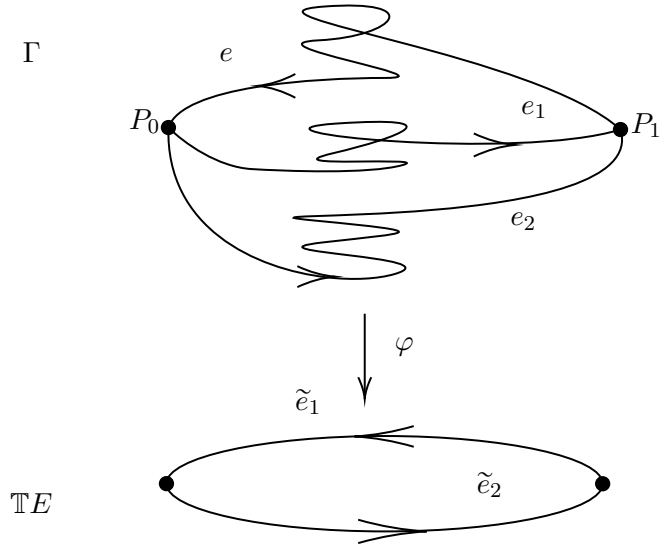
\end{convention}
\begin{lemma}\label{lemma_kernelofpushforwardinTA}
     Let $\Gamma$ be an abstract tropical curve of type "theta" and $\varphi: \Gamma \rightarrow \mathbb{T}E$ a tropical cover of degree $d$. In $\mathbb{T}\mathcal{A}$, the kernel of $\varphi_*$ (i.e. $\Ker(\varphi_*)_0$) is an elliptic curve of length $l_K$ given by the absolute value of the matrix product $vMw^t$, where
\begin{itemize}
    \item $M:=\begin{pmatrix}
l(e) +l(e_2) & l(e_2) \\
l(e_2) & l(e_1) +l(e_2)
\end{pmatrix}$ is the period matrix associated to the choice of basis from Convention \ref{convention_genus2coversofE}.
    \item $v:=(v_1,v_2)$ is any vector whose entries satisfy: $v_2d_e(\varphi)+v_1d_{e_1}(\varphi)=gcd(d_{e}(\varphi),d_{e_1}(\varphi))$.
    \item $w:=(\frac{n_1-n_2}{gcd(n_2-n_1,n+n_2-1)},\frac{n+n_2-1}{gcd(n_2-n_1,n+n_2-1)})$. %is interested in winding
\end{itemize}
\end{lemma}
\begin{proof}
Recalling that $\varphi_*$ is a surjective morphism from a $2$-dimensional to a $1$-dimensional torus, it is not surprising that (in $\mathbb{T}\mathcal{A}$) its kernel is a tav of dimension $1$, in other words an elliptic curve (see Remark \ref{remark_endowingEwithpp}). More interesting is the statement about $l_K$, since the length is an isomorphy-invariant of elliptic curves. So, in order to capture of the geometry of $\Ker(\varphi_*)_0$, we need to look at it in more detail: Its lattice representation (see Definition \ref{definition_(Co-)KernelImage})

\begin{align}
    (\Omega^1_\Gamma/\im({\varphi_*}^\#)^{\text{sat}},\ker({\varphi_*}_\#),[\cdot,\cdot]_K)
\end{align} allows us to access $l_K$ via integration on $\Gamma$, i.e. $l_K=|\int_B\omega|$, where $B$ and $\omega$ are basis of the lattices defining $\Ker(\varphi_*)_0$. To obtain a concrete expression for $l_K$ in terms of covering data, we compute $B$ and $\omega$ explicitly. \\
We start with $\ker({\varphi_*}_\#)$: The  map ${\varphi_*}_\#: \H1(\Gamma,\mathbb{Z})\rightarrow \H1(\mathbb{T}E,\mathbb{Z})$, sends a cycle $\sum e_i$ to the cycle $\sum \varphi(e_i)$. This yields 
\begin{align}
  &{\varphi_*}_\#(B_1)=((n_2-1)+(n-1)+1)\tilde{B} \text{ and }\\
  &{\varphi_*}_\#(B_2)=((n_2-1)-(n_1-1))\tilde{B},
\end{align}
for our choice of basis $(B_1,B_2)$ from Convention \ref{convention_genus2coversofE}. A contracted cycle $B=w_1B_1 +w_2B_2 \in \H1(\Gamma,\mathbb{Z}) $ generates $\ker({\varphi_*}_\#)$, if the pair of coefficients $(w_1,w_2)\in \mathbb{Z}^2$ in addition to satisfying
\begin{align}
  &  0= {\varphi_*}_\#(B) = w_1(n+n_2-1)\tilde{B} + w_2(n_2-n_1)\tilde{B}\\
  &  \Leftrightarrow 0= w_1(n+n_2-1)+ w_2(n_2-n_1),
\end{align}
 is relatively prime. This yields
 \begin{align}
     w:=(w_1,w_2)=\left( \frac{n_1-n_2}{gcd(n_2-n_1,n+n_2-1)},\frac{n+n_2-1}{gcd(n_2-n_1,n+n_2-1)}\right).
 \end{align}
 We continue with $\Omega^1_\Gamma/\im({\varphi_*}^\#)^{\text{sat}}$: The pullback of 1-forms ${\varphi_*}^\#$ is injective as ${\varphi_*}$ is surjective. Hence, $\im({\varphi_*}^\#)$ is generated by ${\varphi_*}^\#(\tilde{\omega})$ whose representation in terms of the basis $(\omega_1,\omega_2)$ from Convention \ref{convention_genus2coversofE} is
 \begin{align}
   {\varphi_*}^\#(\tilde{\omega}=d\tilde{e}_1+d\tilde{e}_2) &=  d_{e_1}(\varphi)de_1 + d_{e_2}(\varphi)de_2 + d_{e}(\varphi)de\\
   &= d_{e}(\varphi)\omega_1 + (d_{e_2}(\varphi)- d_{e}(\varphi))de_2 +d_{e_1}(\varphi)de_1\\
   &=d_{e}(\varphi)\omega_1- d_{e_1}(\varphi)\omega_2.
 \end{align}
 As a result, we get
 \begin{align}
     \im({\varphi_*}^\#)^{\text{sat}}=\langle \omega_Q:=\frac{d_{e}(\varphi)}{gcd(d_{e}(\varphi),d_{e_1}(\varphi))}\omega_1-\frac{d_{e_1}(\varphi)}{gcd(d_{e}(\varphi),d_{e_1}(\varphi))}\omega_2 \rangle,
 \end{align}
 whereby it should be emphasized that $gcd(d_{e}(\varphi),d_{e_1}(\varphi))\neq 0$ holds, even if $\varphi$ is not finite: Suppose otherwise. Then $d_{e}(\varphi)=d_{e_1}(\varphi)=0$ implies $d_{e_2}(\varphi)=0$ due to balancing. So, $\varphi$ cannot be surjective, a contradiction.\\
 Note, the quotient $\Omega^1_\Gamma/\im({\varphi_*}^\#)^{\text{sat}}$ is naturally identified with a complement of $\im({\varphi_*}^\#)^{\text{sat}}$ in $\H1(\Gamma,\mathbb{Z})$. To find one, we need to extend $\omega_Q$ to a basis of $\H1(\Gamma,\mathbb{Z})$ by adding a 1-form $\omega:=v_1\omega_1 + v_2\omega_2$. The requirement that the new pair $(\omega_Q,\omega)$ can be obtained from $(\omega_1,\omega_2)$ by an invertible coordinate change can be rephrased as condition imposed on the coefficients of $\omega$:
 \begin{align}
    v_2d_e(\varphi)+v_1d_{e_1}(\varphi)=gcd(d_{e}(\varphi),d_{e_1}(\varphi)).
 \end{align}
 Any such 1-form $\omega$ generates a complement and is therefore suitable for us (see Remark \ref{remark_integralisindependentfromchoiceofcomplement}).
This concludes our hunt for basis. The desired formula reads
 \begin{align}
     l_K=|\int_B\omega| =|\int_{w_1B_1+w_2B_2} (v_1\omega_1 + v_2\omega_2)|
 \end{align}
and can be compactly written as matrix product $l_K=|vMw^t|$, where $M$ is the period matrix associated to our choice of basis (see Convention \ref{convention_genus2coversofE}).
\end{proof}
\begin{remark}\label{remark_endowingEwithpp}
We interpret the 1-dimensional torus $\Ker(\varphi_*)_0=(\langle \omega \rangle, \langle B \rangle, [\cdot,\cdot]_K)$, where $\omega\in \Omega^1_\Gamma$ and $B\in \H1(\Gamma,\mathbb{Z})$ are lattice basis as in Lemma \ref{lemma_kernelofpushforwardinTA} such that $l_K:=\int_B \omega>0$, as Jacobian of a metric graph $\mathbb{T}E'$ with $l(\mathbb{T}E')=l_K$ via:
\begin{align}
    \Omega^1_{\mathbb{T}E'}=\langle \omega \rangle, \thinspace \H1(\mathbb{T}E',\mathbb{Z})=\langle B \rangle, \text{ and } \zeta_{\mathbb{T}E'},
\end{align}
the usual pp (Definition \ref{definition_JacobianandAbelJacobimap}) and, by abuse of notation, denote $\Jac(\mathbb{T}E')$ again by $\mathbb{T}E'$ (see beginning of Subsection \ref{subsection_curvesofg2coverg1}).
\end{remark}
\begin{remark}\label{remark_integralisindependentfromchoiceofcomplement}
Note that the vector $v$ from Lemma \ref{lemma_kernelofpushforwardinTA} is not uniquely determined: Solutions to
\begin{align}
    v_2d_e(\varphi)+v_1d_{e_1}(\varphi)=gcd(d_{e}(\varphi),d_{e_1}(\varphi))
\end{align} 
differ by elements of $\langle(\frac{d_{e}(\varphi)}{gcd(d_{e}(\varphi),d_{e_1}(\varphi))},\frac{-d_{e_1}(\varphi)}{gcd(d_{e}(\varphi),d_{e_1}(\varphi))})\rangle_\mathbb{Z}$, which corresponds to adding a multiple of the generator of $\im({\varphi_*}^\#)^{\text{sat}}$. The reason for this indeterminacy can be found in the proof of Lemma \ref{lemma_kernelofpushforwardinTA}. It arises from the identification of the quotient $\Omega^1_\Gamma/\im({\varphi_*}^\#)^{\text{sat}}$ with a complement of $\im({\varphi_*}^\#)^{\text{sat}}$ and does not affect the matrix product $vMw^t$, since we have
\begin{align}
    \int_B{\varphi_*}^\#(\omega)=\int_{{\varphi_*}_\#(B)}\omega=0.
\end{align}
\end{remark}
The transition to $\gls{Ab}$ is accompanied by a loss of structure (the kernel is no longer a tav). A trade-off that is necessary to fully understand it.

\begin{lemma}\label{lemma_Kernelofpushforwardabstract}
The group-theoretic kernel of the push-forward map consists of translates of the 1-dimensional torus $\Ker(\varphi_*)_0$ by elements of $\varphi^*\Jac_d(\mathbb{T}E)$ (the image of the $d$-torsion subgroup $\Jac_d(\mathbb{T}E)$ of $\Jac(\mathbb{T}E)$ under $\varphi^*$).
Moreover, the number of connected components is given by $|\ker(\gamma)|$, where $\gamma: \Jac(\Gamma)/{\Ker(\varphi_{*})}_0 \rightarrow \Jac(\mathbb{T}E)$ is the natural map.
\end{lemma}
The beginning of the following proof is similar to the proof of Proposition 6.1. (\cite{MR3782424}). As we work in a different setting, we present it here for the sake of completeness. Compare also to Proposition 4.7 in \cite{röhrle2024tropicalngonalconstruction} for a different characterization of the connected components of $\ker(f)$ for general morphism $f$.
\begin{proof}
The kernel of $\varphi_{*}$ is the image of $\Hom(\varphi^\#_{*})^{-1}\H1(\mathbb{T}E,\mathbb{Z})$ under the universal covering of $\Jac(\Gamma)$. Recalling that $\Hom(\varphi^\#_{*})$ is $\mathbb{R}$-linear we can write the solution set $\Hom(\varphi^\#_{*})^{-1}\H1(\mathbb{T}E,\mathbb{Z})$ as translation of the kernel of $\Hom(\varphi^\#_{*})$ by a set consisting of specific solutions to the equations
\begin{align}
  \Hom(\varphi^\#_{*})(x_C)=C \text{ for }  C\in \H1(\mathbb{T}E,\mathbb{Z}).
\end{align} 
We now compute such a specific solution. The isomorphisms $\H1(\Gamma,\mathbb{Z})\cong \Omega^1_\Gamma$ and $\H1(\mathbb{T}E,\mathbb{Z})\cong \Omega^1_{\mathbb{T}E}$ allow us to identify $\varphi_{*\#}$ with the push-forward of $1-$forms and consider the composition $\varphi_{*\#}\circ \varphi^\#_{*}: \Omega^1_{\mathbb{T}E} \rightarrow \Omega^1_{\mathbb{T}E} $. This is easily seen to be the multiplication-by-$d$ map. Hence, we have 
\begin{align}
    \Hom(\varphi^\#_{*})\circ \Hom(\varphi_{*\#}) = \Hom(\varphi_{*\#}\circ \varphi^\#_{*})= d \cdot id  
\end{align}
for the $\Hom$-duals and $\frac{1}{d}\Hom(\varphi_{*\#})(C)$ as solution of $\Hom(\varphi^\#_{*})(x_C)=C$. \\
We write
\begin{align}
    \Hom(\varphi^\#_{*})^{-1}\H1(\mathbb{T}E,\mathbb{Z})= \ker(\Hom(\varphi^\#_{*})) + \frac{1}{d}\Hom(\varphi_{*\#})\H1(\mathbb{T}E,\mathbb{Z})
\end{align}
and descend to $\Jac(\Gamma)$ to obtain the desired description of $ \ker(\varphi_{*})$. Remembering that $\Hom(\varphi_{*\#})$ is the lift of $\varphi^*$ to a map between the universal coverings (Lemma \ref{lemma_pushandpullaredual}) and using $\Jac_d(\mathbb{T}E)=\frac{\H1(\mathbb{T}E,\mathbb{Z})}{d\H1(\mathbb{T}E,\mathbb{Z})}$, this looks like:
\begin{align}
    \ker(\varphi_{*})&= \Ker(\varphi_{*})_0 + \frac{\Hom(\varphi_{*\#})\H1(\mathbb{T}E,\mathbb{Z})}{d\cdot \Hom(\varphi_{*\#})(\H1(\mathbb{T}E,\mathbb{Z})) \cap \H1(\Gamma,\mathbb{Z})}\\
    & = \Ker(\varphi_{*})_0 + \varphi^*\Jac_d(\mathbb{T}E).
\end{align}
Note that we have just written $ \ker(\varphi_{*})$ as union of cosets of its subgroup $ \Ker(\varphi_{*})_0$. We did even more: We found a \emph{finite} superset for a set of representatives, which is given by $\varphi^*\Jac_d(\mathbb{T}E)$. In other words, the number of connected components is finite and given by the index of $ \Ker(\varphi_{*})_0$ in $ \ker(\varphi_{*})$. We can consider the following commutative exact diagram in the category of abelian groups (see Figure \ref{figure_inlemma_Kernelofpushforwardabstract})

 \begin{figure}[ht]
        \centering
         \begin{tikzpicture}[>=triangle 60]
\matrix[matrix of math nodes,column sep={60pt,between origins},row
sep={60pt,between origins},nodes={asymmetrical rectangle}] (s)
{
&|[name=ka]| 0 &|[name=kb]| 0 &|[name=kc]| \ker(\gamma) \\
 &|[name=A]| \Ker(\varphi_{*})_0 &|[name=B]| \Jac(\Gamma) &|[name=C]| \frac{\Jac(\Gamma)}{{\Ker(\varphi_{*})}_0} &|[name=01]| 0 \\
|[name=02]| 0 &|[name=A']|\ker(\varphi_{*}) &|[name=B']| \Jac(\Gamma) &|[name=C']| \Jac(\mathbb{T}E) &|[name=no]| 0 \\
&|[name=ca]|\frac{\ker(\varphi_{*})}{\Ker(\varphi_{*})_0} &|[name=cb]| 0 &|[name=cc]| 0 \\
};
\draw[->] 
          (ka) edge (A)
          (kb) edge (B)
          (kc) edge (C)
          (A) edge  (B)
          (B) edge (C)
          (C) edge (01)
          (A) edge  (A')
          (B) edge node[auto] {$id$} (B')
          (C) edge node[auto] {$\gamma$} (C')
          (02) edge (A')
          (A') edge  (B')
          (B') edge  (C')
          (A') edge (ca)
          (B') edge (cb)
          (C') edge (cc)
          (C') edge (no)
          
;
\draw[->,gray] (ka) edge (kb)
               (kb) edge (kc)
               (ca) edge (cb)
               (cb) edge (cc)
;
\draw[->,gray,rounded corners] (kc) -| node[auto,text=black,pos=.7]
{\(\partial\)} ($(01.east)+(.5,0)$) |- ($(B)!.35!(B')$) -|
($(02.west)+(-.5,0)$) |- (ca);
\end{tikzpicture}
\caption{Applying the snake Lemma in the proof of Lemma \ref{lemma_Kernelofpushforwardabstract}.}
\label{figure_inlemma_Kernelofpushforwardabstract}
        \end{figure}
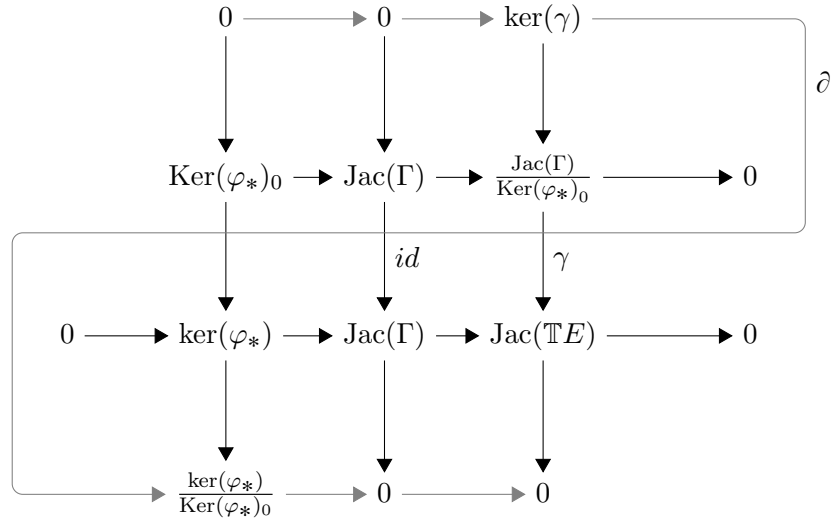 
        and apply the snake Lemma to obtain an isomorphism $\frac{\ker(\varphi_{*})}{\Ker(\varphi_{*})_0} \xrightarrow{\partial} \ker(\gamma)$. In particular we see that the surjective map $\gamma$ has a finite kernel, is therefore an isogeny, whose kernel counts the number of connected components as claimed.
\end{proof}
We give an explicit description of the isogeny $\gamma$ from Lemma \ref{lemma_Kernelofpushforwardabstract} and compute the number of connected components in terms of data from our cover.
\begin{lemma}\label{lemma_explicitdescriptionofgamma}
     The  tav $\Jac(\Gamma)/{\Ker(\varphi_{*})}_0$ is an elliptic curve of length $\tilde{l}=|w^Q M {v^Q}^t|$, where
     \begin{itemize}
    \item $M$ is the period matrix associated to the choice of basis from convention \ref{convention_genus2coversofE}.
    \item $v^Q:=(v^Q_1,v^Q_2)$ is any vector whose entries satisfy: \\
    $v^Q_1(n+n_2-1)-v^Q_2(n_1-n_2)=gcd(n_2-n_1,n+n_2-1)$.
    \item $w^Q:=(-\frac{d_{e_1}(\varphi)}{gcd(d_{e}(\varphi),d_{e_1}(\varphi))},\frac{d_{e}(\varphi)}{gcd(d_{e}(\varphi),d_{e_1}(\varphi))})$.
\end{itemize}
     Under the identifications $\Jac(\mathbb{T}E)\cong \mathbb{R}/l\mathbb{Z}$ and $\Jac(\Gamma)/{\Ker(\varphi_{*})}_0\cong \mathbb{R}/\tilde{l}\mathbb{Z}$ the isogeny $\gamma$ from Lemma \ref{lemma_Kernelofpushforwardabstract} takes the form:
     \begin{align}
          \mathbb{R}/\tilde{l}\mathbb{Z} \rightarrow \mathbb{R}/l\mathbb{Z}, z \mapsto az,
     \end{align}
     where $a=gcd(d_{e}(\varphi),d_{e_1}(\varphi))$.
\end{lemma}
\begin{proof}
    We derive a lattice representation $(\ker(i^\#),\Omega^1_\Gamma/\im(i_\#),[\cdot,\cdot]^t_K)$ of the quotient variety $\Jac(\Gamma)/{\Ker(\varphi_{*})}_0$ from Lemma \ref{lemma_quotientsoftav} (or Remark \ref{remark_latticerepresentationofquotient}) by choosing the natural inclusion $i=(i^\#,i_\#)$ defined by the maps
    \begin{align}
     i^\#: \Omega^1_\Gamma \twoheadrightarrow \Omega^1_\Gamma/\im({\varphi_*}^\#)^{\text{sat}} \text{ and } i_\#: \ker(\varphi_{*\#}) \hookrightarrow \H1(\Gamma,\mathbb{Z})
    \end{align}
    as embedding of ${\Ker(\varphi_{*})}_0 $ into $\Jac(\Gamma)$. Just like in Lemma \ref{lemma_kernelofpushforwardinTA} we determine $\tilde{l}$ by expressing explicit basis $\omega_Q$ and $B_Q$ of the lattices $\ker(i^\#)$ and $\Omega^1_\Gamma/\im(i_\#)$ in terms of our chosen basis from Convention \ref{convention_genus2coversofE}. The equality $\ker(i^\#)=\im({\varphi_*}^\#)^{\text{sat}}$ is immediate from the definition of $i$. Relying on the results of Lemma \ref{lemma_kernelofpushforwardinTA} we set    
    \begin{align}
    \omega_Q:=\frac{d_{e}(\varphi)}{gcd(d_{e}(\varphi),d_{e_1}(\varphi))}\omega_2-\frac{d_{e_1}(\varphi)}{gcd(d_{e}(\varphi),d_{e_1}(\varphi))}\omega_1.
 \end{align}
 Since $i_\#$ is just the inclusion, we have $\im(i_\#)=\ker(\varphi_{*\#})=\langle B \rangle$ (Lemma \ref{lemma_kernelofpushforwardinTA}), where $B=w_1B_1 +w_2B_2 \in \H1(\Gamma,\mathbb{Z}) $ with coefficients
 \begin{align}
     (w_1,w_2)=\left(\frac{n_1-n_2}{gcd(n_2-n_1,n+n_2-1)},\frac{n+n_2-1}{gcd(n_2-n_1,n+n_2-1)}\right).
 \end{align}
 Any cycle $B_Q:=v^Q_1B_1 + v^Q_2B_2$ such that
 \begin{align}
     v^Q_1(n+n_2-1)-v^Q_2(n_1-n_2)=gcd(n_2-n_1,n+n_2-1)
 \end{align}
 holds, generates a complement of $\im(i_\#)$ in $\H1(\Gamma,\mathbb{Z})$. We identify $\langle B_Q \rangle$ with the quotient $\Omega^1_\Gamma/\im(i_\#)$. The length $\tilde{l}$ is then given by the absolute value of the integral
 \begin{align}
     \int_{B_Q}\omega_Q =\int_{w^Q_1B_1+w^Q_2B_2} (v^Q_1\omega_1 + v^Q_2\omega_2),
 \end{align}
which written as matrix product reads $w^Q M {v^Q}^t$.\\%the integral written as a matrix product is wMw
As a closer look at Figure \ref{figure_inlemma_Kernelofpushforwardabstract} reveals, $\gamma$ just sends an equivalence class $[x]$ to the image of a representative under the push-forward, i.e. $\gamma([x])=\varphi_*(x)$. Hence, $\gamma^\#: \Omega^1_{\mathbb{T}E} \rightarrow \ker(i^\#), $ is defined by $\omega \mapsto \varphi^\#_*(\omega)$. We recognize this as multiplication by $a^\#$, where $a^\#$ satisfies $\varphi^\#_*(\omega)=a^\# \omega_Q$ and is therefore equal to $gcd(d_{e}(\varphi),d_{e_1}(\varphi))$. This together with the knowledge of $\tilde{l}$ and $l$ is sufficient to give a concrete description of the isogeny $\gamma$ under the identifications $\Jac(\mathbb{T}E)\cong \mathbb{R}/l\mathbb{Z}$ and $\Jac(\Gamma)/{\Ker(\varphi_{*})}_0\cong \mathbb{R}/\tilde{l}\mathbb{Z}$ (see Lemma \ref{lemma_isogeniesbetweenellipticcurves}):
     \begin{align}
          \mathbb{R}/\tilde{l}\mathbb{Z} \rightarrow \mathbb{R}/l\mathbb{Z}, z \mapsto a^\#z.
     \end{align}
\end{proof}
\begin{remark}\label{remark_acomputationforgamma}
A priory, all arguments, which lead us to $a^\#$, are based on a particular choice of basis (Convention \ref{convention_genus2coversofE}). This should not be the case: The lengths $\tilde{l}$ and $l$ are intrinsic characteristics of the curves, so should the description of $\gamma$ be. And indeed, the balancing condition implies that $gcd(d_{e}(\varphi),d_{e_1}(\varphi))=gcd(d_{e}(\varphi),d_{e_2}(\varphi))$, which in turn guarantees that interchanging the roles of $e_1$ and $e_2$ does not affect $a^\#$.\\
Using Lemma \ref{lemma_isogeniesbetweenellipticcurves} we can also describe the full pair $\gamma=(\gamma^\#,\gamma_\#)$ since
    $\gamma_\#$ is the multiplication by $a_\#:=\frac{\tilde{l}a^\#}{l}$.
\end{remark}
\begin{corollary}\label{corollary_numberofccinkerpushforward}
    The number of connected components of $\ker(\varphi_*)$ is equal to $k+1$, where $k:=max \{ j\in \mathbb{N}: \frac{j}{ a_\#}<1\}$.
    \end{corollary}
    \begin{proof}
        By Lemma \ref{lemma_Kernelofpushforwardabstract} we know that the number of connected components is the order of $\ker(\gamma)$. We can use Lemma \ref{lemma_explicitdescriptionofgamma} and Remark \ref{remark_acomputationforgamma} to write it down explicitly:
        \begin{align}
            \ker(\gamma)=\{x:\thinspace a^\#x\in l\mathbb{Z} \}=\{x:\thinspace \frac{l\cdot a_\#}{\tilde{l}}x\in l\mathbb{Z} \}.
        \end{align}
        and conclude by specifying a set of representatives:
        \begin{align}
            \{0, \frac{\tilde{l}}{a_\#},...,\frac{k\cdot\tilde{l}}{a_\#} \}, \text{ where } k:=max \{ j\in \mathbb{N}: \frac{j}{ a_\#}<1\}.
        \end{align}
    \end{proof}
%    \paragraph{\emph{Collection of Results for the non-finite case.}} 
\subsubsection{Type dumbbell}\label{subsubsection_typedb}
\begin{figure}[H]
      \centering
  \tikzset{every picture/.style={line width=0.5pt}} %set default line width to 0.75pt        

\begin{tikzpicture}[x=0.7pt,y=0.7pt,yscale=-0.75,xscale=0.75]
%uncomment if require:\path (0,300); %set diagram left start at 0, and has height of 300
%Shape: Ellipse [id:dp8660390501968123] 
\draw   (54.34,248.43) .. controls (54.37,235.14) and (103.65,224.37) .. (164.4,224.4) .. controls (225.15,224.42) and (274.37,235.21) .. (274.33,248.51) .. controls (274.3,261.81) and (225.02,272.57) .. (164.27,272.55) .. controls (103.52,272.53) and (54.3,261.73) .. (54.34,248.43) -- cycle ;
%Straight Lines [id:da5945040424239858] 
\draw    (274.33,248.5) ;
\draw [shift={(274.33,248.5)}, rotate = 0] [color={rgb, 255:red, 0; green, 0; blue, 0 }  ][fill={rgb, 255:red, 0; green, 0; blue, 0 }  ][line width=0.75]      (0, 0) circle [x radius= 3.35, y radius= 3.35]   ;
\draw   (151.02,266.65) .. controls (162.54,269.93) and (174.06,271.9) .. (185.59,272.56) .. controls (174.06,273.21) and (162.54,275.18) .. (151.02,278.46) ;
\draw   (179.05,219.39) .. controls (169.08,222.29) and (159.12,224.04) .. (149.15,224.62) .. controls (159.12,225.2) and (169.08,226.94) .. (179.05,229.84) ;
%Straight Lines [id:da24824002703052162] 
\draw    (154.33,162.67) -- (154.33,203.67) ;
\draw [shift={(154.33,205.67)}, rotate = 270] [color={rgb, 255:red, 0; green, 0; blue, 0 }  ][line width=0.75]    (10.93,-3.29) .. controls (6.95,-1.4) and (3.31,-0.3) .. (0,0) .. controls (3.31,0.3) and (6.95,1.4) .. (10.93,3.29)   ;
%Straight Lines [id:da20105247618521016] 
\draw    (273,50.32) -- (275,116.46) ;
\draw [shift={(275,116.46)}, rotate = 88.27] [color={rgb, 255:red, 0; green, 0; blue, 0 }  ][fill={rgb, 255:red, 0; green, 0; blue, 0 }  ][line width=0.75]      (0, 0) circle [x radius= 3.35, y radius= 3.35]   ;
\draw [shift={(273,50.32)}, rotate = 88.27] [color={rgb, 255:red, 0; green, 0; blue, 0 }  ][fill={rgb, 255:red, 0; green, 0; blue, 0 }  ][line width=0.75]      (0, 0) circle [x radius= 3.35, y radius= 3.35]   ;
%Curve Lines [id:da1638898826694648] 
\draw    (68.33,114.95) .. controls (-37.67,133.79) and (277.33,157.9) .. (275,116.46) ;
%Curve Lines [id:da414515407083913] 
\draw    (68.33,114.95) .. controls (128.33,106.66) and (355.33,106.66) .. (213.33,126.25) ;
%Curve Lines [id:da910352681495392] 
\draw    (109.33,99.13) .. controls (-61.67,121.73) and (164.33,131.53) .. (213.33,126.25) ;
%Curve Lines [id:da9531402590250089] 
\draw    (109.33,99.13) .. controls (178.33,93.85) and (276.33,103.65) .. (275,116.46) ;
%Curve Lines [id:da674359026640412] 
\draw    (66.33,48.95) .. controls (-39.67,66.04) and (275.33,87.91) .. (273,50.32) ;
%Curve Lines [id:da3879361777739134] 
\draw    (66.33,48.95) .. controls (126.33,41.44) and (353.33,41.44) .. (211.33,59.21) ;
%Curve Lines [id:da6601218465090932] 
\draw    (107.33,34.6) .. controls (-63.67,55.1) and (162.33,63.99) .. (211.33,59.21) ;
%Curve Lines [id:da8001251255174269] 
\draw    (107.33,34.6) .. controls (180.33,26.4) and (292.33,39.39) .. (273,50.32) ;
\draw   (150.3,132.77) .. controls (160.07,136.3) and (169.82,138.39) .. (179.57,139.03) .. controls (169.84,139.84) and (160.11,142.09) .. (150.41,145.79) ;
\draw   (144.3,64.77) .. controls (154.07,68.3) and (163.82,70.39) .. (173.57,71.03) .. controls (163.84,71.84) and (154.11,74.09) .. (144.41,77.79) ;
%Shape: Ellipse [id:dp17652719783356452] 
\draw   (349.37,251.97) .. controls (349.4,238.67) and (398.68,227.91) .. (459.43,227.93) .. controls (520.18,227.95) and (569.4,238.75) .. (569.36,252.05) .. controls (569.33,265.34) and (520.05,276.11) .. (459.3,276.09) .. controls (398.55,276.06) and (349.33,265.27) .. (349.37,251.97) -- cycle ;
%Straight Lines [id:da22211045753108127] 
\draw    (569.37,252.04) ;
\draw [shift={(569.37,252.04)}, rotate = 0] [color={rgb, 255:red, 0; green, 0; blue, 0 }  ][fill={rgb, 255:red, 0; green, 0; blue, 0 }  ][line width=0.75]      (0, 0) circle [x radius= 3.35, y radius= 3.35]   ;
\draw   (446.05,270.19) .. controls (457.57,273.47) and (469.1,275.44) .. (480.62,276.09) .. controls (469.1,276.75) and (457.57,278.72) .. (446.05,282) ;
\draw   (474.08,222.93) .. controls (464.11,225.83) and (454.15,227.57) .. (444.18,228.15) .. controls (454.15,228.73) and (464.11,230.47) .. (474.08,233.38) ;
%Straight Lines [id:da1835834115478605] 
\draw    (457.37,157.2) -- (457.37,198.2) ;
\draw [shift={(457.37,200.2)}, rotate = 270] [color={rgb, 255:red, 0; green, 0; blue, 0 }  ][line width=0.75]    (10.93,-3.29) .. controls (6.95,-1.4) and (3.31,-0.3) .. (0,0) .. controls (3.31,0.3) and (6.95,1.4) .. (10.93,3.29)   ;
%Curve Lines [id:da9682524857098209] 
\draw    (361.37,52.49) .. controls (255.37,69.58) and (570.37,91.44) .. (568.03,53.86) ;
%Curve Lines [id:da5467395089378191] 
\draw    (361.37,52.49) .. controls (421.37,44.97) and (648.37,44.97) .. (506.37,62.74) ;
%Curve Lines [id:da8222077704639377] 
\draw    (402.37,38.14) .. controls (231.37,58.64) and (457.37,67.53) .. (506.37,62.74) ;
%Curve Lines [id:da5814243319095942] 
\draw    (402.37,38.14) .. controls (475.37,29.94) and (587.37,42.92) .. (568.03,53.86) ;
\draw   (439.33,68.3) .. controls (449.1,71.84) and (458.85,73.93) .. (468.6,74.57) .. controls (458.87,75.38) and (449.14,77.63) .. (439.44,81.33) ;
%Straight Lines [id:da1680480972727767] 
\draw    (568.03,53.86) -- (567.33,120.67) ;
\draw [shift={(567.33,120.67)}, rotate = 90.6] [color={rgb, 255:red, 0; green, 0; blue, 0 }  ][fill={rgb, 255:red, 0; green, 0; blue, 0 }  ][line width=0.75]      (0, 0) circle [x radius= 3.35, y radius= 3.35]   ;
\draw [shift={(568.03,53.86)}, rotate = 90.6] [color={rgb, 255:red, 0; green, 0; blue, 0 }  ][fill={rgb, 255:red, 0; green, 0; blue, 0 }  ][line width=0.75]      (0, 0) circle [x radius= 3.35, y radius= 3.35]   ;
%Shape: Ellipse [id:dp698202279477234] 
\draw   (567.52,169.65) .. controls (562.03,169.67) and (557.53,158.72) .. (557.48,145.2) .. controls (557.43,131.67) and (561.84,120.69) .. (567.33,120.67) .. controls (572.83,120.65) and (577.32,131.59) .. (577.37,145.12) .. controls (577.43,158.65) and (573.01,169.63) .. (567.52,169.65) -- cycle ;
\draw   (559.87,139.44) .. controls (558.51,144.11) and (557.73,148.77) .. (557.53,153.42) .. controls (557.15,148.78) and (556.18,144.16) .. (554.63,139.54) ;

% Text Node
\draw (276,27) node [anchor=north west][inner sep=0.75pt]    {$P_0 $};
% Text Node
\draw (276,120) node [anchor=north west][inner sep=0.75pt]    {$P_1$};
% Text Node
\draw (209.37,241.04) node [anchor=north west][inner sep=0.75pt]    {$\widetilde{e}$};
% Text Node
\draw (298.97,245.86) node [anchor=north west][inner sep=0.75pt]    {$\mathbb{T} E$};
% Text Node
\draw (26.79,27.01) node [anchor=north west][inner sep=0.75pt]    {$e_{1}$};
% Text Node
\draw (26.27,98.24) node [anchor=north west][inner sep=0.75pt]    {$e_{2}$};
% Text Node
\draw (304.44,77.14) node [anchor=north west][inner sep=0.75pt]    {$\Gamma $};
% Text Node
\draw (168,173.4) node [anchor=north west][inner sep=0.75pt]    {$\varphi $};
% Text Node
%\draw (412.79,202.69) node [anchor=north west][inner sep=0.75pt]    {$\widetilde{e_{1}}$};
% Text Node
\draw (504.4,244.58) node [anchor=north west][inner sep=0.75pt]    {$\widetilde{e}$};
% Text Node
\draw (321.82,30.54) node [anchor=north west][inner sep=0.75pt]    {$e_{1}$};
% Text Node
\draw (531.31,130.78) node [anchor=north west][inner sep=0.75pt]    {$e_{2}$};
% Text Node
\draw (471.03,167.94) node [anchor=north west][inner sep=0.75pt]    {$\varphi $};

\end{tikzpicture}
\caption{Schematic picture of two non-finite covers of an elliptic curve by a curve of genus 2 of "dumbell type".}
       \label{figure_dumbbellcover}
   \end{figure}
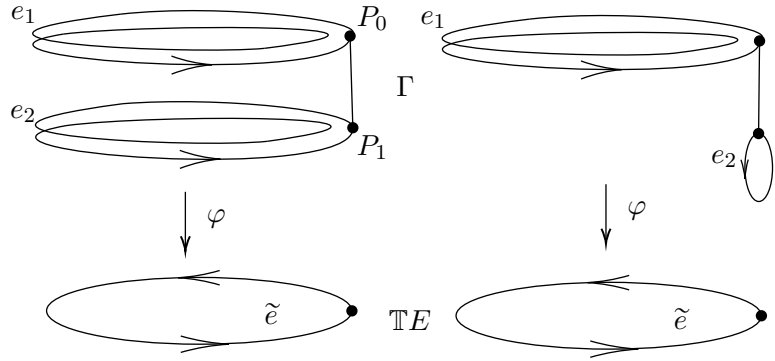

The case where $\varphi$ is finite is already addressed in Subsection \ref{subsubsection_typetheta}. Therefore, working with a curve $\Gamma$ of type "dumbbell" as in Figure \ref{figure_dumbbellcover} automatically leads us to the non-finite case since: 
Whenever $\Gamma$ maps to a curve of genus 1, at least the edge connecting $P_0$ to $P_1$ is contracted (a sketch of all other possibilities is shown in Figure \ref{figure_dumbbellcover}).

To formulate Lemma \ref{lemma_kernelofpushforwardinTA} and \ref{lemma_explicitdescriptionofgamma} in this case, we fix as before:
\begin{itemize}
    \item the homology basis $(B_1,B_2):=(e_1,e_2)$ of $\H1(\Gamma,\mathbb{Z})$ and $\tilde{B}:=\tilde{e}$ of $\H1(\mathbb{T}E,\mathbb{Z})$,
    \item the canonically associated basis $(\omega_1,\omega_2)$ of $\Omega^1_{\Gamma}(\mathbb{Z})$ and $\tilde{\omega}$ of $\Omega^1_{\mathbb{T}E}(\mathbb{Z})$,
    \end{itemize}
    and characterize a cover of degree $d$ by a tuple of \emph{winding numbers} $(n_1,n_2)\in \mathbb{N}_0^2$ together with a tuple of \emph{dilation factors} $(d_{e_1}(\varphi),d_{e_2}(\varphi))\in \mathbb{N}_0^2$ that satisfy $n_1d_{e_1}(\varphi)+n_2d_{e_2}(\varphi)=d$.
    \begin{lemma}\label{lemma_kernelpushforwardinTA_DB}
           Let $ \varphi: \Gamma \rightarrow \mathbb{T}E$ be a degree $d$ cover as in Figure \ref{figure_dumbbellcover}. In $\mathbb{T}\mathcal{A}$, the kernel of $\varphi_*$ is an elliptic curve of length $l_K$ given by the absolute value of the matrix product $vMw^t$, where
\begin{itemize}
    \item $M:=\begin{pmatrix}
l(e_1) & 0\\
0 & l(e_2)
\end{pmatrix}$ is the period matrix associated to the choice of basis.
    \item $v:=(v_1,v_2)$ is any vector whose entries satisfy: $v_2d_{e_1}(\varphi)-v_1d_{e_2}(\varphi)=gcd(d_{e_1}(\varphi),d_{e_2}(\varphi))$.
    \item $w:=(\frac{-n_2}{gcd(n_1,n_2)},\frac{n_1}{gcd(n_1,n_2)})$. 
\end{itemize}
The  tav $\Jac(\Gamma)/{\Ker(\varphi_{*})}_0$ is an elliptic curve of length $\tilde{l}=|w^Q M {v^Q}^t|$, where
 \begin{itemize}
    \item $v^Q:=(v^Q_1,v^Q_2)$ is any vector whose entries satisfy: \\
    $-(n_2v^Q_1+n_1v^Q_2)=gcd(n_1,n_2)$.
    \item $w^Q:=(\frac{d_{e_1}(\varphi)}{gcd(d_{e_1}(\varphi),d_{e_2}(\varphi))},\frac{d_{e_2}(\varphi)}{gcd(d_{e_1}(\varphi),d_{e_2}(\varphi))})$.
\end{itemize}
    \end{lemma}
    \subsection{The pull-back}\label{subsection_pull-back} 
    Changing gears, we now work with divisors to study the dual of $\varphi_*$. Recall that from this perspective $\varphi^{*}$ is simply the pull-back of divisors (see Remark \ref{remark_JacobianandPicardvariety}). For background on divisors on tropical curves see e.g. \cite{MR3782424}, \cite{MR3343870}, \cite{MR2448666} and references therein.
    
    We are still interested in covers of the form $\varphi: \Gamma \rightarrow \mathbb{T}E$, but do not require $\Gamma$ to have genus 2. Just as in Subsection \ref{subsection_push-forward}, we approach $\varphi^{*}$ by understanding its kernel $\ker(\varphi^{*})$. We proceed as follows: The groundwork will be laid in Lemma \ref{lemma_foundationkernelofpullback}, which we formulate for \emph{any} harmonic map $\varphi: \Gamma \rightarrow \tilde{\Gamma}$. This means we find a superset of  $\ker(\varphi^{*})$ whose elements are easily described and allow for concrete manipulation.
    This will be used later in Construction \ref{construction_kernelofpullback}.
\begin{lemma}\label{lemma_foundationkernelofpullback}
Let $\varphi: \Gamma \rightarrow \tilde{\Gamma}$ be a tropical cover of degree $d$. Then $\ker(\varphi^{*}) \subset \Jac_d(\tilde{\Gamma})$ and we have
\begin{align}
    D\in \ker(\varphi^{*})  \text{ if and only if all slopes of } \varphi^{*} f \text{ are divisible by } d,
\end{align}
where $f$ is a rational function satisfying $d\cdot D=div(f)$.
\end{lemma}
\begin{proof}
We obtain a drastic simplification of $\varphi^{*}$ by post-composing it with $\varphi_{*}$: As is easily verified, we have that $\varphi_{*}\varphi^{*}$ is given by the multiplication-by-$d$ map on $\Jac(\tilde{\Gamma})$. This immediately gives us the desired superset for $\ker(\varphi^{*})$ as the kernel of the composition, i.e. $\Jac_d(\tilde{\Gamma})$, and proves the first statement.\\
Since $D$ has order $d$, there exists a piecewise linear function $f: \tilde{\Gamma} \rightarrow \mathbb{R}$ such that $d\cdot D=div(f)$. Suppose $\varphi^{*}D=div(g)$ for $g:\Gamma \rightarrow \mathbb{R}$. Then $div(d\cdot g)=\varphi^{*}(d\cdot D)=div(\varphi^{*} f) $ (recall: $div(\varphi^{*}(f))=\varphi^{*}(div(f))$), shows that $\varphi^{*} f-d \cdot g$ is constant and $g=\frac{1}{d}\varphi^{*} f + c$ for $c\in \mathbb{R}$. Thus, $g$ is only well-defined if $\frac{1}{d}\varphi^{*} f$ is, i.e. all slopes are divisible by $d$.  
\end{proof}
As promised, applying Lemma \ref{lemma_foundationkernelofpullback} to the case where $\tilde{\Gamma}$ is of genus $1$ provides a simple description of divisors that are potentially contained in $\ker(\varphi^{*}) $: By restricting the isomorphism $\mathbb{T}E \rightarrow \Jac(\mathbb{T}E)$ to $\mathbb{T}E[d]$ we see that we can write any $D\in \Jac_d(\mathbb{T}E)$ as $D_P:=P-P_0$, where $P\in \mathbb{T}E[d]$ and $P_0\in \mathbb{T}E$ denotes the identity (i.e. the origin in $\mathbb{T}E\cong \mathbb{R}/l\mathbb{Z}$).\\
%\question{Wenn ich von ellitptischen Kurven spreche schreibe ich sie als Rmod lattice, wenn ich sie als abelsche varietät auffasse dann als tupel von lattices}
\begin{construction}\label{construction_kernelofpullback}
For any point $P\in \mathbb{T}E=\mathbb{R}/l\mathbb{Z}$ let $\overline{PP_0}$ ($\overline{P_0P}$) be the oriented edge of length $l_1$ ($l_2$) connecting $P$ with $P_0$ ($P_0$ with $P$). If $P\in l\mathbb{Q}/l\mathbb{Z}$, we have $ord(P)< \infty$. Hence, there exists a rational function $f$ on $\mathbb{T}E$ whose set of zeroes and poles is $\{P,P_0\}$. We can translate the statement $div(f)=ord(P)\cdot D_P$ into a linear system
    \begin{equation}
\left\{
\begin{alignedat}{4}
 %R & L  &  R & L   
s_1 & +{} & s_2 &  = ord_p(f) & = ord(P) \\
- s_1 & -{} & s_2 &  = ord_{P_0}(f) & = -ord(P), \\
l_1s_1 & -{} & (l-l_1)s_2 & = 0\\ 
\end{alignedat}
\right.
\end{equation} where \begin{itemize}
    \item $l_1$ encodes the position of $p$,
    \item $(s_1,s_2)$ is the data specifying $f$, i.e. $s_1$ ($s_2$) denotes the slope of $f$ on $\overline{PP_0}$ ($\overline{P_0P}$),
\end{itemize}
 and solve for $s_1=\frac{ord(P)(l-l_1)}{l}$ and $s_2=\frac{ord(P)l_1}{l}$.
 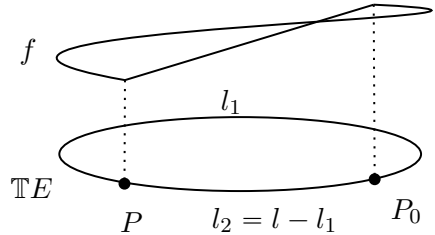
\begin{figure}[H]
     \centering
     \tikzset{every picture/.style={line width=0.75pt}} %set default line width to 0.75pt        

\begin{tikzpicture}[x=0.7pt,y=0.7pt,yscale=-0.8,xscale=0.8]

%Shape: Ellipse [id:dp5518543382197196] 
\draw   (98,211) .. controls (98,197.19) and (152.7,186) .. (220.17,186) .. controls (287.64,186) and (342.33,197.19) .. (342.33,211) .. controls (342.33,224.81) and (287.64,236) .. (220.17,236) .. controls (152.7,236) and (98,224.81) .. (98,211) -- cycle ;
%Straight Lines [id:da9591392563212657] 
\draw    (142,231) ;
\draw [shift={(142,231)}, rotate = 0] [color={rgb, 255:red, 0; green, 0; blue, 0 }  ][fill={rgb, 255:red, 0; green, 0; blue, 0 }  ][line width=0.75]      (0, 0) circle [x radius= 3.35, y radius= 3.35]   ;
%Straight Lines [id:da6392866171227873] 
\draw    (311,228) ;
\draw [shift={(311,228)}, rotate = 0] [color={rgb, 255:red, 0; green, 0; blue, 0 }  ][fill={rgb, 255:red, 0; green, 0; blue, 0 }  ][line width=0.75]      (0, 0) circle [x radius= 3.35, y radius= 3.35]   ;
%Curve Lines [id:da8575691576753361] 
\draw    (142.33,161) .. controls (-63.67,125) and (498.33,119) .. (310.33,110) ;
%Straight Lines [id:da356784853187901] 
\draw    (142.33,161) -- (310.33,110) ;
%Straight Lines [id:da4704050952257315] 
\draw  [dash pattern={on 0.84pt off 2.51pt}]  (142.33,161) -- (142,231) ;
%Straight Lines [id:da792937360888571] 
\draw  [dash pattern={on 0.84pt off 2.51pt}]  (310.33,110) -- (311,228) ;

% Text Node
\draw (320,239.4) node [anchor=north west][inner sep=0.75pt]    {$P_{0}$};
% Text Node
\draw (137,248.4) node [anchor=north west][inner sep=0.75pt]    {$P$};
% Text Node
\draw (205,165.4) node [anchor=north west][inner sep=0.75pt]    {$l_{1}$};
% Text Node
\draw (199,245.4) node [anchor=north west][inner sep=0.75pt]    {$l_{2} =l-l_{1}$};
% Text Node
\draw (69,131.4) node [anchor=north west][inner sep=0.75pt]    {$f$};
% Text Node
\draw (63,225.4) node [anchor=north west][inner sep=0.75pt]    {$\mathbb{T} E$};

\end{tikzpicture}

     \caption{A rational function $f$ on $\mathbb{T}E$ as described in Construction \ref{construction_kernelofpullback}.}
     \label{figure_constructionforanalysisofkernelpullback}
 \end{figure}
 As the pullback of $f$ along $\varphi$ rescales slopes by local dilation factors, the divisor $\varphi^{*} D_P$ is principal if and only if every edge $e\in E(\Gamma)$ satisfies:
\begin{itemize}
    \item If $e \mapsto e'$ for $e' \subset \overline{PP_0}$, then $d_e(\varphi)\frac{(l-l_1)}{l}\in \mathbb{Z}$ holds.
    \item If $e \mapsto e'$ for $e' \subset \overline{P_0P}$, then $d_e(\varphi)\frac{l_1}{l}\in \mathbb{Z}$ holds.
\end{itemize}
We package this information into a continuous function $q_\Gamma: [0,l] \rightarrow \mathbb{R}^{|E(\Gamma)|}$ that assigns to each $l_1 \in [0,l]$ the vector $q_\Gamma(l_1)$ with coordinates 
\begin{align}
    q_\Gamma^e(l_1):=\begin{cases}
        d_e(\varphi)\frac{(l-l_1)}{l}, & \text{ if } e \mapsto e'\subset \overline{PP_0},\\
        d_e(\varphi)\frac{l_1}{l}, & \text{ if } e \mapsto e'\subset \overline{P_0P},
    \end{cases}
    \text{ for } e\in E(\Gamma).
\end{align}
This allows $ \ker(\varphi^{*})$ to be described in a clean and concise manner.

\end{construction}
\begin{proposition}\label{proposition_kernelofpullback}
Let $\varphi: \Gamma \rightarrow \mathbb{T}E$ be a cover of degree $d$, where $g(\Gamma)=g$, and $q_\Gamma: [0,l] \rightarrow \mathbb{R}^{|E(\Gamma)|}$ the function from Construction \ref{construction_kernelofpullback}. Then $q_\Gamma$ decides whether a divisor lives in the kernel of the pullback map $\varphi^{*}$. We have:
\begin{align}
    \ker(\varphi^{*}) \cong q_\Gamma(\mathbb{T}E) \cap \mathbb{Z}^{|E(\Gamma)|},
\end{align}
where, by abuse of notation, we write $q_\Gamma(\mathbb{T}E)$ using that the domain of $q_\Gamma$ parameterizes points of $\mathbb{T}E$. In particular, if $D_P\in  \ker(\varphi^{*})$, then $ord(D_P)$ divides $d$.
\end{proposition}
\begin{proof}
Note that $q^e_\Gamma(l_1)\in \mathbb{Z}$ is true if and only if $l_1\in l\mathbb{Q}/l\mathbb{Z}$. Hence, the isomorphism is just a consequence of Construction \ref{construction_kernelofpullback}. Suppose $D_P\in  \ker(\varphi^{*})$, i.e. $q_\Gamma(l_1)\in \mathbb{Z}^{|E(\Gamma)|}$ where $l_1$ is the length of  $\overline{PP_0}$. Then $l_1$ is equal to $\frac{i}{ord(D_P)}l$ for an $i\in \{0,...,ord(D_P)-1\}$ and
\begin{align}
    q_\Gamma^e(l_1)=\begin{cases}
        d_e(\varphi)(1-\frac{i}{ord(D_P)}), & \text{ if } e \mapsto e'\subset \overline{PP_0},\\
        d_e(\varphi)\frac{i}{ord(D_P)}, & \text{ if } e \mapsto e'\subset \overline{P_0P},
    \end{cases}
    \in \mathbb{Z} \text{ for all } e\in E(\Gamma)
\end{align}
implies
\begin{align}
    ord(D_P)\thinspace | \thinspace d_e(\varphi) \text{ for all } e\in E(\Gamma).
\end{align}
Considering that $\sum_{e\mapsto e'} d_e(\varphi)=d$ holds for every edge $e'$ of $\mathbb{T}E$, $ord(D_P)$ has to divide $d$ as well. 
%Finally, $ord(D_p) \thinspace | \thinspace d$ follows from the additional observation $\sum_{e\mapsto e'} d_e(\varphi)=d$ for an edge $e'$ of $\mathbb{T}E$.
    \end{proof}

\begin{remark}
    Note, that $\varphi^*$ is injective whenever $q_\Gamma(\mathbb{T}E)$ does not contain any integer vector. 
\end{remark}
\subsection{Criteria for Optimality.}\label{subsection_CriteriaforOptimality}
Verifying optimality based on Definition \ref{definition_optimalmap} is tricky, in general. However, should the combinatorial type of $\Gamma$ be the dumbbell-graph, then the structure of (not necessarily optimal) coverings of $\mathbb{T}E$ is particularly simple: They split into two "independent" isogenies. Lemma \ref{lemma_optimalitycharacterization_DB} gives a characterization of optimality in this case.
\begin{lemma}\label{lemma_optimalitycharacterization_DB}
    Let $ \varphi: \Gamma \rightarrow \mathbb{T}E$ be a degree $d$ cover of dumbell type (Figure \ref{figure_dumbbellcover}). Then $\varphi$ factors through a non-trivial isogeny $\phi$ if and only if $gcd(d_{\varphi}(e_1),d_{\varphi}(e_2))> 1$ or $gcd(n_1,n_2)> 1$ holds.
\end{lemma}
\begin{proof}
For the duration of this proof we will assume $\Gamma$ to have the dumbbell-graph as combinatorial type and work in the notation of Lemma \ref{lemma_kernelpushforwardinTA_DB}. We start with a structural observation: The map
\begin{align}
    \mathcal{F}: M:=\{\varphi: \Gamma \rightarrow \mathbb{T}E \} \rightarrow \{(\varphi_1: \mathbb{T}E_1 \rightarrow \mathbb{T}E, \varphi_2: \mathbb{T}E_2 \rightarrow \mathbb{T}E)\}, \varphi \mapsto (\varphi_{|e_1}, \varphi_{|e_2})
\end{align}
is surjective with fibres $\mathcal{F}^{-1}(\varphi_1,\varphi_2)=\{ \varphi_t: \Gamma_t \rightarrow \mathbb{T}E, \thinspace t\in \mathbb{R}_{\geqslant 0} \}$, where $\Gamma_t$ denotes the curve of type dumbbell obtained by connecting $\mathbb{T}E_i$ along an edge of length $t$. This procedure is not unique, but could be made so by forcing $\mathcal{F}$ to remember the positions of the vertices $P_0,P_1\in \Gamma$. To avoid clutter, we do not mark points. \\
Now that we have identified the building blocks of a $\varphi\in M$, it becomes easy to construct a factorization of $\varphi$: Define isogenies
\begin{align}
    \phi: \mathbb{T}E' \rightarrow \mathbb{T}E\text{ and }    \varphi'_i:\mathbb{T}E_i \rightarrow \mathbb{T}E', \text{ for } i=1,2,
\end{align}
where $\mathbb{T}E'$ is the elliptic curve of length $l(e)$ with edge $e$ and $\mathbb{T}E_i$ the one with edge $e_i$, by 
\begin{align}
   & n(\phi):= gcd(n_1,n_2), d_\phi(e):= gcd(d_{\varphi}(e_1),d_{\varphi}(e_2)), l(e):=\frac{n(\phi)l(\mathbb{T}E)}{d_\phi(e)},\\
  & n_i(\varphi'_i):= \frac{n_i(\varphi)}{n(\phi)}, d_{\varphi'_i}(e_i):=\frac{d_{\varphi}(e_i)}{d_\phi(e)}.
\end{align}
Set $\varphi'\in \mathcal{F}^{-1}(\varphi'_1,\varphi'_2)$ to be a cover whose source curve is $\Gamma$. We have $\varphi=\phi \circ \varphi' $ with $deg(\phi)=n(\phi)\cdot d_\phi(e) $, that is $\phi$ is non-trivial if and only if either $gcd(d_{\varphi}(e_1),d_{\varphi}(e_2))> 1$ or $gcd(n_1,n_2)> 1$ holds. We conclude the proof by observing that this factorization is maximal (in the sense that $deg(\phi)$ is maximal by construction), so that any other factorization of $\varphi$ must factor through this one.

\end{proof}
We could attempt a similar approach for the theta-graph and characterize the absence of optimality by the existence of a non-trivial solution to a system of equations analogous to, but more complex than, the one behind the previous proof. The following Proposition is the choice of a different path.
\begin{proposition}\label{proposition_characterizationofoptimalmaps}
    A cover $\varphi$ is tropically optimal if and only if $\ker(\varphi_{*})$ is connected.
\end{proposition}
\begin{proof}
We know from Lemma \ref{lemma_facorizationgtropicalmaps} that factorizations of $\varphi$ correspond to factorizations of the push-forward $\varphi_{*}$. Thus, consider
\begin{center}
     \begin{tikzpicture}
\matrix(m)[matrix of math nodes,
row sep=3em, column sep=2.8em,
text height=1.5ex, text depth=0.25ex]
{\Jac(\Gamma)& \mathbb{T}\tilde{E} & \mathbb{T}E,\\};
\path[->]
(m-1-1) edge node[above] {$\tilde{\mu}$} (m-1-2)
(m-1-2) edge node[above] {$\phi$} (m-1-3)
(m-1-1) edge [bend right] node[below] {$\varphi_{*}$} (m-1-3)
;
\end{tikzpicture}
    
\end{center}
where $\tilde{\mu}$ is a morphism of tori and $\phi$ an isogeny. We make two observations:
\begin{enumerate}
    \item Without loss of generality we can assume $\ker(\tilde{\mu})$ to be connected: Otherwise, decompose
$\tilde{\mu}$ further as
    \begin{center}
     \begin{tikzpicture}
\matrix(m)[matrix of math nodes,
row sep=3em, column sep=2.8em,
text height=1.5ex, text depth=0.25ex]
{\Jac(\Gamma)& \Jac(\Gamma)/{\Ker(\tilde{\mu})}_0 & \mathbb{T}\tilde{E} & \mathbb{T}E\\};
\path[->]
(m-1-1) edge node[above] {$\tilde{\tilde{\mu}}$} (m-1-2)
(m-1-2) edge node[above] {$\tilde{\phi}$} (m-1-3)
(m-1-3) edge node[above] {$\phi$} (m-1-4)
(m-1-1) edge [bend right] node[below] {$\varphi_{*}$} (m-1-4)
(m-1-1) edge [bend right] node[above] {$\tilde{\mu}$} (m-1-3)
;
\end{tikzpicture}
    
\end{center}
    and replace $\tilde{\mu}$ by $\tilde{\tilde{\mu}}$ and $\phi$ by the isogeny $\tilde{\phi} \circ \phi$.
    \item Assuming (1) (i.e. connectedness of the kernel) we know that ${\ker(\tilde{\mu})}\subset \ker(\varphi_*)$ is a torus. More precisely, it is a subtorus of $\Ker(\varphi_*)_0$ of the same dimension. Since every l-dimensional real torus is simple, we actually have an equality: ${\ker(\tilde{\mu})}= \Ker(\varphi_*)_0$.
\end{enumerate}
   Combining observations (1) and (2) we see that any factorization of $\varphi_*$ factors through the minimal map, $\pi_{\varphi_*}: \Jac(\Gamma) \rightarrow \Jac(\Gamma)/{\Ker(\varphi_*)}_0$ of the Stein factorization of $\varphi_*$ (Lemma \ref{lemma_tropicalSteinfactorization}). Back to the statement of Proposition \ref{proposition_characterizationofoptimalmaps}:\\
    If $\ker(\varphi_*)$ is connected, then $\Jac(\Gamma)/{\Ker(\varphi_*)}_0\cong \mathbb{T}E$ and $\pi_{\varphi_*}=\varphi_*$. In this case, we see that the minimal map is also maximal (Definition \ref{definition_minimalmaximalmaps}).  
    Hence, $\varphi$ must be optimal (again by Lemma \ref{lemma_facorizationgtropicalmaps}).\\
    Otherwise, the map $\Jac(\Gamma) \rightarrow \Jac(\Gamma)/{\Ker(\varphi_*)}_0$ gives rise to a non-trivial factorization of $\varphi_*$ that contradicts optimality. This finishes the proof.
\end{proof}
\begin{remark}\label{remark_minimalmap}
The proof of Proposition \ref{proposition_characterizationofoptimalmaps} yields the following reformulation: A cover $\varphi$ is tropically optimal if and only if the minimal factorization of $\varphi_*$ is also maximal. 
\end{remark}

As a direct consequence we obtain:
\begin{corollary}\label{corollary_anycoverfactorsthroughanoptimalone}
    Any tropical cover $\varphi: \Gamma \rightarrow \mathbb{T}E$, where $g(\Gamma)=2$, factors through an optimal cover.
\end{corollary}

This means: Optimal covers are not uncommon, but an integral part of every cover. While the first map, $\pi_{\varphi_*}$, in the Stein factorization of $\varphi_*$ corresponds to the optimal cover described in Corollary \ref{corollary_anycoverfactorsthroughanoptimalone}, the second map, $\phi_{\varphi_*}$, measures how far $\varphi$ is from being optimal by counting the number of connected components of $\ker(\varphi_*)$.
\section{From Covers to Curves: Tropical split Jacobians}\label{section_tropicalsplitJacobians}
We have seen: In genus $1$, the curve and its Jacobian are (in some sense) indistinguishable and as $1$-dimensional pptav easy to understand. As an increase in genus, however, entails an increase in complexity of the associated Jacobian, it seems natural to look for situations in which it is built of simpler objects.
\begin{definition}\label{definition_splitJacobians}
    Let $\Gamma$ be a tropical curve of genus $2$. We say that $\Jac(\Gamma)$ \emph{splits}, if $\Jac(\Gamma)$ is isogeneous to the coproduct of two elliptic curve $\mathbb{T}E\oplus \mathbb{T}E'$. In this case we call $\mathbb{T}E'$ a \emph{complement} of $\mathbb{T}E$ and vice versa.
\end{definition}

This is a phenomenon in $\mathbb{T}\mathcal{A}$ that takes on the following form in $\mathbb{T}\mathcal{C}$.
\begin{theorem}\label{theorem_Jacsplitsiffcoverexists}
    Let $\Gamma$ be a tropical curve of genus $2$. Then $\Jac(\Gamma)$ splits if and only if $\Gamma$ covers an elliptic curve. 
\end{theorem}
We will leave the proof of Theorem \ref{theorem_Jacsplitsiffcoverexists} to Subsection \ref{subsection_Proof}. At this point we only would like to suggest \emph{why} this characterization might not be satisfactory: We can perturb a splitting of $\Jac(\Gamma)$ by isogenies of elliptic curves, leading to non-uniqueness of the complement or non-uniqueness of the cover, depending on the direction considered in Theorem \ref{theorem_Jacsplitsiffcoverexists}.
\begin{example}
    Suppose $\Jac(\Gamma)$ is isogeneous to $\mathbb{T}E\oplus \mathbb{T}E'$ via $\phi_1$. Let $\mathbb{T}\tilde{E}$ be an elliptic curve of length $3\cdot l_{\mathbb{T}E'}$ and $\phi_2:\mathbb{T}E' \rightarrow \mathbb{T}\tilde{E}$ the corresponding isogeny of degree $3$. Then $(id \oplus \phi_2) \circ\phi: \Jac(\Gamma) \rightarrow \mathbb{T}E\oplus \mathbb{T}\tilde{E}$ is (as composition of isogenies) an isogeny and $\mathbb{T}\tilde{E}$ another complement of $\mathbb{T}E$.
\end{example}
\subsection{The canonical complement}\label{subsection_canonicalcomplement}
Optimal covers improve the situation considerably and will set the scene for Subsections \ref{subsection_canonicalcomplement}, \ref{subsection_complementarycover} and \ref{subsection_fromcoverstocurvesultimately}. The exact setting is described in Theorem \ref{theorem_Jacsplits1} and will be maintained throughout.
\begin{theorem}\label{theorem_Jacsplits1}
If $\varphi: \Gamma \rightarrow \mathbb{T}E$ is an optimal cover of degree $d$, then $\mathbb{T}E':=\ker(\varphi_{*})$ is a tropical elliptic curve and $\Jac(\Gamma)$ splits, i.e. there exists an isogeny $\phi:\mathbb{T}E' \bigoplus \mathbb{T}E \rightarrow \Jac(\Gamma)$ whose kernel is isomorphic to $\Jac_d(\mathbb{T}E)$.
\end{theorem}
\begin{proof}
We start by observing that optimal maps generate short exact sequences: Since $\varphi$ is optimal, we know that $\mathbb{T}E':=\ker(\varphi_*)$ is a tropical abelian variety of dimension 1 (see Lemma \ref{lemma_kernelofpushforwardinTA}) and even more, it is a tropical elliptic curve with principal polarization $\zeta_{\mathbb{T}E'}$ (see Remark \ref{remark_endowingEwithpp}). This means that we have an exact sequence of tropical abelian varieties
    \begin{align}
       0 \rightarrow \mathbb{T}E' \xrightarrow{i} \Jac(\Gamma) \xrightarrow{\varphi_*} \Jac(\mathbb{T}E) \cong \mathbb{T}E \rightarrow 0,   
    \end{align}
    where we used that for tropical elliptic curves the embedding $\mathbb{T}E \hookrightarrow \Jac(\mathbb{T}E), P \mapsto P-Q_0$ for the choice of a point $Q_0\in \mathbb{T}E$ is an isomorphism of rational polyhedral spaces (\cite{röhrle2024tropicalngonalconstruction}, Proposition 4.13). Picking $Q_0$ to be identity element, turns it into a morphism of tori. The dual sequence
    \begin{align}
       0 \rightarrow \widecheck{\mathbb{T}E} \xrightarrow{\widecheck{\varphi_*}} \widecheck{\Jac(\Gamma)} \xrightarrow{\widecheck{i}} \widecheck{\mathbb{T}E'} \rightarrow 0
    \end{align}
    is also exact (see Lemma \ref{lemma_dualizingexactsequences}). Identifying $\widecheck{\mathbb{T}E'}, \widecheck{\Jac(\Gamma)}$ and $\widecheck{\mathbb{T}E}$ with their respective duals yields two exact sequences:
\begin{diagram}\label{diagram_pushandpullseq}
    &0 \rightarrow \mathbb{T}E' \xrightarrow{i} \Jac(\Gamma) \xrightarrow{\varphi_*}  \mathbb{T}E \rightarrow 0  \\
   &   0 \rightarrow \mathbb{T}E \xrightarrow{\varphi^*} \Jac(\Gamma) \xrightarrow{g}  \mathbb{T}E' \rightarrow 0 .
\end{diagram}
Consider the coproduct $\mathbb{T}E' \bigoplus \mathbb{T}E$ in the category of tropical abelian varieties (see Definition \ref{definition_productandcoproducts}) and let $\phi$ be the unique map making the diagram 
 \begin{center} 
    \begin{tikzcd}[row sep=huge]
        & \Jac(\Gamma) & \\
        \mathbb{T}E' \ar[ur,"i",sloped] \ar[r, swap] & \mathbb{T}E' \bigoplus \mathbb{T}E \ar[u,dashed,"{\phi}" description] &\mathbb{T}E \ar[ul,"\varphi^{*}",sloped] \ar[l]
    \end{tikzcd}
         
 \end{center}
    commute. We claim that the morphism $\phi$  which is defined by $(P',P)\mapsto i(P')+ \varphi^{*}(P)$ is an isogeny, i.e. surjective and finite. From the two sequences (\ref{diagram_pushandpullseq}) we see that $\Jac(\Gamma) $ is the sum (not the direct sum!) of two elliptic curves:
    \begin{align}
        \Jac(\Gamma) \cong \Jac(\Gamma)/\ker(\varphi_*) + \ker(\varphi_*) \cong \im(\varphi^{*}) + \im(i).
    \end{align}
This yields another exact sequence (now of abelian groups!):
\begin{diagram}\label{diagram_intheoremJacsplits_sequenceforphisurjective}
    &0 \rightarrow \im(\varphi^{*}) \cap \im(i) \xrightarrow{f_1}  \im(\varphi^{*}) \bigoplus \im(i) \xrightarrow{f_2} \Jac(\Gamma) \rightarrow 0,
\end{diagram}
where the map $f_1$ is given by $x \mapsto (x,-x)$ and $f_2$ is the usual addition. Recalling that $i$ and $\varphi^{*}$ are injective, identifies sequence (\ref{diagram_intheoremJacsplits_sequenceforphisurjective}) with the exact sequence associated to $\phi$  (as group homomorphism in the category of abelian groups). This proves surjectivity. Finally, since 
\begin{align}
        \ker(\phi) \cong \im(\varphi^{*}) \cap \im(i) = \im(\varphi^{*}) \cap \ker(\varphi_{*}) \cong \ker (\varphi_{*}\varphi^{*})=\Jac_d(\mathbb{T}E)
    \end{align} is finite, $\phi$ is finite.
\end{proof}
To avoid working exclusively in the abstract category of tav, we translate the proof into the explicit language of matrices.
\begin{example}\label{example_isogenyarisingfromoptimalcover}
    Let $\varphi$ be the optimal cover from Example \ref{example_optimalandnotoptimal}. Recall that the kernel of $\varphi_{*}$ is an elliptic curve of length $1$ with lattice representation given in terms of the base choice from Convention \ref{convention_genus2coversofE} (see also Lemma \ref{lemma_kernelofpushforwardinTA})
    \begin{align}
        (\langle\omega:= \omega_1-\omega_2 \rangle, \thinspace \langle B:= B_2 \rangle,\thinspace [\cdot,\cdot]_K).
    \end{align}
    Set $\mathbb{T}E':=\ker(\varphi_{*})$. We already have a coordinate description of $\varphi_{*}$, let us now determine one for the natural inclusion $i: \mathbb{T}E' \hookrightarrow \Jac(\Gamma)$: As a map on the underlying tori, $i$ is induced by 
    \begin{align}
        \Hom(i^\#): \Hom(\langle\omega\rangle, \mathbb{R}) \rightarrow \Hom(\langle\omega_1,\omega_2\rangle, \mathbb{R}),
    \end{align}
    where $i^\#: \langle\omega_1,\omega_2\rangle \twoheadrightarrow \langle\omega_1,\omega_2\rangle/\langle2\omega_1-\omega_2\rangle \cong\langle\omega\rangle$ is the quotient map. For $f \in \Hom(\langle\omega\rangle, \mathbb{R})$ we have $f\circ i^\#(\omega_1)= f(\omega) $ and $f\circ i^\#(\omega_2)= 2f(\omega) $ (since $\omega_1= \omega_1 -(2\omega_1-\omega_2)= \omega$ in the quotient, analogously for $\omega_2$), which makes the first exact sequence in (\ref{diagram_pushandpullseq}) of the preceding proof explicit:
    \begin{diagram}[ampersand replacement=\&]
     % \begin{tikzcd}
 0 \arrow[r] 
\&  \mathbb{T}E' \arrow[r, "i"] \arrow[d, "\cong"] \& \Jac(\Gamma) \arrow[d, "\cong "] \arrow[r, "\varphi_* "] \& \mathbb{T}E \arrow[d, "\cong "] \arrow[r] \& 0 \\
 0 \arrow[r]  \&  \mathbb{R}/\mathbb{Z} \arrow[r, "{\begin{pmatrix}
     1\\
     2
 \end{pmatrix}}"]\& \mathbb{R}^2/ \begin{pmatrix}
     2 & 1\\
     1 & 2
 \end{pmatrix}\mathbb{Z}^2 \arrow[r,, "{\begin{pmatrix}
     2 & -1
 \end{pmatrix}}"] \& \mathbb{R}/3\mathbb{Z} \arrow[r]  \&  0.
%\end{tikzcd}  
\end{diagram}
  We proceed in the same way for the second and obtain
  \begin{diagram}[ampersand replacement=\&]
     % \begin{tikzcd}
 0 \arrow[r] 
\&  \mathbb{T}E \arrow[r, "\varphi^*"] \arrow[d, "\cong"] \& \Jac(\Gamma) \arrow[d, "\cong "] \arrow[r, "g "] \& \mathbb{T}E' \arrow[d, "\cong "] \arrow[r] \& 0 \\
 0 \arrow[r]  \&  \mathbb{R}/3\mathbb{Z} \arrow[r, "{\begin{pmatrix}
     1\\
     0
 \end{pmatrix}}"] \& \mathbb{R}^2/ \begin{pmatrix}
     2 & 1\\
     1 & 2
 \end{pmatrix}\mathbb{Z}^2 \arrow[r,, "{\begin{pmatrix}
     0 & 1
 \end{pmatrix}}"] \& \mathbb{R}/\mathbb{Z} \arrow[r]  \&  0 .
%\end{tikzcd}  
  \end{diagram} 
  By joining the left parts of the both sequences we obtain a coordinate description for the isogeny $\phi$:
  \begin{diagram}[ampersand replacement=\&]
     % \begin{tikzcd}
\mathbb{T}E' \oplus \mathbb{T}E \arrow[rr, "\phi"] \arrow[d, "\cong"] \& \& \Jac(\Gamma) \arrow[d, "\cong "]  \\
 \mathbb{R}/\mathbb{Z}  \oplus \mathbb{R}/3\mathbb{Z} \arrow[rr, "{\begin{pmatrix}
     1\thinspace & 1\\
     2\thinspace & 0
 \end{pmatrix}}"] \& \& \mathbb{R}^2/ \begin{pmatrix}
     2 & 1\\
     1 & 2
 \end{pmatrix}\mathbb{Z}^2 
%\end{tikzcd}  
  \end{diagram} 
  whose kernel $\ker(\phi)=\{\begin{pmatrix}
     0 \\
     0 
 \end{pmatrix}, \begin{pmatrix}
     \frac{1}{2} \\
     \frac{3}{2}
 \end{pmatrix} \} \cong \mathbb{R}/3\mathbb{Z}[2]$.
  
\end{example}
\subsection{The complementary cover}\label{subsection_complementarycover}
The proof of Theorem \ref{theorem_Jacsplits1} reveals a certain kind of symmetry: In terms of exact sequences $\mathbb{T}E'$ is related to $\Jac(\Gamma)$ in the same way as $\mathbb{T}E$ is. Hence, it is natural to ask, whether this symmetry does not hold at a more fundamental level, i.e. whether the sequences 
\begin{diagram}\label{diagram_subsectioncomplcover:pushandpullseq}
    &0 \rightarrow \mathbb{T}E' \xrightarrow{i} \Jac(\Gamma) \xrightarrow{\varphi_*}  \mathbb{T}E \rightarrow 0  \\
   &   0 \rightarrow \mathbb{T}E \xrightarrow{\varphi^*} \Jac(\Gamma) \xrightarrow{g}  \mathbb{T}E' \rightarrow 0 
\end{diagram}
are generated by an optimal map as well. To this end, consider:
\begin{construction}\label{construction_complementarycover}
Let $\Gamma$ be of genus $2$ with labeling as fixed in Figure \ref{figure_generalgenus2cover} and Figure \ref{figure_dumbbellcover}, depending on whether $\Gamma$ is of type theta or of type dumbbell. Let us define $\varphi': \Gamma \rightarrow \mathbb{T}E'$ as the composition 
   \begin{align}
     \Gamma  \xrightarrow{\Phi_{P_0}} \Jac(\Gamma) \overset{f_{\zeta_\Gamma}}{\cong}\widecheck{\Jac(\Gamma)} \xrightarrow{\widecheck{i}}   \widecheck{\mathbb{T}E'} \overset{f_{\zeta_{\mathbb{T}E'}}}{\cong} \mathbb{T}E',
   \end{align} 
    where $\Phi_{P_0}$ denotes the Abel-Jacobi map with base point $P_0\in \Gamma$.
\end{construction}
The following Lemma paves the way for a reinterpretation of the sequences (\ref{diagram_subsectioncomplcover:pushandpullseq}) in terms of $\varphi'$.
\begin{lemma}\label{lemma_thecomplementarycover}
    The map $\varphi': \Gamma \rightarrow \mathbb{T}E'$ from Construction \ref{construction_complementarycover} is a tropical cover. 
\end{lemma}
\begin{proof}[Proof of Lemma \ref{lemma_thecomplementarycover}.]
We prove the statement in two steps. In step 1 we deal with $\mathbb{T}E'$ in its incarnation as the tav $\ker(\varphi_*)$ whose lattice representation is (see Definition \ref{definition_(Co-)KernelImage})
\begin{align}
    (K,K',[\cdot,\cdot]_K):=(\Omega^1_\Gamma/\im({\varphi_*}^\#)^{\text{sat}},\ker({\varphi_*}_\#),[\cdot,\cdot]_K)
\end{align}
and show that $\varphi'$ (now viewed as map to a tropical torus) is \emph{tropical}, i.e. satisfies the following properties:
    \begin{enumerate}
        \item $\varphi'$ is continuous and restricts to an affine function on each edge.
        \item the differential $D\varphi_P'$ at $P\in \Gamma$ takes integral tangent vectors to integral tangent vectors.
        \item $\varphi'$ is harmonic.
    \end{enumerate}
    For step 2, we turn to $\mathbb{T}E'$ as metric graph and recognize in $\varphi'$ an actual morphism of graphs, a tropical cover.\\
  Note that the first part of (1), $\varphi'$ is continuous, is an immediate consequence of the continuity of the Abel-Jacobi map $\Phi_{P_0}$ (Theorem 4.1, \cite{MR2772537}) since $f_{\zeta_\Gamma},f_{\zeta_{\mathbb{T}E'}}$ and $\widecheck{i}$ are obviously continuous as morphisms. For the second part we utilize the local description of the lift $\tilde{\Phi}_{P_0}$ of the Abel-Jacobi map to $\Hom(\Omega^1_\Gamma,\mathbb{R})$ given in \cite{MR2772537} to determine a description of the lift $\tilde{\varphi}'$ of $\varphi'$ to $\Hom(K,\mathbb{R})$. That is the map given by the upper path in the following diagram: 

   %\begin{diagram}
      \begin{tikzcd}[scale cd=0.85]
 \Gamma \arrow[r, "\widetilde{\Phi}_{P_0} "] \arrow[dr, "\Phi_{P_0} "] 
&  \Hom(\Omega^1_\Gamma,\mathbb{R}) \arrow[r, "\Hom(\zeta_\Gamma) "] \arrow[d, "\pi_1 "] & \Hom(\H1(\Gamma, \mathbb{Z}),\mathbb{R}) \arrow[d, "\pi_2 "] \arrow[r, "\Hom(i_\#) "] & \Hom(K',\mathbb{R}) \arrow[d, "\pi_3 "] \arrow[r, "\Hom(\zeta_{\mathbb{T}E'}) "] & \Hom(K,\mathbb{R}) \arrow[d, "\pi_4 "] \\
& \Jac(\Gamma) \arrow[r,"f_{\zeta_\Gamma}"]  &  \widecheck{\Jac}(\Gamma) \arrow[r,"\widecheck{i}"] &  \widecheck{\ker}(\varphi_*) \arrow[r,"f_{\zeta_{\mathbb{T}E'}}"]  &  \ker({\varphi_*}) \\
\end{tikzcd}  
 % \end{diagram}

  where $\pi_i$ for $i=1,...,4$ are the canonical projections. Let $\hat{e}$ be an oriented edge with parametrization 
  \begin{align}
      [0,l(\hat{e})] \ni t \mapsto S + t\cdot v_{\hat{e}},
  \end{align}
 where $S$ is the initial vertex and $v_{\hat{e}}$ a primitiv integral tangent vector at some interior point $P\in \hat{e} $. Then $\widetilde{\Phi}_{P_0} $ is given by the affine map
 \begin{align}
      [0,l(\hat{e})] \ni t \mapsto \widetilde{\Phi}_{P_0}(S) + t\cdot (\frac{1}{l(\hat{e})}\int_{\hat{e}})\in \Hom(\Omega^1_\Gamma,\mathbb{R})
  \end{align}
  and $\tilde{\varphi}'$ by
  \begin{align}
      [0,l(\hat{e})] \ni t \mapsto \widetilde{\Phi}_{P_0}(S)\circ \zeta_\Gamma \circ i_{\#}\circ \zeta_{\mathbb{T}E'} + t \cdot (\frac{1}{l(\hat{e})}\int_{\hat{e}} \zeta_\Gamma \circ i_{\#}\circ \zeta_{\mathbb{T}E'})\in \Hom(K,\mathbb{R}).
  \end{align}
  We (need to) show that (the image of) $(\frac{1}{l(\hat{e})}\int_{\hat{e}} \zeta_\Gamma \circ i_{\#}\circ \zeta_{\mathbb{T}E'})$ (under $\pi_4$) is an integral tangent vector of $\ker(\varphi_*)$, i.e. (that the lift is) an element of the dual lattice $\Hom(K,\mathbb{Z})$. Since $\zeta_{\mathbb{T}E'}:K' \rightarrow K$ is an isomorphism, it suffices to verify that
  \begin{align}
      \frac{1}{l(\hat{e})}\int_{\hat{e}} \zeta_\Gamma \circ i_{\#}(\cdot) \in \Hom(K',\mathbb{Z})
  \end{align}
  holds. Recalling that $K'=\langle B \rangle_\mathbb{Z}$ with $B= w_1B_1 +w_2B_2$ and $w_1,w_2\in \mathbb{Z}$ yields:
  \begin{diagram}[ampersand replacement = \&]\label{diagram_inproofofthecomplentarycover}
      \frac{1}{l(\hat{e})}\int_{\hat{e}} \zeta_\Gamma \circ i_{\#}(B)=\frac{1}{l(\hat{e})}\int_{\hat{e}}(w_1(de+de_2)+w_2(de_2-de_1)\\
      =\begin{cases}
          w_1, & \hat{e}=e\\
          -w_2, & \hat{e}=e_1\\
          w_1+w_2, & \hat{e}=e_2
      \end{cases},
  \end{diagram}
  whenever $\Gamma$ is a theta-graph (Lemma \ref{lemma_kernelofpushforwardinTA} and notation therein), and
   \begin{diagram}[ampersand replacement = \&]\label{diagram_inproofofthecomplentarycover_DB}
      \frac{1}{l(\hat{e})}\int_{\hat{e}} \zeta_\Gamma \circ i_{\#}(B)=\frac{1}{l(\hat{e})}\int_{\hat{e}}(w_1(de_1)+w_2(de_2)\\
      =\begin{cases}
          w_1, & \hat{e}=e_1\\
          w_2, & \hat{e}=e_2\\
         0, & \hat{e}=e
      \end{cases},
  \end{diagram}
   whenever $\Gamma$ is a dumbbell-graph  (Lemma \ref{lemma_kernelpushforwardinTA_DB} and notation therein). The maps $\frac{1}{l(\hat{e})}\int_{\hat{e}} \zeta_\Gamma \circ i_{\#}(\cdot)$ for $\hat{e} \in E(\Gamma)$ take integer values on a lattice base of $K'$ and are therefore elements of $\Hom(K',\mathbb{Z})$. This proves (1) and (2).\\
  Computations (\ref{diagram_inproofofthecomplentarycover}) and (\ref{diagram_inproofofthecomplentarycover_DB}) also provide a coordinate description of the differential
  \begin{align}
      D\varphi_P': T_P(\Gamma) \rightarrow T_{\varphi'(P)}(\ker(\varphi_*))
  \end{align} of $\varphi'$ at an interior point $P$ of $\hat{e}$. Consider, for example, (\ref{diagram_inproofofthecomplentarycover}): A generator $v_{\hat{e}}$ of $ T_P(\Gamma)$ is sent to $w_1,-w_2$ or $w_1+w_2$, depending on whether $ \hat{e}$ is (equal to) $e,e_1$ or $e_2$. With this one easily checks that $\varphi'$ is balanced at the vertices of $\Gamma$, exemplarily for $P_1$ we have:
  \begin{align}
      \sum_\lambda D\varphi_{P_1}'(\lambda)= w_1 + w_2- (w_1+w_2)=0,
  \end{align}
  where the first sum is over all primitive integral tangent vectors $\lambda$ at $P_1$, that is $\lambda=v_{e}, -v_{e_1},-v_{e_2}$. The verification of the balancing condition (using (\ref{diagram_inproofofthecomplentarycover_DB})), when $\Gamma$ is a dumbbell-graph, is analogous. \\
  For step 2 we turn $\varphi'$ into a morphism of metric graphs by viewing $\ker(\varphi_*')$ as elliptic curve $\mathbb{T}E'$ with vertex set $V(\mathbb{T}E'):=\{\varphi'(P_0),\varphi'(P_1)\}$. Since $\varphi'$ is non-constant and harmonic at all vertices of $\Gamma$, it must be surjective (see \cite{MR3375652}, Remark 2.7) and therefore a tropical cover.
\end{proof}
\begin{remark}
    Note that $\varphi'$ is not necessarily finite (see Example \ref{example_complementarycover}).
\end{remark}
Lemma \ref{lemma_thecomplementarycover} provides us with a second cover. From this perspective Theorem \ref{theorem_Jacsplits1} reads as follows:
\begin{theorem}\label{theorem_Jacsplits2}
In the setting of Theorem \ref{theorem_Jacsplits1}, there exists another cover 
 $\varphi': \Gamma \rightarrow \mathbb{T}E'$ such that 
 \begin{itemize}
     \item $\varphi'$ also generates the exact sequences (\ref{diagram_subsectioncomplcover:pushandpullseq}), i.e. $i=\varphi'^*$ and $g=\varphi'_*$. 
     \item $\varphi'$ is optimal and of degree $d$.
 \end{itemize}
 Moreover, the isogeny $\phi:\mathbb{T}E' \bigoplus \mathbb{T}E \rightarrow \Jac(\Gamma)$ satisfies $\ker(\phi)\cong Jac_d(\mathbb{T}E')$ and is polarized with respect to $f_{\zeta_\Gamma}$ and $m_d\circ f_{\zeta_{\mathbb{T}E' \bigoplus \mathbb{T}E}}$ (see Definition \ref{definition_induced/pf/pbpolarizationandpolarizedisogeny}), where $m_d: \mathbb{T}E' \bigoplus \mathbb{T}E \rightarrow \mathbb{T}E' \bigoplus \mathbb{T}E$ denotes the multiplication-by-$d$ map and $\gls{zetaGamma}$, respectively $\zeta_{\mathbb{T}E' \bigoplus \mathbb{T}E}$, are the usual principal polarization on $\Jac(\Gamma)$ and $\mathbb{T}E' \bigoplus \mathbb{T}E$ (see Definition \ref{definition_productandcoproducts}).
\end{theorem}
\begin{proof}
Let $\varphi'$ be the cover from Construction \ref{construction_complementarycover}. The notions of push-forward and pull-back morphisms that naturally come with it, make it worthwhile to re-examine the exact sequences (\ref{diagram_subsectioncomplcover:pushandpullseq}). 
First, note that their labelling with $\varphi^*$ and $\varphi_*$ is not quite precise: We should actually replace $\mathbb{T}E$ with $\Jac(\mathbb{T}E)$, but have omitted the isomorphism $\mathbb{T}E \cong \Jac(\mathbb{T}E)$ to avoid cluttering notation. The same applies to the desired reinterpretation of $i$ and $g$ in terms $\varphi'$. We specify an embedding of $\mathbb{T}E'$ into its Jacobian via
\begin{align}
    j: \mathbb{T}E' \rightarrow \Jac(\mathbb{T}E'),\thinspace P \mapsto P- 0_{\Jac(\mathbb{T}E')},
\end{align}
where $0_{\Jac(\mathbb{T}E')}$ is the identity element. Again, by abuse of notation we call $g$ the map between the respective Jacobians to get a commutative diagram
 \begin{diagram}
     % \begin{tikzcd}
 \Gamma \arrow[r, "{\Phi}_{P_0} "] \arrow[d, "\varphi' "] 
&  \Jac(\Gamma) \arrow[d, "g"] \\
\mathbb{T}E' \arrow[r, hook ,"j"] & \Jac(\mathbb{T}E'),   \\
%\end{tikzcd}  
  \end{diagram}
  which proves $g=\varphi'_*$. Turning towards Lemma \ref{lemma_pushandpullaredual} next, provides us with the right interpretation for $i$ as pull-back of $\varphi'$ and seeing that $\ker(\varphi'_*)$ is connected (exactness of the second sequence in (\ref{diagram_subsectioncomplcover:pushandpullseq})), ultimately shows that $\varphi'$ is optimal (Lemma \ref{proposition_characterizationofoptimalmaps}).\\
  As tropical cover $\varphi'$ satisfies
  \begin{align}
      \varphi'_{*}\circ {\varphi'}^*= deg(\varphi')id.
  \end{align}
 Then $deg(\varphi')$ can be computed as $|\ker(\varphi'_{*}\circ \varphi^{'*})|$. Using the reinterpretation of $i$ and $g$ in terms of $\varphi'$ to interchange the roles of $\varphi$ and $\varphi'$ in the proof of Theorem \ref{theorem_Jacsplits1}, we obtain $\ker(\phi)\cong \ker(\varphi'_{*}\circ \varphi^{'*}) $ and thus $deg(\varphi')=d$ and $\ker(\phi)\cong \Jac_d(\mathbb{T}E') $, simultaneously.\\
We conclude by investigating the behaviour of $\phi$ with respect to the pp $\gls{zetaGamma}$ and $\zeta_{\mathbb{T}E' \bigoplus \mathbb{T}E}$: Since $\phi$ is finite, it induces a polarization $\phi^*(\zeta_\Gamma)$ on $\mathbb{T}E' \bigoplus \mathbb{T}E$ that makes the diagram 
\begin{diagram}
    %\begin{tikzcd}
\H1(E', \mathbb{Z}) \bigoplus \H1(E, \mathbb{Z}) \arrow[r, "\phi_{\#} "] \arrow[d, "\phi^*(\zeta_\Gamma)"]
& \H1(\Gamma, \mathbb{Z}) \arrow[d, "\zeta_\Gamma"] \\
\Omega^1_{E'} (\mathbb{Z}) \bigoplus \Omega^1_{E} (\mathbb{Z})
& \arrow[l,"\phi^{\#}"] \Omega^1_{\Gamma}(\mathbb{Z})
%\end{tikzcd}
\end{diagram}
commute. Using that
\begin{itemize}
    \item $\phi_{\#} = \varphi'_{*,\#} \oplus \varphi_{*,\#}$ is the composition of $(\varphi'_{*,\#} , \varphi_{*,\#})$ with the addition-of-components map,
    \item $\phi^{\#}=(\varphi'^{*,\#}, \varphi^{*,\#})$,
    \item $\varphi'_{*}\circ \varphi^{'*}$ and $\varphi^{*} \circ \varphi_{*}$ are multiplication-by-$d$ maps on $\mathbb{T}E' $ and $ \mathbb{T}E$,
\end{itemize} 
together with exactness of (\ref{diagram_subsectioncomplcover:pushandpullseq}), yields
\begin{align}
    \phi^*(\zeta_\Gamma)(\sum_{e'}a_{e'} e',\sum_{e}a_{e} e)= d\cdot (\sum_{e'}a_{e'} de',\sum_{e} a_e de)
\end{align}
as claimed.
\end{proof}
\begin{example}\label{example_complementarycover}
Let us retrace Subsection \ref{subsection_complementarycover} by continuing with Example \ref{example_isogenyarisingfromoptimalcover}. The natural principal polarization $\zeta_\Gamma$ provides us with an identification of $\Jac(\Gamma)$ with its dual. Due to our choice of basis for $\H1(\Gamma,\mathbb{Z})$ and $\Omega^1_\Gamma$ (Convention \ref{convention_genus2coversofE}), its representation matrix is simply $I_2$, analogously for $f_{\zeta_{\mathbb{T}E'}}$. Then $\varphi'$ (as a tropical map, here we view its target as a torus and not as a curve) is given by the composition:
\begin{diagram}[ampersand replacement = \&]
     % \begin{tikzcd}
 \Gamma \arrow[r, hook, "\Phi_{P_0} "] \& \Jac(\Gamma) \arrow[r,"f_{\zeta_\Gamma}"] \arrow[d, "\cong "] \&  \widecheck{\Jac}(\Gamma) \arrow[r,"\widecheck{i}"] \arrow[d, "\cong "] \&  \widecheck{\mathbb{T}E'} \arrow[r,"f_{\zeta_{\mathbb{T}E'}}"] \arrow[d, "\cong "]  \& \mathbb{T}E' \arrow[d, "\cong ."] \\
\& \mathbb{R}^2/ \begin{pmatrix}
     2 & 1\\
     1 & 2
 \end{pmatrix}\mathbb{Z}^2  \arrow[r,"I_2"]  \&  \mathbb{R}^2/ \begin{pmatrix}
     2 & 1\\
     1 & 2
 \end{pmatrix}\mathbb{Z}^2  \arrow[r,"{\begin{psmallmatrix}
     0 & 1
 \end{psmallmatrix}}"] \&  \mathbb{R}/\mathbb{Z} \arrow[r,"1"]  \&  \mathbb{R}/\mathbb{Z}
%\end{tikzcd}  
  \end{diagram}
  We want to understand $\varphi'$ as harmonic morphism of metric graphs next. By fixing $P_0$ as in Figure \ref{figure_phi'astropicalcandgraphc} (on the right), we obtain a representation of $\Gamma$ inside its Jacobian (see Figure \ref{figure_phi'astropicalcandgraphc} on the left). Equipped with coordinates for the points of $\Gamma$ we can simply compute $\varphi'(V(\Gamma))$ by
  \begin{align}
      \varphi'(P_0)=0 \text{ and } \varphi'(P_1)=\begin{pmatrix}
          0 & 1
      \end{pmatrix} \begin{pmatrix}
          1\\
          1
      \end{pmatrix}= 0\in \mathbb{R}/\mathbb{Z}
  \end{align} and set $V(\mathbb{T}E'):=\{P:=\varphi'(P_0)\}$. To determine the local descriptions of $\varphi'$ we proceed as in the proof of Lemma \ref{lemma_thecomplementarycover}. Since $B=B_2$, we have:
 \begin{align}
     \varphi'_{|e}: [0,1] \ni t \mapsto P , \varphi'_{|e_1}: [0,1] \ni t \mapsto P -t \text{ and }\varphi'_{|e_2}: [0,1] \ni t \mapsto P +t.
 \end{align}
 These encode \emph{how} $\Gamma$ covers $\mathbb{T}E'$, so we see that $\varphi'$
 \begin{itemize}
     \item contracts $e$, i.e. that $d_e(\varphi')=0$.
     \item maps $e_i$ to $e$ without any dilation (but reversing the orientation of $e_1$).
     \item is of degree $2$.
 \end{itemize}
 \begin{figure}[H]
     \centering
    \begin{tikzpicture}[x=0.75pt,y=0.75pt,yscale=-0.85,xscale=0.85]

%Straight Lines [id:da18000378658315386] 
\draw  [dash pattern={on 0.84pt off 2.51pt}]  (51.4,250.12) -- (201.33,173.33) ;
\draw [shift={(51.4,250.12)}, rotate = 332.88] [color={rgb, 255:red, 0; green, 0; blue, 0 }  ][fill={rgb, 255:red, 0; green, 0; blue, 0 }  ][line width=0.75]      (0, 0) circle [x radius= 3.35, y radius= 3.35]   ;
%Straight Lines [id:da701595484295165] 
\draw [color={rgb, 255:red, 145; green, 217; blue, 64 }  ,draw opacity=1 ]   (128.33,98.33) -- (51.4,250.12) ;
%Straight Lines [id:da4788926723840894] 
\draw  [dash pattern={on 0.84pt off 2.51pt}]  (128.33,98.33) -- (275.33,23.33) ;
\draw [shift={(275.33,23.33)}, rotate = 332.97] [color={rgb, 255:red, 0; green, 0; blue, 0 }  ][fill={rgb, 255:red, 0; green, 0; blue, 0 }  ][line width=0.75]      (0, 0) circle [x radius= 3.35, y radius= 3.35]   ;
\draw [shift={(128.33,98.33)}, rotate = 332.97] [color={rgb, 255:red, 0; green, 0; blue, 0 }  ][fill={rgb, 255:red, 0; green, 0; blue, 0 }  ][line width=0.75]      (0, 0) circle [x radius= 3.35, y radius= 3.35]   ;
%Straight Lines [id:da9746966587966581] 
\draw  [dash pattern={on 0.84pt off 2.51pt}]  (275.33,23.33) -- (201.33,173.33) ;
%Straight Lines [id:da7840168828408347] 
\draw    (128.33,173.33) -- (51.4,250.12) ;
\draw [shift={(94.82,206.78)}, rotate = 135.05] [color={rgb, 255:red, 0; green, 0; blue, 0 }  ][line width=0.75]    (10.93,-4.9) .. controls (6.95,-2.3) and (3.31,-0.67) .. (0,0) .. controls (3.31,0.67) and (6.95,2.3) .. (10.93,4.9)   ;
\draw [shift={(128.33,173.33)}, rotate = 135.05] [color={rgb, 255:red, 0; green, 0; blue, 0 }  ][fill={rgb, 255:red, 0; green, 0; blue, 0 }  ][line width=0.75]      (0, 0) circle [x radius= 3.35, y radius= 3.35]   ;
%Straight Lines [id:da45837369766206004] 
\draw    (128.33,98.33) -- (128.33,173.33) ;
\draw [shift={(128.33,141.83)}, rotate = 270] [color={rgb, 255:red, 0; green, 0; blue, 0 }  ][line width=0.75]    (10.93,-4.9) .. controls (6.95,-2.3) and (3.31,-0.67) .. (0,0) .. controls (3.31,0.67) and (6.95,2.3) .. (10.93,4.9)   ;
%Straight Lines [id:da2069419068736087] 
\draw    (201.33,173.33) -- (128.33,173.33) ;
\draw [shift={(171.83,173.33)}, rotate = 180] [color={rgb, 255:red, 0; green, 0; blue, 0 }  ][line width=0.75]    (10.93,-4.9) .. controls (6.95,-2.3) and (3.31,-0.67) .. (0,0) .. controls (3.31,0.67) and (6.95,2.3) .. (10.93,4.9)   ;
\draw [shift={(201.33,173.33)}, rotate = 180] [color={rgb, 255:red, 0; green, 0; blue, 0 }  ][fill={rgb, 255:red, 0; green, 0; blue, 0 }  ][line width=0.75]      (0, 0) circle [x radius= 3.35, y radius= 3.35]   ;
%Shape: Axis 2D [id:dp05226806402163109] 
\draw  (23.52,250.12) -- (302.33,250.12)(51.4,5.33) -- (51.4,277.32) (295.33,245.12) -- (302.33,250.12) -- (295.33,255.12) (46.4,12.33) -- (51.4,5.33) -- (56.4,12.33) (126.4,245.12) -- (126.4,255.12)(201.4,245.12) -- (201.4,255.12)(276.4,245.12) -- (276.4,255.12)(46.4,175.12) -- (56.4,175.12)(46.4,100.12) -- (56.4,100.12)(46.4,25.12) -- (56.4,25.12) ;
\draw   (133.4,262.12) node[anchor=east, scale=0.75]{1} (208.4,262.12) node[anchor=east, scale=0.75]{2} (283.4,262.12) node[anchor=east, scale=0.75]{3} (48.4,175.12) node[anchor=east, scale=0.75]{1} (48.4,100.12) node[anchor=east, scale=0.75]{2} (48.4,25.12) node[anchor=east, scale=0.75]{3} ;
%Shape: Ellipse [id:dp14336277412102638] 
\draw  [color={rgb, 255:red, 126; green, 211; blue, 33 }  ,draw opacity=1 ] (353.44,237.09) .. controls (353.46,228.95) and (384.86,222.36) .. (423.57,222.37) .. controls (462.28,222.38) and (493.64,229) .. (493.62,237.14) .. controls (493.6,245.29) and (462.2,251.88) .. (423.49,251.87) .. controls (384.78,251.85) and (353.41,245.24) .. (353.44,237.09) -- cycle ;
%Straight Lines [id:da8921489481323863] 
\draw [color={rgb, 255:red, 126; green, 211; blue, 33 }  ,draw opacity=1 ]   (353.44,237.09) ;
\draw [shift={(353.44,237.09)}, rotate = 0] [color={rgb, 255:red, 126; green, 211; blue, 33 }  ,draw opacity=1 ][fill={rgb, 255:red, 126; green, 211; blue, 33 }  ,fill opacity=1 ][line width=0.75]      (0, 0) circle [x radius= 3.35, y radius= 3.35]   ;
\draw  [color={rgb, 255:red, 126; green, 211; blue, 33 }  ,draw opacity=1 ] (415.05,248.25) .. controls (422.39,250.26) and (429.73,251.47) .. (437.07,251.87) .. controls (429.73,252.27) and (422.39,253.48) .. (415.05,255.49) ;
\draw   (408.19,99.9) .. controls (414.73,99.41) and (421.16,99.57) .. (427.47,100.36) .. controls (421.28,98.93) and (415.2,96.85) .. (409.25,94.12) ;
%Straight Lines [id:da057219205946186236] 
\draw    (351.65,65) -- (351.65,134.72) ;
%Curve Lines [id:da11667775833077876] 
\draw    (351.65,65) .. controls (378.44,44.99) and (500.76,65.29) .. (494.51,102.07) ;
%Curve Lines [id:da05397321221656548] 
\draw    (351.65,134.72) .. controls (369.51,163.85) and (486.48,144.72) .. (494.51,102.07) ;
%Curve Lines [id:da9129942626283467] 
\draw    (351.65,65) .. controls (349.87,102.95) and (509.69,105.6) .. (467.73,107.36) ;
%Curve Lines [id:da5076202627312301] 
\draw    (351.65,134.72) .. controls (365.94,108.24) and (448.98,107.36) .. (467.73,107.36) ;
\draw   (409.69,144.72) .. controls (416.57,145.5) and (423.3,145.49) .. (429.89,144.69) .. controls (423.45,146.28) and (417.16,148.65) .. (411.02,151.81) ;
%Straight Lines [id:da4811503087969162] 
\draw    (423.08,160.61) -- (423.65,203.62) ;
\draw [shift={(423.68,205.62)}, rotate = 269.24] [color={rgb, 255:red, 0; green, 0; blue, 0 }  ][line width=0.75]    (10.93,-3.29) .. controls (6.95,-1.4) and (3.31,-0.3) .. (0,0) .. controls (3.31,0.3) and (6.95,1.4) .. (10.93,3.29)   ;
%Straight Lines [id:da32740956847029334] 
\draw    (351.65,134.72) ;
\draw [shift={(351.65,134.72)}, rotate = 0] [color={rgb, 255:red, 0; green, 0; blue, 0 }  ][fill={rgb, 255:red, 0; green, 0; blue, 0 }  ][line width=0.75]      (0, 0) circle [x radius= 3.35, y radius= 3.35]   ;
%Straight Lines [id:da6153634439503337] 
\draw    (351.65,65) ;
\draw [shift={(351.65,65)}, rotate = 0] [color={rgb, 255:red, 0; green, 0; blue, 0 }  ][fill={rgb, 255:red, 0; green, 0; blue, 0 }  ][line width=0.75]      (0, 0) circle [x radius= 3.35, y radius= 3.35]   ;
\draw   (348.08,108.54) .. controls (350.04,103.21) and (351.23,97.88) .. (351.66,92.53) .. controls (351.99,97.88) and (353.1,103.23) .. (354.96,108.6) ;

% Text Node
\draw (-1.23,253.28) node [anchor=north west][inner sep=0.75pt]    {$\textcolor[rgb]{0.49,0.83,0.13}{P} =P_{0}$};
% Text Node
\draw (65.9,137.52) node [anchor=north west][inner sep=0.75pt]  [color={rgb, 255:red, 126; green, 211; blue, 33 }  ,opacity=1 ]  {$\mathbb{T} E'$};
% Text Node
\draw (163.58,129.84) node [anchor=north west][inner sep=0.75pt]    {$\Gamma $};
% Text Node
\draw (133.88,151.38) node [anchor=north west][inner sep=0.75pt]    {$P_{1}$};
% Text Node
\draw (73.53,75.72) node [anchor=north west][inner sep=0.75pt]    {$\textcolor[rgb]{0.49,0.83,0.13}{P} =P_{0}$};
% Text Node
\draw (187.69,156.47) node [anchor=north west][inner sep=0.75pt]    {$e$};
% Text Node
\draw (130.33,101.73) node [anchor=north west][inner sep=0.75pt]    {$e_{1}$};
% Text Node
\draw (85.42,184.6) node [anchor=north west][inner sep=0.75pt]    {$e_{2}$};
% Text Node
\draw (503.38,234.07) node [anchor=north west][inner sep=0.75pt]  [color={rgb, 255:red, 126; green, 211; blue, 33 }  ,opacity=1 ]  {$\mathbb{T} E'$};
% Text Node
\draw (336.43,84.17) node [anchor=north west][inner sep=0.75pt]    {$e$};
% Text Node
\draw (376.44,37.92) node [anchor=north west][inner sep=0.75pt]    {$e_{1}$};
% Text Node
\draw (375.55,70.69) node [anchor=north west][inner sep=0.75pt]    {$e_{2}$};
% Text Node
\draw (319.98,43.45) node [anchor=north west][inner sep=0.75pt]    {$P_{0}$};
% Text Node
\draw (322.66,129.95) node [anchor=north west][inner sep=0.75pt]    {$P_{1}$};
% Text Node
\draw (330.17,230.68) node [anchor=north west][inner sep=0.75pt]  [color={rgb, 255:red, 126; green, 211; blue, 33 }  ,opacity=1 ]  {$P$};
% Text Node
\draw (509.65,86.82) node [anchor=north west][inner sep=0.75pt]    {$\Gamma $};
% Text Node
\draw (443.25,169.78) node [anchor=north west][inner sep=0.75pt]    {$\varphi '$};

\end{tikzpicture}

     \caption{$\varphi'$ as tropical cover on the left and as harmonic morphism of graphs on the right.}
     \label{figure_phi'astropicalcandgraphc}
 \end{figure}
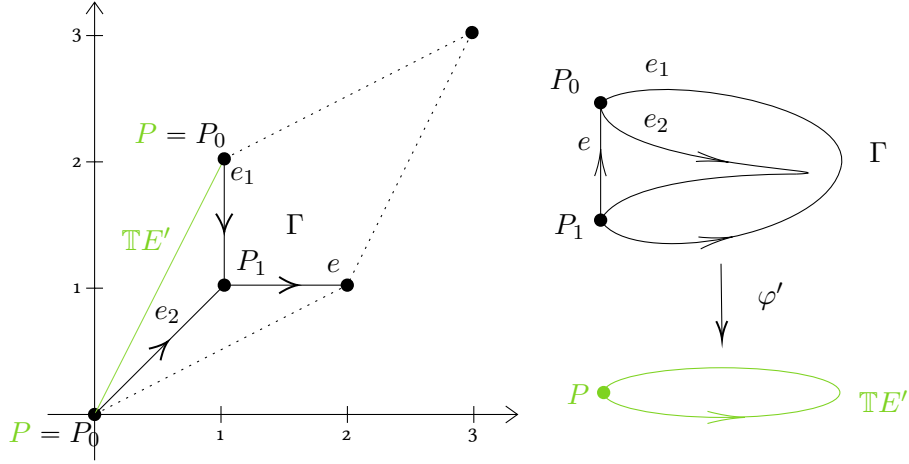
 Recall from Example \ref{example_isogenyarisingfromoptimalcover} that $\ker(\phi)=\{\begin{pmatrix}
     0 \\
     0 
 \end{pmatrix}, \begin{pmatrix}
     \frac{1}{2} \\
     \frac{3}{2}
 \end{pmatrix} \} $. By projecting onto the first coordinate we see $\ker(\phi)\cong \mathbb{R}/\mathbb{Z}[2]$ and $ \ker(\phi) \cong \mathbb{R}/3\mathbb{Z}[2]$ by projecting onto the second. 
\end{example}
\subsection{Algorithm for \ensuremath{\varphi'}}\label{subsection_algorithmforcomplementarycover}
 Identify a cover with the discrete data it is determined by, i.e. if $\Gamma$ is of  
 \begin{itemize}
    \item Type theta with: a triple of \emph{winding numbers} $(n,n_1,n_2)\in \mathbb{N}_0^3$ and \emph{dilation factors} $(d_e(\varphi),d_{e_1}(\varphi),d_{e_2}(\varphi))\in \mathbb{N}_0^3$ such that
    \begin{align}
   & nd_e(\varphi)+(n_1-1)d_{e_1}(\varphi)+(n_2-1)d_{e_2}(\varphi)=d,\\
   & (n-1)d_e(\varphi)+n_1d_{e_1}(\varphi)+n_2d_{e_2}(\varphi)=d.
\end{align}
  
    \item Type dumbbell with: a tuple of \emph{winding numbers} $(n_1,n_2)\in \mathbb{N}_0^2$ and \emph{dilation factors}
    $(d_{e_1}(\varphi),d_{e_2}(\varphi))\in \mathbb{N}_0^2$ such that 
   \begin{align}
    n_1d_{e_1}(\varphi)+n_2d_{e_2}(\varphi)=d.
\end{align}
\end{itemize}
 Given an optimal cover, we can compute the complementary cover:
 \begin{algorithm}\label{algorithm_complementarycover}
     Input: An optimal cover $\varphi: \Gamma \rightarrow \mathbb{T}E$.\\
     Output:  The complementary cover $\varphi': \Gamma \rightarrow \mathbb{T}E'$.
    \begin{enumerate}
         \item Determine $\ker(\varphi_*)$ as described in Lemma \ref{lemma_kernelofpushforwardinTA}, if $\Gamma$ is of type theta, or Lemma \ref{lemma_kernelpushforwardinTA_DB}, if $\Gamma$ is of type dumbbell. Set $\mathbb{T}E':=\ker(\varphi_*)$.
         \item Determine local dilation factors as described in Lemma \ref{lemma_thecomplementarycover} in computation (\ref{diagram_inproofofthecomplentarycover}), if $\Gamma$ is of type theta, or computation (\ref{diagram_inproofofthecomplentarycover_DB}), if $\Gamma$ is of type dumbbell.
         \item Determine winding numbers by using that $\varphi'$ is locally integer  affine linear: Each edge gives rise to an equation involving dilation factors, winding numbers, and the metric data of $\Gamma$ and $\mathbb{T}E'$.
     \end{enumerate}
 \end{algorithm}
\subsection{From Covers to Curves ultimately}\label{subsection_fromcoverstocurvesultimately}
Bringing together the two perspectives, that of $\varphi$ and that of $\varphi'$, finally reveals the symmetry hidden in the proof of Theorem \ref{theorem_Jacsplits1}: Given an optimal cover $\varphi$, we have constructed another optimal cover, $\varphi'$, whose interaction with $\varphi$ (through push-forward and pull-back of divisors) generates two so-called \emph{complementary} exact sequences
\begin{align}
 &0 \rightarrow \mathbb{T}E' \xrightarrow{\varphi'^{*}} \Jac(\Gamma) \xrightarrow{\varphi_*}  \mathbb{T}E \rightarrow 0  \\ 
   &   0 \rightarrow \mathbb{T}E \xrightarrow{\varphi*} \Jac(\Gamma) \xrightarrow{\varphi_*'}  \mathbb{T}E' \rightarrow 0.
\end{align}
The pair of maps, $(\varphi'^{*},\varphi^{*})$ and $(\varphi'_{*},\varphi_{*})$, in turn give rise to isogenies, $\phi$ and $\tilde{\phi}$, by utilizing the coproduct, respectively the product, property of $\mathbb{T}E' \bigoplus \mathbb{T}E$ (recall that finite products coincide with finite coproducts in $\mathbb{T}\mathcal{A}$). We capture the pattern of their interactions in the diagram below:
\begin{diagram}\label{diagram_fromctocultimately}
     %\begin{tikzcd}[row sep=huge]
     \mathbb{T}E' \ar[dd,"m_d"] \ar[dr,"\varphi^{'*}",sloped] \ar[r,"\iota_1", hook] & \mathbb{T}E' \bigoplus \mathbb{T}E \ar[d,dashed,"{\phi}" description] &\mathbb{T}E \ar[dl,"\varphi^{*}",sloped] \ar[l,"\iota_1", hook', swap] \ar[dd,"m_d"]\\
        & \ar[dl,"\varphi'_{*}",sloped] \Jac(\Gamma) \ar[d,dashed,"{\tilde{\phi}}" description] \ar[dr,"\varphi_{*}" ] &  \\
        \mathbb{T}E'  & \ar[l, "p_1"] \mathbb{T}E' \bigoplus \mathbb{T}E \ar[r, "p_2"]  &\mathbb{T}E  
   % \end{tikzcd}
\end{diagram}
 where $\iota_i$ and $p_i$ are the canonical injections, respectively projections. The diagonals are formed by our exact sequences and a small diagram-chase shows that $\tilde{\phi}\circ \phi$ is the componentwise multiplication-by-$d$ map.
 \begin{example}\label{example_interactionpattern}
     The following diagram summarizes the interaction of the two optimal covers $\varphi$ and $\varphi'$ from Example \ref{example_complementarycover} (see Diagram \ref{diagram_fromctocultimately}):
 \begin{diagram}[ampersand replacement = \&]
     %\begin{tikzcd}[row sep=huge]
     \mathbb{R}/\mathbb{Z} \ar[dddd,"\cdot 2",swap] \ar[ddr,"{\begin{psmallmatrix}
     1 \\
     2
 \end{psmallmatrix}}",sloped] \ar[r,"{\begin{psmallmatrix}
     1\\
     0
 \end{psmallmatrix}}", hook] \& \mathbb{R}/\mathbb{Z} \bigoplus \mathbb{R}/3\mathbb{Z} \ar[dd,"{\begin{psmallmatrix}
     1 &1 \\
     2 & 0
 \end{psmallmatrix}}" description] \& \mathbb{R}/3\mathbb{Z} \ar[ddl,"{\begin{psmallmatrix}
     1  \\
      0
 \end{psmallmatrix}}",sloped] \ar[l,"{\begin{psmallmatrix}
     0 \\
     1
 \end{psmallmatrix}}", hook', swap] \ar[dddd,"\cdot 2"]\\
  \& \& \& \\
        \& \ar[ddl,"{\begin{psmallmatrix}
     0 & 1 \\
 \end{psmallmatrix}}",sloped] \mathbb{R}^2/ \begin{pmatrix}
     2 & 1\\
     1 & 2
 \end{pmatrix}\mathbb{Z}^2 \ar[dd,"{\begin{psmallmatrix}
     0 &1 \\
     2 & -1
 \end{psmallmatrix}}"description ] \ar[ddr,"{\begin{psmallmatrix}
     2 & -1
 \end{psmallmatrix}}" ,sloped] \&  \\
 \& \& \& \\
        \mathbb{R}/\mathbb{Z}  \& \ar[l, "{\begin{psmallmatrix}
     1 & 0
 \end{psmallmatrix}}"] \mathbb{R}/\mathbb{Z} \bigoplus \mathbb{R}/3\mathbb{Z} \ar[r, "{\begin{psmallmatrix}
     0 &1 
 \end{psmallmatrix}}", swap]  \& \mathbb{R}/3\mathbb{Z}  
   % \end{tikzcd}
\end{diagram}
 \end{example}

 \subsection{Proof of Theorem \ref{theorem_Jacsplitsiffcoverexists}}\label{subsection_Proof}
 We use methods and results from Subsection \ref{subsection_canonicalcomplement} and \ref{subsection_complementarycover} to complete the proof of Theorem \ref{theorem_Jacsplitsiffcoverexists}. 
 \begin{proof}[Proof of Theorem \ref{theorem_Jacsplitsiffcoverexists}]
 Suppose $\Jac(\Gamma)$ splits, i.e. there exists an isogeny 
 \begin{align}
     \phi: \Jac(\Gamma)\rightarrow \mathbb{T}E\oplus \mathbb{T}E'.
 \end{align}
 We have to show that $\Gamma$ covers an elliptic curve. To this end, define $\varphi: \Gamma \rightarrow \mathbb{T}E$ (in analogy to the classical case, see \cite{thesisDoercksen}) as the composition 
   \begin{align}
       \Gamma  \xrightarrow{\Phi_{P_0}} \Jac(\Gamma) \overset{\phi}{\rightarrow} \mathbb{T}E\oplus \mathbb{T}E' \xrightarrow{p_1}    \mathbb{T}E,
    \end{align} 
    where $p_1$ is the projection onto the first factor. To realize that $\varphi$ is a tropical cover we proceed exactly as described in the proof of Lemma \ref{lemma_thecomplementarycover}. We only mention the points where adjustments have to be made. A lift $\tilde{\varphi}$ of $\varphi$ is given by the upper path in the following diagram: 
   \begin{diagram}
     % \begin{tikzcd}
 \Gamma \arrow[r, "\widetilde{\Phi}_{P_0} "] \arrow[dr, "\Phi_{P_0} "] 
&  \Hom(\Omega^1_\Gamma,\mathbb{R}) \arrow[r, "\Hom(\phi^\#) "] \arrow[d, "\pi_1 "] & \Hom(\Omega^1_{\mathbb{T}E}\oplus \Omega^1_{\mathbb{T}E'},\mathbb{R}) \arrow[d, "\pi_2 "] \arrow[r, "\Hom(p_1^\#) "] & \Hom(\Omega^1_{\mathbb{T}E},\mathbb{R}) \arrow[d, "\pi_3 "]\\
& \Jac(\Gamma) \arrow[r,"\phi"]  &  \mathbb{T}E\oplus \mathbb{T}E' \arrow[r,"p_1"] &  \mathbb{T}E\\
%\end{tikzcd}  
  \end{diagram} 
  where $\pi_i$ for $i=1,2,3$ are the canonical projections. To verify that
  \begin{align}
      \frac{1}{l(\hat{e})}\int_{\hat{e}} \phi^\# \circ p_1^\#(\cdot) \in \Hom(\Omega^1_{\mathbb{T}E},\mathbb{Z})
  \end{align}
  holds, let $\omega$ be a generator of $\Omega^1_{\mathbb{T}E}$ and let
  \begin{align}
     \phi^\# \circ p_1^\#(\omega)= \phi^\#(\omega,0)= \tilde{w}_1 \omega_1 + \tilde{w}_2 \omega_2
  \end{align}
  for some $\tilde{w}_1,\tilde{w}_2\in \mathbb{Z}$ be basis representation. Finally, replace $w_i$ in computations (\ref{diagram_inproofofthecomplentarycover}) and (\ref{diagram_inproofofthecomplentarycover_DB}) by $\tilde{w}_i$ and note that this does not make any difference to the remaining arguments.\\ 
  For the other direction, suppose $\Gamma$ covers an elliptic curve $\mathbb{T}E$. By Corollary \ref{corollary_anycoverfactorsthroughanoptimalone} we know that any cover $\varphi: \Gamma \rightarrow \mathbb{T}E$ factors through an optimal cover, say $\varphi''$:
  \begin{center}
     \begin{tikzpicture}
\matrix(m)[matrix of math nodes,
row sep=3em, column sep=2.8em,
text height=1.5ex, text depth=0.25ex]
{\Gamma& \mathbb{T}E'' & \mathbb{T}E.\\};
\path[->]
(m-1-1) edge node[above] {$\varphi''$} (m-1-2)
(m-1-2) edge node[above] {$\phi''$} (m-1-3)
(m-1-1) edge [bend right] node[below] {$\varphi$} (m-1-3)
;
\end{tikzpicture}
    
\end{center}
 Then, by Theorem \ref{theorem_Jacsplits1} $\Jac(\Gamma)$ splits. 
 \end{proof}
% \newpage
 %\printbibliography
 \bibliographystyle{plain} 
\bibliography{SplitJacobians} 
 \end{document}